\input ./style/arxiv-general.cfg
\documentclass[aop,MSNbibl,seceqn,dvips]{arximspdf}
\makeatletter
   \@ifpackageloaded{graphicx}{}{\usepackage{graphicx}}
\makeatother
\usepackage{mathbh}

\doi{10.1214/14-AOP986}
\volume{44}
\issue{1}
\pubyear{2016}
\firstpage{684}
\lastpage{738}
\docsubty{FLA}

\makeatletter
\newcommand{\rrvert}{\vert}
\newcommand{\rrVert}{\Vert}
\newcommand{\llvert}{\vert}
\newcommand{\llVert}{\Vert}
\renewcommand{\mid}{|}
\newcommand{\mathbbm}{\mathbh}
\renewcommand{\P}{\mathbb{P}}
\renewcommand{\a}{\alpha}
\renewcommand{\b}{\beta}
\renewcommand{\l}{\lambda}
\renewcommand{\O}{\Omega}
\newcommand{\mcG}{\mathcal{G}}
\newcommand{\mcF}{\mathcal{F}}
\newcommand{\mcM}{\mathcal{M}}
\newcommand{\mcS}{\mathcal{S}}
\newcommand{\N}{\mathbb{N}}
\newcommand{\g}{\gamma}
\newcommand{\vp}{\varphi}
\newcommand{\ve}{\varepsilon}
\newcommand{\supp}{\operatorname{supp}}
\newcommand{\td}{\tilde}
\newcommand{\implies}{\Longrightarrow}
\newcommand{\s}{\sigma}
\newcommand{\R}{\mathbb{R}}
\newcommand{\Z}{\mathbb{Z}}
\newcommand{\D}{\Delta}
\newcommand{\E}{\mathbb{E}}
\newcommand{\mcD}{\mathcal{D}}
\newcommand{\eqref}[1]{(\ref{#1})}
\newtheorem{theorem}{Theorem}[section]
\newtheorem{condition}[theorem]{Condition}
\newtheorem{corollary}[theorem]{Corollary}
\newproclaim{definition}[theorem]{Definition}
\newproclaim{example}[theorem]{Example}
\newtheorem{lemma}[theorem]{Lemma}
\newtheorem{proposition}[theorem]{Proposition}
\newproclaim{remark}[theorem]{Remark}
\makeatother

\begin{document}
\begin{frontmatter}

\title{The Jain--Monrad criterion for rough paths and applications to
random Fourier series and non-Markovian H\"{o}rmander theory}
\runtitle{Jain--Monrad for rough paths}

\begin{aug}
\author[A]{\fnms{Peter K.}~\snm{Friz}\corref{}\ead[label=e1]{friz@math.tu-berlin.de}\ead[label=e11]{friz@wias-berlin.de}\thanksref{T1}},
\author[B]{\fnms{Benjamin}~\snm{Gess}\thanksref{T1}\ead[label=e2]{gess@math.tu-berlin.de}},
\author[C]{\fnms{Archil}~\snm{Gulisashvili}\ead[label=e3]{gulisash@ohio.edu}}
\and
\author[D]{\fnms{Sebastian}~\snm{Riedel}\thanksref{T1,T2}\ead[label=e4]{riedel@math.tu-berlin.de}}
\runauthor{Friz, Gess, Gulisashvili and Riedel}
\affiliation{Technische Universit\"at Berlin,
Humboldt-Universit\"at zu Berlin,\\
Ohio University
and
Technische Universit\"at Berlin}
\address[A]{P.~K. Friz\\
Institut f\"ur Mathematik\\
Technische Universit\"at Berlin\\
Strasse des 17. Juni 136\\
10623 Berlin\\
Germany\\
and\\
Weierstra{\ss}-Institut f\"ur Angewandte Analysis\\
\quad und Stochastik\\
Mohrenstrasse 39\\
10117 Berlin\\
Germany\\
\printead{e1}\\
\phantom{E-mail:\ }\printead*{e11}}
\address[B]{B. Gess\\
Department of Mathematics\hspace*{10pt}\\
University of Chicago\\
5734 S. University Avenue\\
Chicago, Illinois 60637\\
USA\\
\printead{e2}}
\address[C]{A. Gulisashvili\\
Department of Mathematics\\
Ohio University\\
Morton Hall 321\\
1 Ohio University\\
Athens, Ohio 45701\\
USA\\
\printead{e3}}
\address[D]{S. Riedel\\
Institut f\"ur Mathematik\\
Technische Universit\"at Berlin\\
Strasse des 17. Juni 136\\
10623 Berlin\\
Germany\\
\printead{e4}}
\end{aug}
\thankstext{T1}{Supported by the European Research Council under the European Union's Seventh Framework
Programme (FP7/2007-2013)/ERC Grant agreement nr. 258237.}
\thankstext{T2}{Supported in part by a scholarship from the Berlin Mathematical School (BMS).}

%
\received{\smonth{11} \syear{2013}}
%
\revised{\smonth{7} \syear{2014}}

%
\begin{abstract}
We discuss stochastic calculus for large classes of Gaussian processes,
based on rough path analysis. Our key condition is a covariance measure
structure combined with a classical criterion due to Jain and Monrad
[\textit{Ann. Probab.} \textbf{11} (1983) 46--57].
This condition is verified in many examples, even in absence of explicit
expressions for the covariance or Volterra kernels. Of special interest are
random Fourier series, with covariance given as Fourier series itself,
and we
formulate conditions directly in terms of the Fourier coefficients. We also
establish convergence and rates of convergence in rough path metrics of
approximations to such random Fourier series. An application to SPDE is
given. Our criterion also leads to an embedding result for
Cameron--Martin paths and complementary Young regularity (CYR) of the
Cameron--Martin space and Gaussian sample paths. CYR is known to imply
Malliavin regularity and also
It\^{o}-like probabilistic estimates for stochastic integrals (resp.,
stochastic differential equations)
despite their (rough) pathwise construction. At last, we
give an application in the context of non-Markovian H\"{o}rmander theory.
\end{abstract}

%
\end{frontmatter}

\setcounter{footnote}{1}

\section*{Introduction}

There\label{intsec1} is a lot of interest, from financial mathematics to nonlinear SPDE
theory, in having a stochastic calculus for nonsemimartingales. In the
past, much emphasis was laid upon stochastic integration (resp., stochastic
differential equations) driven by fractional Brownian motion (fBm), and then
general Volterra processes; cf., for example \cite{N06}, Section~5,
\cite{D05}.
More recently, an effort was made to dispense with the Volterra structure
(cf. \cite{KRT07,KR10}) leading to a key condition of finite planar
(or 2D)
variation of the covariance.
A completely different approach was started by Lyons \cite{L98}; cf. also
\cite{LQ02,LCL07,FV10,FH14}. In essence, it suffices to have a.s.
enough $%
p $-variation regularity of sample paths $X_{\cdot} ( \omega
) $
\textit{and} \textit{existence of stochastic area}(\textit{s}), also subject
to some
variation-type regularity. The problem is then shifted away from developing
a general stochastic integration theory to the (arguably) much simpler task
of constructing the first few iterated (stochastic) integrals; the rest then
follows from deterministic rough path integration theory.

In the case of Gaussian sample paths, a general sufficient condition
for the existence of stochastic areas was introduced in \cite{FV10}.
Namely, it was shown that if the covariance of the underlying process
is sufficiently regular in terms of finite two-dimensional $\rho$-variation, the process can be enhanced with stochastic areas in a
canonical way. The point is that uniform $%
L^{2}$-estimates on the stochastic \mbox{areas} (more precisely, smooth
approximations thereof) are possible, thanks to two-dimensional Young
estimates, as long
as $\rho<\rho^{\ast}=2 $. It is then fairly
straightforward and carried out in detail in \cite{FV10}, Chapter~15 (cf.
also \cite{FH14}) to construct a \mbox{(random)} rough path $\mathbf{X}$ associated
to $X$. This setup has proven rather useful, applications include
non-Markovian H\"ormander theory (\cite{CF10,CHLT12}, more below)
and Hairer's construction \cite{H11,H12} of a spatial rough path
associated to the stochastic heat equation
(in one space dimension) which laid the foundation to prove
well-posedness of certain nonlinear SPDEs. However, finding bounds for
the \mbox{$\rho$-}variation of the covariance of a stochastic process in
concrete examples is not an easy task, and checkable conditions have
been dearly missing in the literature.\footnote{The situation is
easier when $\rho=1$.
In this case, the covariance has finite $1$-variation
if and only if its mixed distributional derivative is a finite signed
measure. In the fBm case this means precisely $H\geq1/2$.} Providing
such conditions is the first main contribution of the present work.

These conditions immediately apply to known examples such as fractional
Brownian motion with Hurst parameter $H$. In this case, it is known
that $\rho=1/ ( 2H ) \vee1$ and the critical $\rho< 2$
corresponds to
$H>1/4$; sharpness of this condition follows from the well-documented
divergence of the L\'{e}vy area for $H^{\ast}=1/4$.


Knowing the precise parameter $\rho$ also has other benefits: it was
shown (cf. \cite{FV10-2}) that finite $\rho$-variation of the
covariance of a Gaussian process implies that the Cameron--Martin space
$\mathcal{H}$ can be continuously embedded in the space of paths with
finite $\rho$-variation; in other words,
\[
\mathcal{H} \hookrightarrow C^{\rho\mbox{-}\mathrm{var}}
\]
holds. In the case $\rho< 3/2$, this embedding assures that the mixed
iterated integral between a Gaussian sample path and a Cameron--Martin
path can be defined via Young's integration theory, and we thus speak
of ``complementary Young regularity'' (CYR) here. CYR has many
consequences: for instance, it allows for a
Malliavin calculus \cite{CFVTAMS,CF10}, \cite{FV10}, Chapters~15, 20,
w.r.t. Gaussian rough paths. In fact, SDE solutions---by
which we mean solutions to rough differential equations driven by
$\mathbf{X}%
( \omega) $ for a.e. $\omega$---will a.s. be Fr\'{e}%
chet-smooth in Cameron--Martin directions as long as CYR holds.
This led to the development of non-Markovian H\"{o}%
rmander theory \cite{CF10,CHLT12}, a significant extension of previous work
\cite{BH07} specific to fBm with $H>1/2$. CYR is important
also for other reasons. It is the condition
under which one has Stroock--Varadhan-type support theorems (see \cite{FV10}, Chapter~19,
and the references therein). It is also the key to good probabilistic estimates
in (Gaussian) rough path theory. To appreciate this, note that the
available pathwise estimates in rough path theory are ill-suited to see the
probabilistic cancellations which are the heart of the It\^{o} theory. It
was only recently understood that Gaussian isoperimetry (in the form of the
Borell--Sudakov--Tsirelson inequality) can bridge this gap (cf.~\cite
{CLL13} and also \cite{FO10}):
in the generic setting of $\rho=1$, if applied to stochastic integrals
(cf. \cite{FR13})
of $\operatorname{Lip}$ 1-forms (as it is typical in rough path
integration theory),
one obtains identical (Gaussian) moment estimates as in the It\^{o}
theory. This deteriorates as $\rho$ increases, but exponential
integrability---and even better, depending on $\rho$---remains
true.\footnote{%
Such integrability properties can be crucial in SPDE theory \cite
{H11,FR13,DFSforth} and in robust filtering theory \cite{CDFO12,DO13}.}
A natural question is whether one can extend CYR to processes which
have finite $\rho$-variation for $\rho\geq3/2$. In the case of
fractional Brownian motion, a direct analysis of its Cameron--Martin
paths (using the Volterra structure of fBm) reveals that in this
special case the stronger embedding
\[
\mathcal{H}^H \hookrightarrow C^{q\mbox{-}\mathrm{var}}\qquad\mbox{for any }
q > \frac{1}{H + (1/2)}
\]
holds (cf. \cite{FV06}) which implies CYR for all $H > 1/4$. Another
contribution of the present work is to show that this stronger
embedding holds in much greater generality and, in particular, even in
absence of a Volterra structure of the process under consideration,
which readily implies CYR for all $\rho< 2$ and thus closes this gap.

The structure of our article is as follows. In Section~\ref{sec1} we
answer in the affirmative the following question: given a
multidimensional Gaussian process with covariance of finite $\rho
$-variation, $\rho< 2$, does CYR hold? The caveat here is that the
$\rho$ is not related anymore to the $\rho$-variation of the covariance
but instead to finite \textit{mixed} $ (1,\rho)$-variation, a
mild strengthening that we prove not to be restrictive at all in
applications. The usefulness of such a result stands and falls with
one's ability to verify this condition in concrete cases. The situation is
aggravated by the examples from random Fourier series (rFs) where the
covariance itself is not known explicitly, but only given as a Fourier
series in its own right.
A general and checkable condition for finite mixed $(1,\rho)$%
-variation is the main result of Section~\ref{sec2}; see Theorem~\ref
{teosuffcritmixed2dvar}.
Loosely speaking, our condition
is a combination of a classical criterion for Gaussian processes to
have $p$%
-variation sample paths due to Jain--Monrad, with a covariance measure
structure condition (the distributional mixed derivative is assumed to
be Radon away
from the diagonal). We then run through a (long) list of examples (see
Examples~\ref{exFBM}--\ref{exspectralmeasure}) which illustrate
the wide
applicability of our criterion. (This way, we also recover from general
principles previously-known results on fBM, such as \cite{FV06}.)
In Section~\ref{sec3} we apply the results of
Section~\ref{sec2} to study rFs in greater depth. In particular, once we have
established finite $\rho$-variation for the covariance of rFs and
therefore the existence of
associated (random) rough paths, we ask for convergence (with rates in rough
path metrics) of natural approximations given in terms of Fourier
multipliers.\footnote{%
Rough path convergence of piecewise linear-, mollifier, Karhunen--Loeve
approximation follows from general Gaussian rough path theory \cite{FV10}
and requires no further discussion.} The best rates are obtained by
considering the rough paths under consideration as $p$-rough paths with
large $p$, which also means that one has to go beyond level $2,3 $
considerations. Thankfully, we can rely here on general results for Gaussian
rough paths established in \cite{FR12a}. The main results in Section~\ref{sec3} are
Theorems~\ref{teorfs}~and~\ref{teodistliftRFSasRP}. In
Sections \ref{sec4}~and~\ref{sec5}, we discuss some concrete rFs (resp., random Fourier
transforms) arising from (fractional) stochastic heat equations in the study
of the stochastic Burgers's \cite{H11} and the KPZ \cite{H12} equation.
Namely, we show how to regard a [fractional, with dissipative term
$- (
-\partial_{xx} ) ^{\alpha}, \alpha\leq1$] heat equation with
space--time white noise, on bounded intervals subject to various boundary
conditions (resp., the entire real line) as an evolution in rough path space.
The key here is spatial covariance of finite $\rho$-variation, where $%
2\alpha=1+1/\rho$. Note $\rho=1$ if and only if $\alpha=1$ and that
$\alpha>3/4$
is handled by our theory.\footnote{%
The covariance structure, including local decorrelation as measured by
mixed variational regularity, of the fractional SHE in the space variables
is similar to fBm with $H=\alpha-1/2$.} This type of spatial rough
path was
first used by Hairer (with $\alpha=1$, and periodic boundary
conditions) to
analyze the stochastic Burgers equation \cite{H11}; a similar construction
with other boundary conditions (incl. those we handle here) was left as open
(technical) problem in \cite{H11}. In a recent preprint, Gubinelli et al.
\cite{GIP12} consider the fractional stochastic Burgers equation, also
with periodic boundary conditions, when $%
\alpha> 5/6$ based on a direct spatial rough path
construction.\footnote{%
In absence of $\rho$-variation estimates, no conclusions toward CYR
and its
numerous consequences are drawn in \cite{GIP12}, nor do the results
allow one to
use the general body of Gaussian rough path approximation theory \cite%
{FV10-2,FV10,FR12a} based on uniform $\rho$-variation estimates. That said,
the overall aim of \cite{GIP12} was quite different.} Finally, in
Section~\ref{sec6} we illustrate (by the example of a driving rFs) how our
results can
be used to check the technical conditions put forward in \cite{CHLT12}
(cf. also \cite{CF10,HP11}), under which differential equations driven
by such
Gaussian signals and along H\"ormander vector fields possess a smooth
density at positive times.

\subsection*{Notation}
\label{secnotation} Let $I = [S,T] \subset\R$ be a closed interval. We
define the simplex by $\D_I:= \{(s,t) \mid s \le t \in I\}$. A
dissection $D$ of an interval $I=[S,T]$ is of the form
\[
D= ( S=t_{0}\leq t_{1}\leq\cdots\leq t_{n}=T
),
\]
and we write $\mathcal{D} ( I ) $ for the family of all such
dissections.

We will now very briefly recall the elements of rough paths theory used in
this paper. For more details we refer to \cite{FV10}. Let $T^N(\R^d)
=\R%
\oplus\R^d \oplus(\R^d \otimes\R^d) \oplus\cdots\oplus(\R
^d)^{\otimes
N}$ be the truncated step-$N$ tensor algebra. For paths in $T^N(\R^d)$
starting at the fixed point $e:= 1 + 0 + \cdots+ 0$, one may define
$\beta$%
-H\"older and $p$-variation metrics, extending the usual metrics for paths
in $\R^d$ starting at zero: the \textit{homogeneous} $\beta$-H\"
older and $p$%
-variation metrics will be denoted by $d$\tsub{$\beta$-H\"ol} and
$d$\tsub{$p$-var}, the \textit{inhomogeneous} ones by $\rho
$\tsub{$\beta$-H\"ol}
and $\rho$\tsub{$p$-var}, respectively. Note that both
\mbox{$\beta$-}H\"older
and $p$-variation metrics induce the same topology on the path spaces.
Corresponding norms are defined by $\llVert \cdot\rrVert$\tsub{$\beta$-H\"ol}${} =
d\mbox{\tsub{$\beta$-H\"ol}}(\cdot,0)$ and $\llVert \cdot\rrVert
\mbox{\tsub{$p$-var}} = d\mbox{\tsub{$p$-var}}(\cdot,0)$ where~$0$ denotes the constant
$e$-valued path.

A geometric $\beta$-H\"older rough path $\mathbf{x}$ is a path in $%
T^{\lfloor1/\beta\rfloor}(\R^d)$ which can be approximated by lifts of
smooth paths in the $d$\tsub{$\beta$-H\"ol} metric; geometric $p$-rough
paths are defined similarly. Given a rough path $\mathbf{x}$, the projection
on the first level is an $\R^d$-valued path and will be denoted by
$\pi_1(%
\mathbf{x})$. It can be seen that rough paths actually take values in the
smaller set $G^N(\R^d) \subset T^N(\R^d)$, where $G^N(\R^d)$ denotes the
free step-$N$ nilpotent Lie group with $d$ generators. The
Carnot--Caratheodory metric turns $(G^N(\R^d),d)$ into a metric space.
Consequently, we denote by
\[
C^{0,\beta\mbox{-H\"ol}}_0\bigl(I,G^{\lfloor1/\beta\rfloor}\bigl(\R^d
\bigr)\bigr)\quad\mbox{and}\quad C^{0,p\mbox{-}\mathrm{var}}_0
\bigl(I,G^{\lfloor p
\rfloor}\bigl(\R^d\bigr)\bigr)
\]
the rough paths spaces where $\beta\in(0,1]$ and $p \in[1,\infty)$. Note
that both spaces are Polish spaces.

\section{Complementary Young regularity under mixed \texorpdfstring{$(1,\rho)$}{(1,rho)}-variation assumption}\label{sec1}

Let $X\dvtx[0,T]\to\R$ be a real-valued, centered, continuous Gaussian process
with covariance
\[
R_{X}(s,t)=\E X_{s}X_{t}.
\]
We will denote the associated Cameron--Martin space by $\mathcal{H}$.
It is
well known that $\mathcal{H}\subset C([0,T],\R)$ and each $h\in
\mathcal{H}$
is of the form $h_{t}=\E ZX_{t}$ with $Z$ being an element of the
$L^{2}$%
-closure of $\operatorname{span}\{X_{t}\mid t\in[0,T]\}$, a Gaussian
random variable. If
$h_{t}=\E ZX_{t}$, $h'_{t}=\E Z'X_{t}$, $\langle
h,h' \rangle_{%
\mathcal{H}}=\E ZZ'$.

For any function $h\dvtx[0,T]\to\R$ we define $h_{s,t}:=h_{t}-h_{s}$ for
all $s,t%
\in[0,T]$. We recall the definition of mixed right $(\gamma,\rho)$-variation
given in \cite{Tow02}: for $\gamma,\rho\geq1$ let
%
\begin{eqnarray}\label{eqndefmixedrhovar}
&& V_{\gamma,\rho}\bigl(R_{X};[s,t]\times
[u,v]\bigr)
\nonumber\\[-8pt]\\[-8pt]\nonumber
&&\qquad :=
\mathop{\sup_{(t_{i})\in
\mathcal{D}([s,t])}}_{(t_{j}')\in\mathcal{D} (
[u,v ] )
} \biggl(\sum
_{t'_{j}} \biggl(\sum_{t_{i}}\biggl|
R_{X}\pmatrix{ t_{i},t_{i+1}
\cr
t_{j}',t_{j+1}'}\biggr|
^{\gamma} \biggr)^{\rho/\gamma} \biggr)^{1/\rho},
\end{eqnarray}
where $\mcD([s,t])$ denotes the set of all dissections of $[s,t]$ and
\[
R_{X}\pmatrix{ t_{i},t_{i+1}
\cr
t_{j}',t_{j+1}'}=\E
X_{t_{i},t_{i+1}}X_{t_{j}',t_{j+1}'}.
\]

The notion of the 2D $\rho$-variation is recovered as $V_{\rho
}=V_{\rho,\rho}
$. Recall that \mbox{$V_{\rho}$-}regularity plays a key role in Gaussian
rough path
theory \cite{FV10,FV10-2,FH14} and in particular yields a stochastic
integration theory for large classes of multidimensional Gaussian processes.
Clearly, $V_{\gamma\vee\rho}(R;A)\leq V_{\gamma,\rho}(R;A)\leq
V_{\gamma\wedge\rho}(R;A)$ for all rectangles $A\subseteq
[0,T]^{2}$. As the
main result of this section, we present the following embedding theorem
for the
Cameron--Martin space.

%
\begin{theorem}\label{teorefinedCMembedding} Assume that the covariance $R_{X}$ has
finite mixed $ (1,\rho)$-variation in 2D sense. Then there
is a
continuous embedding
\[
\mathcal{H}\hookrightarrow C^{q\mbox{-}\mathrm{var}}\qquad\mbox{with }q=
\frac{1}{1/(2\rho)+1/2}<2.
\]
More precisely,
\begin{eqnarray*}
\llVert h\rrVert_{q\mbox{-}\mathrm{var};[s,t]}\leq\llVert h\rrVert
_{\mathcal{H}}
\sqrt{V_{1,
\rho}\bigl(R_{X};[s,t]^{2}\bigr)}\qquad
\forall[s,t]\subseteq[0,T].
\end{eqnarray*}
\end{theorem}

The following is then immediate.

%
\begin{corollary}\label{cor2}
Assume $\rho\in[1,2)$. Then complementary Young regularity holds,
that is, we can choose $p>2\rho$ small enough such that $X$ has a.s.
$p$-variation sample paths, $h \in\mathcal{H}$ has finite
$q$-variation with $1/p+1/q>1$.
\end{corollary}

We shall in see in Section~\ref{sec2} (as one of many examples) that the
assumption of mixed $(1,\rho)$-variation is met in the case of fBm in
the rough regime $H\le1/2$ with $\rho= 1/(2H)$. (E.g.,
Example~\ref{StInII} applies with $k=0$ and in fact gives a neat
criterion for processes with stationary increments.) It then follows
that fractional Cameron--Martin paths enjoy finite $%
q=\frac{1}{H+1/2}$-variational regularity, which is
consistent (and in fact a mild sharpening) of $q>\frac{1}{H+1/2}$, previously obtained in~\cite{FV06} with methods
specific to fBm.
Let us also note that, for the sole purpose of Theorem~\ref{teorefinedCMembedding}, it would have been enough to consider identical
dissections $(t_{i})\equiv(t_{j}')$ in the definition of mixed
variation $V_{\gamma,\rho}$ in \eqref{eqndefmixedrhovar}. The
criteria in Theorem~\ref{teomain} below would then allow for a mildly simplified proof.
On the
other hand, this criteria derived in Theorem~\ref{teomain} below are also
sufficient (and interesting) for finite \mbox{$\rho$-}variation $%
V_{\rho}=V_{\rho,\rho}$ which is the key condition for the
construction of Gaussian rough paths
needed later on,
hence the additional generality of different vertical and horizontal
dissections.

%
\begin{remark}
Let $X\colon[0,T]\to\R^{d}$ be a multidimensional centered Gaussian process.
Then every path $h$ in the associated Cameron--Martin space $\mathcal
{H}$ is
of the form $h_{t}=\E ZX_{t}$ with $Z$ being an element of the $L^{2}$%
-closure of $\operatorname{span}\{X_{t}^{i}\mid t\in I, i=1,\ldots,d\}
$ and $\llVert h\rrVert _{%
\mathcal{H}}=\llVert Z\rrVert _{L^{2}}$. The $q$-variation of $h$ is
finite if and only
if the \mbox{$q$-}variation of every $h_{\cdot}^{i}=\E ZX_{\cdot}^{i}$ is finite,
and we obtain the bound
\[
\llVert h\rrVert_{q\mbox{-}\mathrm{var};[s,t]}\leq C\llVert h\rrVert
_{\mathcal{H}}\max
_{i=1,\ldots,d}\sqrt{%
V_{1,\rho}\bigl(R_{X^{i}};[s,t]^{2}
\bigr)},
\]
where $C$ is a constant depending only on the dimension $d$.
\end{remark}

We now give the proof of Theorem~\ref{teorefinedCMembedding}. In
fact, having identified
the importance of mixed variation, the proof is pleasantly short.

\begin{pf*}{Proof of Theorem \ref{teorefinedCMembedding}}
Let $h=\E ZX_{\cdot}\in\mathcal{H}$. Fix a dissection $D=
(t_{j} )%
\subset[s,t]$, write $h_{j}\equiv h_{t_{j},t_{j+1}},X_{j}=X_{t_{j},t_{j+1}}$
and also $\llVert h\rrVert _{q}^{q}:=\sum_{j}\llvert
h_{j}\rrvert ^{q}$. Let $q'$ and $\rho'$ be the
conjugate exponents of $q$ and $\rho$. An easy calculation shows that $
\rho'=q'/2$. By duality,
\[
\llVert h\rrVert_{q}=\sup_{\beta\dvtx\llVert \beta\rrVert
_{q'}\leq1}\sum
\beta_{j}h_{j}=\sup_{\beta\dvtx\llVert
\beta\rrVert _{q'}\leq1}\E\biggl(Z
\sum_{j}\beta_{j}X_{j}
\biggr),
\]
and so by Cauchy--Schwarz
\[
\llVert h\rrVert_{q}^{2}\leq\llVert h\rrVert
_{\mathcal{H}%
}^{2}\sup_{\beta\dvtx\llVert \beta\rrVert
_{q'}\leq1}\sum
_{j,k}\beta_{j}\beta_{k}{
\E}X_{j}X_{k}.
\]
Set $R_{j,k}=\E X_{j}X_{k}$. Then, using the symmetry of $R$ and H\"older's
inequality,
\begin{eqnarray*}
\sum_{k,j}\beta_{j}
\beta_{k}R_{k,j} & \le&\frac{1}{2}\sum
_{j,k}%
\beta_{j}^{2}\llvert
R_{j,k}\rrvert+\frac{1}{2}\sum_{j,k}
\beta_{k}^{2}\llvert R_{j,k}\rrvert
\\
& =&\sum_{j}\beta_{j}^{2}
\sum_{k}\llvert R_{k,j}\rrvert
\\
& \leq&\llVert\beta\rrVert_{2\rho^{\prime}}^{2} \biggl(\sum
_{j} \biggl(\sum_{k}\llvert
R_{i,k}\rrvert\biggr)^{\rho} \biggr)^{1/\rho}
\\
& \leq& V_{1,\rho}\bigl(R;[s,t]^{2}\bigr)
\end{eqnarray*}
when $\llVert \beta\rrVert _{2\rho'}=\llVert \beta\rrVert _{q'}\leq1$
which shows
the claim.
\end{pf*}

\section{Jain--Monrad revisited}\label{sec2}

\subsection{Preliminaries and motivation from fBm}\label{sec21}

Let $I\subset\R$ be a compact interval and $R\colon I \times I \to\R
$ be a
symmetric, continuous function. We set $T = \llvert I\rrvert $,
%
\begin{equation}
D_{h}:= \bigl\{ (s,t)\in I^{2}\dvtx \llvert s-t\rrvert\leq
h \bigr\} \label{defDh}
\end{equation}
and let $D:=D_{0}$ be the diagonal of $I^{2}$. In this section we will
give conditions under which $R$ has finite
$\rho$%
-variation on $I^2 = I \times I$. For a rectangle $[s,t] \times
[u,v]\subseteq I^2$, we define the rectangular increment by
\[
R\pmatrix{ s,t
\cr
u,v} =R(s,u)-R(s,v)-R(t,u)+R(t,v),
\]
and we set
%
\begin{equation}
\sigma^{2}(s,t):=R\pmatrix{ s,t
\cr
s,t} =R(s,s)+R(t,t)-2R(s,t),
\label{eqnsigma}
\end{equation}
where symmetry of $R$ was used in the last step. Note that
\[
\partial_{s,t}\sigma^{2}=-2\partial_{s,t}R
\]
whenever these mixed derivatives make sense. In many applications $R$
is the
covariance function of a zero mean stochastic process $X$, that is,
$R(s,t)=E X_{s}X_{t}$,
and in this case $\sigma^{2}(s,t)=\operatorname{Var}(X_{t}-X_{s})\geq
0$ is the
variance of increments. However, it will be important to conduct the present
discussion in a generality that goes beyond covariance functions.

Given a dissection $ ( t_{i} ) $ of $I = [ 0,T ]
$, the
square $ [ 0,T ] ^{2}$ can be decomposed into little squares
$\bigcup_{j} [ t_{i},t_{i+1} ] ^{2}$ and off-diagonal rectangles,
say $%
\{ Q_{j} \} $. Then
\[
\sum_{i}\sigma^{2} (
t_{i},t_{i+1} ) +\sum_{j}R (
Q_{j} ) =R\pmatrix{ 0,T
\cr
0,T} =\sigma^{2} ( 0,T ) <
\infty,
\]
and the right-hand side is independent of the dissection. Depending on the
behavior of $\sigma^{2} ( s,t ) $, we can or cannot\vspace*{1pt} ignore the
on-diagonal contributions in the limit $\operatorname{mesh} (
t_{i} )
\rightarrow0 $. For instance, if $\sigma^{2} ( s,t )
=\llvert
t-s\rrvert ^{2H}$with $H>1/2$, then
\[
\lim_{\operatorname{mesh}(t_i)\rightarrow0}\sum_{i}
\sigma^{2} ( t_{i},t_{i+1} ) =0
\]
and with $R ( Q_{j} ) \approx\partial_{s,t}R.\Delta_{j}$
for small $Q_{j}$, or by direct calculus, we find%
%
\begin{eqnarray}\label{sigmaSquareAsDoubleInt}
\sigma^{2} ( 0,T ) &=& T^{2H}=-\frac{1}{2}\int _{0}^{T}\!\!\int_{0}^{T}%
\partial_{s,t}\llvert t-s\rrvert^{2H} \,ds \,dt
\nonumber\\[-8pt]\\[-8pt]\nonumber
&=& H ( 2H-1 ) \int
_{0}^{T}\!\!\int_{0}^{T}
\llvert t-s\rrvert^{2H-2} \,ds \,dt,
\end{eqnarray}
noting that $\llvert t-s\rrvert ^{2H-2}=\llvert t-s\rrvert
^{-1+2 ( H-1/2 ) }$ is integrable at the diagonal (and then
everywhere on $ [ 0,T ] ^{2}$) if and only if $H>1/2$. When
$H=1/2$ this
computation fails. Indeed, the prefactor $2H-1=0$ combined with the
diverging integral effectively leaves us with $0\cdot\infty$. The reason
of course is that $R ( Q_{j} ) =0$ in this case (Brownian
increments are uncorrelated), and everything hinges on the (nonvanishing)
on-diagonal contribution
\[
\sum_{i}\sigma^{2} (
t_{i},t_{i+1} ) =\sum_{i} (
t_{i+1}-t_{i} ) =T.
\]
As a Schwartz distribution $\partial_{s,t}R=\partial_{s,t}\min(
s,t ) =\delta_{ \{ s=t \} }$ is a ``Dirac'' on the
diagonal and
indeed with this interpretation as a measure,
\[
\sigma^{2} ( 0,T ) =R\pmatrix{ 0,T
\cr
0,T} =\int
_{0}^{T}\!\!\int_{0}^{T}
\delta_{ \{ s=t \} } \,ds \,dt = T.
\]
When $H<1/2$, $\sigma^{2} ( s,t ) =\llvert t-s\rrvert ^{2H}$, the
on-diagonal contributions are not only nonvanishing but
divergent [as
the mesh of $ ( t_{i} ) $ goes to zero]. That is,
\[
\sigma^{2} ( 0,T ) =T^{2H}=\underbrace{\sum
_{i}\sigma^{2} ( t_{i},t_{i+1}
) }_{\rightarrow\infty}%
+\sum_{j}R (
Q_{j} )
\]
and so, necessarily, $\sum_{j}R ( Q_{j} ) \rightarrow
-\infty$.
Translated to the calculus setting, this causes (\ref
{sigmaSquareAsDoubleInt}%
) to fail. Indeed, ignoring the infinite contribution from the diagonal
leaves us with
\[
T^{2H}\neq H ( 2H-1 ) \underbrace{%
\int_{0}^{T}\!\!\int_{0}^{T}\llvert t-s\rrvert^{2H-2} \,ds
\,dt}_{=+\infty}=-\infty\qquad\mbox{for }H<1/2.
\]
Let us remark that, with our standing assumption $R\in C ( [
0,T%
] ^{2} ) $ the (distributional) mixed derivative $\partial
_{s,t}R $ always exists, that is,
\[
\langle\partial_{s,t}R,\varphi\rangle:=\int_{0}^{T}\!\!\int_{0}^{T}R ( s,t )\, \partial_{s,t}
\varphi( s,t ) \,ds \,dt\qquad\forall\varphi\in C_{c}^{\infty
} \bigl( (
0,T ) ^{2} \bigr).
\]
One can ask if, or when, $\partial_{s,t}R$ is given by a signed and finite
(i.e., of finite total variation) Borel measure $\mu$ on $ [
0,T ]
^{2}$, say
\[
\langle\partial_{s,t}R,\varphi\rangle=\int_{ [
0,T ]
^{2}}
\varphi \,d\mu,
\]
with associated Hahn--Jordan decomposition $\mu=\mu_+-\mu_-$. When $H>1/2$,
the answer is affirmative with $\mu=\mu_{+}=$ $H ( 2H-1 )
\llvert t-s\rrvert ^{2H-2}\,ds \,dt$. For $H=1/2$, the answer is also
affirmative with $\mu=\mu_{+}=\delta_{ \{ s=t \} }$. For $H<1/2$,
the answer is negative.

However, for all values of $H\in(0,1)$ it is possible to define a
(signed) $\sigma$-finite measure by
\[
\mu( A ):=\int_{A}H ( 2H-1 ) \llvert t-s\rrvert
^{2H-2} \,ds \,dt
\]
which we shall regard as a signed Radon measure on $ ( 0,T )
^{2}\setminus D$. Note
%
\begin{eqnarray}
\mu\equiv\mu_{+},\qquad
\mu\equiv0,\qquad\mu\equiv-\mu_{-}\nonumber
\\
\eqntext{\mbox{for }H>1/2, H=0, H<1/2,\mbox{ respectively}.}
\end{eqnarray}
In general, as seen when $H<1/2$, $\mu$ does not need to be a finite
measure on $ ( 0,T ) ^{2}\setminus D$. On the other hand, its
restriction to any compact in $ ( 0,T ) ^{2}\setminus D$ is finite
so that $\mu$ defines a signed Radon measure on $ ( 0,T )
^{2}\setminus D$. Hence, for all values of $H\in( 0,1 )
$ the
(distributional) mixed derivative $\partial_{s,t}R$ on $ (
0,T )
^{2}\setminus D$ is given by the Radon measure $\mu$. (This was certainly
observed previously, e.g., in~\cite{KR10}.)

Care is necessary, for important information has been lost by the
restriction to $ ( 0,T ) ^{2}\setminus D$. For instance, nothing
was left of Brownian motion ($\mu=0$). It follows that when $H\leq1/2$,
and in particular in the case $H<1/2$ where $\llvert \mu\rrvert =\mu
_{-}$ has infinite mass on $ ( 0,T ) ^{2}\setminus D$, the
on-diagonal information must be captured differently. We shall achieve this
by a somewhat classical condition due to Jain--Monrad \cite{JM83,DN98}
which imposes ``on-diagonal'' $\rho$%
-variation of $\s^{2}$ by
\[
v_{\rho}\bigl(\s^{2};[s,t]\bigr):=\sup_{D=(t_{i})\in\mcD([s,t])}
\biggl( \sum_{i}\bigl\llvert\s^{2}(t_{i},t_{i+1})
\bigr\rrvert^{\rho} \biggr) ^{1/\rho}<\infty.
\]
Clearly $\rho=1/2H\geq1$ in the fBm example with $H\leq1/2$, but the
concept is much more general.

\subsection{Main result of the section}\label{sec22}

Throughout we work on some closed interval $I \subset\R$ with length
$T =
\llvert I\rrvert $.

%
\begin{condition}[(Jain--Monrad)]
Let $\rho\geq1$ and $\omega\colon\Delta_I \to\R_+$ be a super additive
function [i.e., $w(s,r)+w(r,t)\le w(s,t)$ for all $s\le r\le t$]. We say
that $(JM)_{\rho,\omega}$ holds if
\begin{eqnarray*}
&&\bigl\llvert\sigma^2(s,t)\bigr\rrvert\leq\omega(s,t)^{1/\rho}
\end{eqnarray*}
holds for all $s<t$.
\end{condition}

If $v_\rho(\s^2;I)<\infty$, we can always set $\omega(s,t) = v_\rho
(\s%
^2;[s,t])^\rho$. Conversely, if $(JM)_{\rho,\omega}$ holds, we have
$v_\rho(%
\s^2;[s,t]) \leq\omega(s,t)^{1/\rho}$ for all
$[s,t]\subseteq I$.

Recall\vspace*{1pt} the definition of mixed right $(\gamma,\rho)$-variation given
in \eqref{eqndefmixedrhovar}, noting in particular the triangle
inequality: for all rectangles $A\subseteq
I^{2}$,\vspace*{1pt}
%
\begin{equation}
V_{\gamma,\rho}(R_{1}+R_{2};A)\leq V_{\gamma,\rho
}(R_{1};A)+V_{\gamma,\rho}(R_{2};A).
\label{eqnsubadd}
\end{equation}

Recall that a signed Radon measure $\mu$ is a locally finite signed Borel
measure with decomposition $\mu=\mu_{+}-\mu_{-}$ where $\mu_{\pm}$ are
locally finite, nonnegative Borel measures, one of which has finite mass.
For a finite measure $\mu$ on $(0,T)^2\setminus D$ we will consider its
extension to $[0,T]^2$ by $\mu(A):=\mu(A\cap(0,T)^2\setminus D)$ without
further notice.
We now give the main theorem of this section. For simplicity, we only
formulate it for the case $I = [0,T]$.

%
\begin{theorem}
\label{teosuffcritmixed2dvar}\label{teomain} Let $R\colon
[0,T]^{2}\rightarrow\R$ be a symmetric, continuous function and $\s$
as in %
\eqref{eqnsigma}. Assume that the (Schwartz) distributional mixed
derivative $\mu:=\frac{\partial^{2}R}{\partial_{t}\,\partial
_{s}}=-\frac{1%
}{2}\frac{\partial^{2}\s^{2}}{\partial_{t}\,\partial_{s}}$ is a Radon
measure on $ ( 0,T ) ^{2}\setminus D$ with decomposition
\mbox{$\mu=\mu_+ - \mu_-$}.
\begin{longlist}
\item[\textit{Part} A.] Assume that:
\begin{longlist}[(A.ii)]
\item[(A.i)] $\mu_{-}$ has finite mass and a continuous distribution
function.

\item[(A.ii)] There exists an $h>0$ such that $\sigma^{2}(s,t)\geq0$
whenever $\llvert t-s\rrvert \leq h$.\newcounter{fnnumber}\footnote{%
Automatically true if $R$ is a covariance function.} \setcounter
{fnnumber}{%
\thefootnote}
\end{longlist}

Then
%
\begin{eqnarray}
V_{1}\bigl(R;[s,t]\times [ u,v]\bigr)\leq R\pmatrix{ s,t
\cr
u,v} +2
\mu_{-}\bigl([s,t]\times [ u,v]\bigr)\nonumber
\\
\eqntext{\forall [ s,t]\times [
u,v]\subseteq[0,T]^{2}.}
\end{eqnarray}

\item[\textit{Part} B.] Assume that:
\begin{longlist}[(B.iii)]
\item[(B.i)] $\mu_{+}$ has finite mass and a continuous distribution
function.

\item[(B.ii)] There exists an $h>0$ such that\footnote{%
With the exception of bi-fBm, Example~\ref{examplebifbm}, we typically
check (B.ii) by simply showing that $\tau\mapsto\sigma^{2}(\tau,t+\tau)$,
respectively, $\sigma^{2}(t-\tau,t)$ are\vspace*{1pt} nondecreasing for all $t$ and
$\tau<h$.
In particular, in stationary situations where $\sigma^{2} (
s,t )
=F ( t-s ) $ this amounts for $F$ to be nondecreasing on
$ [
0,h ] $; conversely it is not hard to see (\ref{Bii}) implies $F$
nondecreasing on $[0,h/2]$.}
%
\begin{eqnarray}\label{Bii}
2R\pmatrix{ s,t
\cr
u,v} =\sigma^{2}(s,v)-\sigma^{2}(s,u)+
\sigma^{2}(u,t)-\sigma^{2}(v,t)\geq0\nonumber
\\
\eqntext{\forall [ u,v]
\subseteq{}[ s,t]\subseteq I\mbox{ s.t. }\llvert t-s\rrvert\leq h.}
\end{eqnarray}

\item[(B.iii)] $(JM)_{\rho,\omega}$ holds.
\end{longlist}

Then for all $[s,t]^{2}\subset D_{h}$, as defined in (\ref{defDh}),
we have
%
\begin{equation}
\label{eqnon-diag-est} V_{1,\rho}\bigl(R;[s,t]^{2}\bigr)\leq C \bigl(
\omega^{1/\rho}(s,t)+\mu_{+}\bigl([s,t]^{2}\bigr)
\bigr),
\end{equation}
for some constant $C=C(\rho)$.

If, in addition, $R\colon{}[0,T]^{2}\rightarrow\R$ satisfies a
Cauchy--Schwarz inequality\footnote{%
That is, $\llvert R{
s,t \choose
u,v} \rrvert \leq\llvert R{
s,t \choose
s,t} \rrvert ^{1/2}\llvert R{
u,v \choose
u,v} \rrvert ^{1/2}$, for all $[s,t]\times [
u,v]\subseteq I^{2}$, which is automatically true if $R$ is a covariance
function.} then, more generally, there is a constant $C=C(\rho,h,T)$ such
that
\begin{eqnarray}\label{eqnoff-diag-est-1rho}
&& V_{1,\rho}\bigl(R;[s,t]\times [ u,v]\bigr)
\nonumber\\[-8pt]\\[-8pt]\nonumber
&&\qquad \leq C \bigl(
\omega^{1/(2\rho)}(s,t)\omega^{1/(2\rho)}(u,v)+\mu_{+}\bigl([s,t]\times [ u,v]\bigr) \bigr), \nonumber
\end{eqnarray}
for all rectangles $[s,t]\times [ u,v]\subset{}[0,T]^{2}$.
\end{longlist}
\end{theorem}

The interest in Theorem~\ref{teomain} is two-fold. First, it has
far-reaching conclusions: mixed $(1,\rho)$-variation controls $\rho$%
-variation which, if applied (componentwise) to the covariance of a Gaussian
process (multidimensional, with independent components), is the key quantity
for the existence of associated rough paths; here one needs $\rho<2$ (which
corresponds to $H>1/4$; cf. Example~\ref{fbmrough} below).

Let us state the consequence in terms of rough paths construction
specifically as a corollary.

\begin{corollary}
Assume $ ( X_{t}\dvtx0\leq t\leq T ) $ is a $d$%
-dimensional, centered Gaussian process with independent components.
For each component $X^i$, assume that either the assumptions of part~\textup{A}
of Theorem~\ref{teomain}
are satisfied, in which case we set $\rho_i =1 $, or those of part~\textup{B}
for some $\rho_i<2$.
Set $\rho:= \max_{i=1,\dots,d} \rho_i <2$. Then, for any $p> 2 \rho
$, it follows that
$X$ admits a ``canonical'' lift ${\mathbf X}={\mathbf X}(\omega)$ to a
random geometric $p$-rough path.%
\footnote{By ``canonical'' we mean that ${\mathbf X}$ is the limit, in
probability and $p$-variation rough path metric, of standard
approximations procedures including piecewise linear, mollifications and
of Karhunen--Loeve type. We also note that the estimates of Theorem~\ref{teosuffcritmixed2dvar} allow us to show, under natural assumptions
on the quantities appearing on the right-hand side, that the
covariances of $X$ have \mbox{finite}  ``H\"older-controlled'' $\rho
$-variation, thereby
allowing us to conclude that ${\mathbf X}$ is a random geometric
\mbox{$\alpha$-}H\"older rough path, for $\alpha< \frac{1}{2\rho}$. See
\cite{FV10,FH14} for more details.}
\end{corollary}

Moreover, mixed $%
(1,\rho)$-variation was seen in Section~\ref{sec1} to imply \textit{complementary
Young regularity},
an extremely important property leading to good
probabilistic estimates of rough integrals, as explained in the
\hyperref[intsec1]{Introduction}. It is also required for Stroock--Varadhan-type support
theorems and
is one of the key conditions for the applicability of
Malliavin calculus and then non-Markovian H\"ormander theory; cf. \cite%
{CF10,CHLT12}.

Secondly, the theorem is practical because its conditions are easy to check
and widely applicable. To illustrate this we now run through a list of
examples. Roughly speaking, part~\textup{A} handles situations similar or nicer than
Brownian motion, whereas part~\textup{B} handles situations similar or worse than
Brownian motion. The finite measure $m=\mu_{-}$ (resp., $\mu_{+}$) in
part~\textup{A}
(resp., \textup{B}) should be considered as (harmless) perturbation which adds some
extra flexibility. Typically $m$ is given by a density, that is, by the
(integrable) negative (resp., positive) part of some locally integrable
function. Continuity of the distribution function is then trivial. In
fact, $%
m=0$ in many interesting examples.

\subsection{Examples}\label{sec23}

\subsubsection{Examples handled by part~\textup{A}}\label{sec231}

%
\begin{example}[(Fractional Brownian motion $H\geq1/2$)]
\label{exFBM} Consider a (standard) fractional Brownian motion $B^{H}$,
with $\sigma^{2} ( s,t ) =\llvert t-s\rrvert
^{2H}$ in the
regime $H>1/2$.
We have, as a measure on $ [ 0,T ] ^{2}\setminus D$,
\begin{eqnarray*}
 \mu&=&\mu_{+}=H ( 2H-1 ) \llvert t-s\rrvert^{2H-2} \,ds \,dt
\geq0\qquad\mbox{if }H>1/2,
\\
 \mu&=&0\qquad\mbox{if }H=1/2,
\end{eqnarray*}
which\vspace*{1pt} clearly yields a Radon measure on $ [ 0,T ]
^{2}\setminus D$
(and even a finite Borel measure on $ [ 0,T ] ^{2}$). Note
that $%
\mu_{-}\equiv0$ in the\vspace*{1pt} decomposition $\mu=\mu_{+}-\mu_{-}$; hence~(A.i)
holds trivially. Also, since $R ( s,t ) =\frac{1}{2} (
s^{2H}+t^{2H}-\llvert t-s\rrvert ^{2H} ) $ is a genuine
covariance function, (A.ii) comes for free. It follows that $R$ has
finite ``H%
\"{o}lder controlled'' $1$-variation, in the sense that
\[
V_{1} \bigl( R; [ s,t ] ^{2} \bigr) \leq R\pmatrix{ s,t
\cr
s,t} =\llvert t-s\rrvert^{2H}=O \bigl( \llvert t-s\rrvert\bigr).
\]
\end{example}

%
\begin{example}[(Brownian bridge)]
Given a standard Brownian motion $B$, the \textit{Brownian bridge}
over $%
[ 0,T ] $ can be defined as
\[
X_{t}=B_{t}-\frac{t}{T}B_{T}\quad\implies\quad R (
s,t ) =\min( s,t ) -st/T.
\]
It follows that $\mu=\partial_{s,t}R$, as a measure on $ [
0,T ]
^{2}\setminus D$, decomposes into $\mu_{+}=0$ and $\mu_{-}$ with
(constant) density $1/T$. Part~\textup{A} applies and immediately gives ``H\"{o}lder
controlled'' $1$-variation, that is, $V_{1} ( R; [ s,t
] ^{2} )
=O ( \llvert t-s\rrvert ) $.
\end{example}

%
\begin{example}[(Stationary increments I, Brownian and better regularity)]
\label{StInI}Consider a process with stationary increments in the
sense that
the variance of its increments is given by
\[
\sigma^{2} ( s,t ) =F \bigl( \llvert t-s\rrvert\bigr) \geq0,
\]
for some $F\in C^{2} ( [0,T] ) $. A concrete (Gaussian)
example is
the \textit{stationary Ornstein--Uhlenbeck process} with $F (
x )
=1-e^{-x}$. In any case, we may expand
\[
F ( h ) =F' ( 0 ) h+F^{\prime\prime
} ( 0 ) h^{2}/2+o
\bigl( h^{2} \bigr).
\]
We compute
\[
\partial_{s,t}\sigma^{2} ( s,t ) =-F^{\prime\prime} \bigl(
\llvert t-s\rrvert\bigr) +F' ( 0 ) 2\delta( t-s )
\]
so that
\[
\frac{\partial^{2}R}{\partial s\,\partial t}=-\frac{1}{2}\frac
{\partial^{2}\s%
^{2}}{\partial s\,\partial t}=\frac{1}{2}F^{\prime\prime}
\bigl( \llvert t-s\rrvert\bigr)\qquad\mbox{on } ( 0,T ) ^{2}
\setminus D.
\]
It then follows that (A.i) holds with
\[
\mu( A ) =\frac{1}{2}\int_{A}F^{\prime\prime}
\bigl( \llvert t-s\rrvert\bigr) \,ds \,dt,
\]
and we immediately obtain finite (H\"{o}lder controlled) $1$-variation,
\[
V_{1} \bigl( R; [ s,t ] ^{2} \bigr) \leq
\sigma^{2} ( s,t ) +\bigl\llvert F^{\prime\prime}\bigr\rrvert
_{\infty}\llvert t-s\rrvert^{2}=O \bigl( \llvert t-s\rrvert
\bigr).
\]
For a concrete $F$, of course, one can compute $\mu_{-}$ and obtain sharper
conclusions. This may also be possible if we are in a ``better than Brownian''
setting, namely $F' ( 0 ) =0$, in which case
$\sigma
^{2} ( s,t ) =O(\llvert t-s\rrvert ^{2})$. Note that in this case
$F^{\prime
\prime} ( 0 ) >0$, unless $F$ is trivial.\footnote{%
Indeed, if $F' ( 0 ) =F^{\prime\prime} (
0 ) =0$, then $\llVert X_{t}-X_{s}\rrVert _{L^{2}}=o ( t-s
) $
which is enough to conclude that $X_{t}$ is a constant in $L^{2}$, but
then $%
\sigma^{2} ( s,t ) =\llVert X_{t}-X_{s}\rrVert
_{L^{2}}^{2}=0$.} It follows that, in a neighborhood of the diagonal,
$\mu
>0$, and so $\mu_{-}\equiv0$. We then have
\[
V_{1} \bigl( R; [ s,t ] ^{2} \bigr) \leq
\sigma^{2} ( s,t ) =O\bigl(\llvert t-s\rrvert^{2}\bigr),
\]
for $\llvert t-s\rrvert \leq\sup\{h>0\dvtx F^{\prime\prime
}(h)>0\}$.
\end{example}

%
\begin{example}[(Volterra processes I; Brownian and better regularity)]
Assume $X_{t}=\int_{0}^{t}K ( t,r ) \,dB_{r}$ where $K (
t,\cdot) $ is assumed to be square-integrable. For $s<t$, we have
\begin{eqnarray*}
X_{s,t} &=&\int_{0}^{t} \bigl( K (
t,r ) -K ( s,r ) 1_{\{r\leq
s\}} \bigr) \,dB_{r},
\\
\sigma^{2} ( s,t ) &=& \E X_{s,t}^{2}=\int
_{0}^{s} \bigl( K ( t,r ) -K ( s,r ) \bigr)
^{2}\,dr+\int_{s}^{t}K ( t,r )
^{2}\,dr.
\end{eqnarray*}
We assume a regular situation, by which we shall mean here that $K$ is
continuous on the simplex $ \{ 0\leq s\leq t\leq T \} $, and
assuming suitable differentiability properties of $K$, one computes
\[
\partial_{s,t}R=K ( s,s )\, \partial_{t}K ( t,s ) +\int
_{0}^{s}\partial_{s}K ( s,r )\,
\partial_{t}K ( t,r ) \,dr=:f ( s,t ).
\]
If $\mu:=f ( s,t ) \,ds \,dt$ defines a Radon measure on
$ [ 0,T%
] ^{2}\setminus D$, with $\mu_{-}$ having finite mass, part~\textup{A} is
applicable. Rather than imposing technical conditions on $K$, we verify this
in the model case of \textit{Volterra fBm}, $K ( t,s )
= (
t-s ) ^{H-1/2}, H>1/2$ (As above, there is nothing to do in the
Brownian case $H=1/2$ since then $f\equiv0$ and so $\mu\equiv0$.)
Specializing the above formula for $\partial_{s,t}R$, we have
\[
\partial_{s,t}R= ( H-1/2 ) ^{2}\int_{0}^{s}
( t-r ) ^{H-3/2} ( s-r ) ^{H-3/2}\,dr=:f ( s,t ) \geq0.
\]
Since $f$ remains bounded away from the diagonal, it clearly defines a
(nonnegative!) Radon measure. Trivially, $\mu_{-}\equiv0$, and so thanks
to part~\textup{A},
\[
V_{1} \bigl( R; [ s,t ] ^{2} \bigr) \leq
\sigma^{2} ( s,t ) = O\bigl(\llvert t-s\rrvert\bigr).
\]
\end{example}

\subsubsection{Examples handled by part~\textup{B}}\label{sec232}

%
\begin{example}[(Fractional Brownian motion $H\leq1/2$)]
\label{fbmrough} Consider a (standard) fractional Brownian motion $B^{H}$,
with $\sigma^{2} ( s,t ) =\llvert t-s\rrvert
^{2H}$ in the
regime $H\leq1/2$.
We compute $\mu=\partial_{s,t}R= ( -1/2 )\, \partial
_{s,t}\sigma
^{2}$ away from the diagonal and find
\[
\mu=-\mu_{-}=-H ( 1-2H ) \llvert t-s\rrvert^{2H-2} \,ds \,dt
\leq0
\]
which clearly yields a Radon measure on $ [ 0,T ]
^{2}\setminus D$. Note that $\mu_{+}\equiv0$ in the decomposition $\mu
=\mu_{+}-\mu
_{-}$. Conditions (B.ii) and (B.iii) with $\rho=1/ ( 2H ) $,
$\omega
( s,t ) =t-s$ are clear. It follows that the fBm covariance
function, $R ( s,t ) =\frac{1}{2} (
s^{2H}+t^{2H}-\llvert
t-s\rrvert ^{2H} ) $, has finite ``H\"{o}lder controlled''
mixed $%
( 1,\rho) $-variation, in the sense that
\[
V_{1,\rho} \bigl( R; [ s,t ] ^{2} \bigr) \leq O \bigl( \llvert
t-s\rrvert^{1/\rho}\bigr).
\]
\end{example}

%
\begin{example}[(Stationary increments II, Brownian and worse regularity)]
\label{StInII}Consider the case
\[
\sigma^{2} ( s,t ) =F \bigl( \llvert t-s\rrvert\bigr) \geq0,
\]
with $F$ continuous, nonnegative and with $F(0)=0$. A simple condition
on $%
F $ which generalizes at once the above fBm example and the previous
Example %
\ref{StInI} is \textit{semi-concavity}, that is,
\[
F^{\prime\prime}\leq k\qquad\mbox{in distributional sense on } ( 0,T
)\mbox{ for some } k\in\R,
\]
which is tantamount to say that $-F^{\prime\prime}+k$ is a (nonnegative)
Radon measure on $ ( 0,T ) $, which in turn induces a signed Radon
measure on $ [0,T ]^{2}\setminus D$, given by
\[
A\mapsto\int_{A}\bigl(-F^{\prime\prime} \bigl( \llvert t-s
\rrvert\bigr) +k\bigr) \,ds \,dt-k\lambda(A),
\]
where $\lambda$ is the two-dimensional Lebesgue measure. Then $\mu
:=\partial_{s,t}R= -\frac{1}{2} \partial_{s,t}\sigma^{2}$ is also a
signed Radon measure, with $\mu_{+}\leq\frac{k}{2}\lambda$.
Clearly, there will always be some $h>0$ (depending on $F$) such that (B.ii)
holds. Under the additional assumption $F ( t ) =O (
t^{1/\rho
} ) $ for some $\rho\geq1$, we then have (B.iii), with $\omega
(
s,t ) =C ( t-s ) $ and conclude that, with changing constants,
\[
V_{1,\rho} \bigl( R; [ s,t ] ^{2} \bigr) \leq C \biggl(
\llvert t-s\rrvert^{1/\rho}+\frac{k}{2}\llvert t-s\rrvert
^{2} \biggr) \le O \bigl( \llvert t-s\rrvert^{1/\rho}\bigr).
\]
\end{example}

%
\begin{example}[(Sums of fBm)]
In the previous example, $F^{\prime\prime}$ was bounded, as a Schwartz
distribution, by an $L^{\infty}$-function on $ [ 0,T ] ^{2}$,
namely by the constant~$k$. But $L^{1}$ would be enough. Consider $%
X=B^{H_{1}}+B^{H_{2}}$, a sum of two independent fBm with Hurst
parameters $%
H_{1}\geq1/2\geq H_{2}$.
A look at our two previous fBm examples reveals that
\[
\mu=\underbrace{H_{1} ( 2H_{1}-1 ) \llvert t-s\rrvert
^{2H_{1}-2} \,ds \,dt}_{=:\mu_{+}}- \underbrace{H_{2} (
1-2H_{2} ) \llvert t-s\rrvert^{2H_{2}-2} \,ds \,dt}_{=:\mu_{-}}.
\]
We easily check all conditions, in particular (B.iii) holds with $\rho
=1/ ( 2H_{2} ) \geq1$ and $\omega( s,t )
=t-s$. As a
consequence,
\begin{eqnarray*}
V_{1,\rho}\bigl(R;[s,t]^{2}\bigr) &\leq&C \biggl( \llvert t-s
\rrvert^{1/\rho
}+H_{1} ( 2H_{1}-1 ) \int
_{ [ s,t ] ^{2}}\bigl\llvert t'-s'\bigr
\rrvert^{2H_{1}-2} \,ds' \,dt^{\prime
} \biggr)
\\
&\le&C \bigl( \llvert t-s\rrvert^{1/\rho}+\llvert t-s\rrvert
^{2H_{1}} \bigr) =O \bigl( \llvert t-s\rrvert^{1/\rho
} \bigr).
\end{eqnarray*}
(Of course, the same conclusion can be obtained from our previous fBm
examples, using $R_{X}=R_{B^{H_{1}}}+R_{B^{H_{2}}}$ and then the triangle
inequality for the semi-norm $V_{1,\rho}$.)
\end{example}

%
\begin{example}[(Volterra processes II)]
Volterra fBm with $H<1/2$, that is, singular kernel $K ( t,s
) = (
t-s ) ^{H-1/2}$ is also covered by part \textup{B}. More generally, it is
possible (thanks to the robustness of the conditions of part \textup{B}), if tedious,
to give technical assumptions on $K$ which guarantee that (B.i)--(B.iii) are
satisfied. We note that $\rho\geq1$ of condition (B.iii) is determined
from the blow-up behavior of $K$ near the diagonal.
\end{example}

\subsubsection{Further examples handled by part \textup{B}}\label{sec233}

(This section may be skipped at first reading. In particular,
the reader may want to read Section~\ref{sec3} on random Fourier
series before looking in detail at the ``Fourier-based'' examples
below. Related applications to SPDEs are discussed in Section~\ref{secapplication}.)


%
\begin{example}[(Bifractional Brownian motion)]
\label{examplebifbm} Consider a \textit{bifractional Brownian
motion} (cf., e.g., \cite{HV03,RT06,KRT07}), that is, a centered
Gaussian process $B^{H,K}$
on $[0,T]$ with covariance function given by\footnote{%
As pointed out, for example, in \cite{KRT07} this process does not fit
in the
Volterra framework.}
\begin{eqnarray*}
&&R(s,t) = \frac{1}{2^K} \bigl( \bigl(s^{2H} + t^{2H}
\bigr)^K - \llvert t-s\rrvert^{2HK} \bigr),
\end{eqnarray*}
for some $H \in(0,1)$ and $K\in(0,1]$.
It is known (cf. \cite{HV03}, Proposition 3.1) that whenever $s<t$,
%
\begin{equation}
\label{eqnquasihelixpropbifbm} 2^{-K}\llvert t-s\rrvert^{2HK} \leq
\sigma^2(s,t) \leq2^{1-K} \llvert t-s\rrvert^{2HK}.
\end{equation}
We claim that the case $HK \geq\frac{1}{2}$ (resp., $\leq\frac
{1}{2}$) is
handled by part \textup{A} (resp., \textup{B}) of Theorem~\ref{teomain}. To this end, first
note that
\begin{eqnarray*}
\partial_{s,t} R(s,t) &=& \frac{(2H)^2K(K-1)}{2^K} \frac{s^{2H-1}
t^{2H-1}}{%
(s^{2H} + t^{2H})^{2-K}}
\\
&&{} +
\frac{2HK(2HK - 1)}{2^K} \llvert t-s\rrvert^{2HK - 2}.
\end{eqnarray*}
The measure
\begin{eqnarray*}
&&\nu:= - \frac{(2H)^2K(K-1)}{2^K} \frac{s^{2H-1} t^{2H-1}}{(s^{2H} +
t^{2H})^{2-K}} \,ds \,dt
\end{eqnarray*}
has finite mass. Indeed, it is enough to show that
\begin{eqnarray*}
&&\int_{B_{\delta}(0)} \frac{\llvert st\rrvert ^{2H-1}}{(\llvert
s\rrvert ^{2H} +
\llvert t\rrvert ^{2H})^{2-K}} \,ds \,dt
\end{eqnarray*}
is finite for some $\delta> 0$, where $B_{\delta}(0)$ denotes the closed
ball around $0$ with radius~$\delta$. Introducing polar coordinates, this
integral equals
%
\begin{eqnarray}
\label{eqnestimatekernelbfbm}
&& \int_{0}^{\delta} \int
_{0}^{2\pi} r^{2HK - 1} \frac{\llvert \sin(\theta)
\cos(\theta)\rrvert ^{2H - 1}}{(\llvert \sin(\theta)\rrvert ^{2H} +
\llvert \cos(\theta)\rrvert ^{2H})^{2-K}} \,d
\theta \,dr
\nonumber\\[-8pt]\\[-8pt]\nonumber
&&\qquad
\leq2^{1 - 2H} \int_{0}^{2\pi} \bigl
\llvert\sin(2\theta)\bigr\rrvert^{2H
- 1} \,d\theta\int_0^{\delta}
r^{2HK - 1} \,dr
\end{eqnarray}
and both integrals are finite for $H,K > 0$. Note that estimate \eqref
{eqnestimatekernelbfbm} also implies that $\nu([s,t]^2) \leq
C\llvert t-s\rrvert ^{2HK}$ for some constant $C$ depending of $H$,
$K$ and $T$.

Hence we obtain that $\partial_{s,t} R(s,t):= \mu$ is a Radon
measure on $%
(0,T)^2 \setminus D$. If $HK \geq\frac{1}{2}$, we have the
decomposition $%
\mu= \mu_+ - \mu_{-}$ with $\mu_{-} = \nu$, and we have already
seen that
(A.i) holds. (A.ii) is trivially satisfied, and with %
\eqref{eqnquasihelixpropbifbm} we may conclude that
%
\begin{eqnarray}
V_1\bigl(R;[s,t]^2\bigr) \leq2^{1 - K} \llvert
t-s\rrvert^{2HK} + 2\nu\bigl([s,t]^2\bigr)
\leq C\llvert t-s\rrvert^{2HK}\nonumber
\\
\eqntext{\mbox{for all } [s,t] \subseteq[0,T].}
\end{eqnarray}
If $HK \leq\frac{1}{2}$, $\mu_+ \equiv0$ on $(0,T)^2 \setminus D$,
thus (B.i) is satisfied in both cases. (B.ii) is also easy to see. Indeed,
since $B^{H,K}$ is a self-similar process with index $HK$, one can use
scaling to see that it is enough to show that for all $t_0 \in\R_+$
and $%
h_0 \in[0,1]$, the function
\begin{eqnarray*}
&& h \mapsto R\pmatrix{ t_0,t_0 + 1
\cr
t_0 + h_0, t_0 + h_0 + h}
=: \phi(h)
\end{eqnarray*}
is nonnegative on $[0, 1-h_0]$. Since $\phi(0) = 0$, it is enough to show
that $\phi'\geq0 $ on $(0, 1-h_0)$ which follows by a simple
calculation. Finally, from~\eqref{eqnquasihelixpropbifbm} we see that
(B.iii) holds with $\rho= \frac{1}{2HK}$ and $\omega(s,t) = \llvert
t-s\rrvert $,
therefore
\begin{eqnarray*}
&& V_{1,\rho}\bigl(R;[s,t]^2\bigr) = O \bigl( \llvert t-s\rrvert
^{1/\rho} \bigr).
\end{eqnarray*}
\end{example}

%
\begin{example}[(Random Fourier series I: stationary)]
\label{exrfsstationary} Consider a stationary random Fourier series%
\footnote{%
We may ignore the (constant, random) zero-mode in the series since we are
only interested in properties of the increments of the process.}
\[
\Psi( t ) =\sum_{k=1}^{\infty}
\a_kY^k \sin( kt ) + \a%
_{-k}Y^{-k}
\cos( kt ),\qquad t\in[ 0,2\pi],
\]
with zero-mean, independent Gaussians $ \{ Y^{k} \mid k\in\Z
\} $
with unit variance. We compute
\begin{eqnarray*}
R ( s,t ) &=&\sum\a_k^2 \sin(ks)\sin(kt) +
\alpha_{-k}^{2}\cos(ks)\cos(kt)
\\
&=&\frac{1}{2}\sum\bigl( \alpha_{k}^{2}+
\a_{-k}^{2} \bigr) \cos\bigl(k(t-s)\bigr)+ \bigl(
\alpha_{k}^{2}-\a_{-k}^{2} \bigr) \cos
\bigl(k(t+s)\bigr)
\end{eqnarray*}
and note that $\alpha_{k}^{2}\equiv\a_{-k}^{2}$ due to the assumed
stationarity of $\Psi$. This leaves us with
\begin{eqnarray*}
R ( s,t ) &=&K \bigl( \llvert t-s\rrvert\bigr),
\\
\sigma^{2} ( s,t ) &=&2 \bigl( K ( 0 ) -K \bigl( \llvert t-s\rrvert
\bigr) \bigr) =:F \bigl( \llvert t-s\rrvert\bigr),
\end{eqnarray*}
where
\[
K(t):=\sum_{k=1}^{\infty}\alpha_{k}^{2}
\cos( kt ).
\]
In special situations, for example, when $\alpha_{k}^{2}=1/k^{2}$, one
can find $%
K\in\break C^{2} ( [ 0,2\pi] ) $ in closed form, which
brings us back to Example~\ref{StInI}. This is not possible in
general, but
in view of Example~\ref{StInII} above, it would suffice to know that $K$
is convex and $1/\rho$-H\"{o}lder. Conditions on the Fourier-coefficients
for this to hold true are known from Fourier analysis (recalled in
detail in
Section~\ref{secrfs} below). For instance, given (eventually)
decreasing $%
( \alpha_{k}^{2} ) $, $K$ is $1/\rho$-H\"{o}lder if and
only if $\alpha
_{k}^{2}=O ( k^{-(1+1/\rho)} ) $. In particular,
in the
model case
\[
\alpha_{k}^{2}=\frac{1}{k^{2\alpha}},
\]
the desired decay holds true if and only if
\[
2\alpha=1+1/\rho\leftrightarrow\rho=\frac{1}{2\alpha-1}\geq1\qquad
(\mbox{for }
\alpha\leq1).
\]
Convexity also holds true here and we conclude that for all $ [ s,t
] \subset[ 0,2\pi] $,
\[
V_{1,\rho} \bigl( R; [ s,t ] ^{2} \bigr) = O \bigl( \llvert
t-s\rrvert^{1/\rho} \bigr).
\]
\end{example}

%
\begin{example}[(Random Fourier series II: nonstationary, general
case)]\label{exrfsnon-stationary}
As seen in the previous example,
the covariance may be written as%
%
\begin{eqnarray}
\label{eqnfourierseriesdecomp} R ( s,t ) & =&K \bigl( \llvert t-s\rrvert
\bigr) +K \bigl( \llvert
t+s\rrvert\bigr) +\td K \bigl( \llvert t-s\rrvert\bigr) -\td K \bigl(
\llvert
t+s\rrvert\bigr)
\\
& =:& R^{-} ( s,t ) +R^{+} ( s,t ) +\td R^{-} (
s,t ) -\td R^{+} ( s,t ),
\end{eqnarray}
where $R^{\pm}$ and $K$ are as before and
\[
\td K(t):=\sum_{k=1}^{\infty}
\alpha_{-k}^{2}\cos( kt ).
\]
Under the assumption that $K,\td K$ are convex and $1/\rho$-H\"
{o}lder, the
cases $R\in\{R^{-},\td R^{-}\}$
were already handled in the previous example, where we established
\[
V_{1,\rho} \bigl( R; [ s,t ] ^{2} \bigr) =O \bigl( \llvert t-s
\rrvert^{1/\rho
} \bigr).
\]
We claim that $R^{+}$ %
can be handled with part \textup{A}. $\td R^{+}$ may then be treated analogously.
Condition (A.i) is simple: using convexity of $K$,
\[
\partial_{s,t}R^{+}=K^{\prime\prime} ( t+s ) \geq0\qquad
\mbox{on } [ 0,T ] ^{2}\setminus D,
\]
so that $\mu:=\partial_{s,t}R^{+}=\mu_{+}$ is a nonnegative (but in
general not finite) Radon-measure on $ [ 0,T ]
^{2}\setminus D$.
Unlike in previous examples, condition (A.ii) is not trivial, since $R^{+}$
is not a covariance function in general. Nonetheless, we have
%
\begin{eqnarray}
R^{+}\pmatrix{ s,t
\cr
s,t} &=&K ( 2t ) +K ( 2s ) -\bigl(2K ( t+s )
\bigr)
\\
&=&2 \biggl( \frac{K ( 2t ) +K ( 2s ) }{2}-K \biggl( \frac{%
2t+2s}{2} \biggr) \biggr)
\\
&\geq&0\qquad\forall0\leq s\leq t\leq\pi,\label{eqncondAiiforrfs}
\end{eqnarray}
thanks to convexity of $K$ on $ [ 0,2\pi] $. This settles
condition (A.ii). We conclude that $R^{+}$ has finite
$1$-variation,\footnote{%
The situation here is reminiscent of absolutely continuous paths
$x=x (
t ) $ on $ [ 0,T ] $ with $\dot{x}\in L^{p}$ where
$1/\rho
+1/p=1$. Indeed, as may be seen from H\"{o}lder's inequality, the
$L^{1}$%
-norm of $\dot{x}| _{ [ s,t ]
}$, which equals
the $1$-variation of $x$ over $ [ s,t ] $, is finite and of
order $%
\llvert t-s\rrvert ^{1/\rho}$.}
%
\begin{eqnarray}\label{eqnpospartrfshoelder}
V_{1} \bigl( R^{+}; [ s,t ] ^{2}
\bigr) &\leq& R \pmatrix{ s,t
\cr
s,t} =K ( 2t ) +K ( 2s ) -2K ( t+s )
\nonumber\\[-8pt]\\[-8pt]\nonumber
&=& O \bigl(\llvert t-s\rrvert^{1/\rho} \bigr).
\end{eqnarray}
Since $R=R^{-}+R^{+}+\td R^{-}-\td R^{+}$, we can now conclude with $%
V_{1,\rho}\leq V_{1}$ and the triangle inequality to see that $R$ has
(H%
\"{o}lder controlled) mixed $ ( 1,\rho) $-variation, in
the sense
that
\[
V_{1,\rho} \bigl( R; [ s,t ] ^{2} \bigr) =O \bigl( \llvert t-s
\rrvert^{1/\rho} \bigr),
\]
for all $ [ s,t ] \subset[ 0,\pi] $. (The
extension of
this estimate to $[0,2\pi]$ is not difficult.\footnote{%
Considering the Fourier series with argument shifted by $\pi$, gives the
same estimate on $[\pi,2\pi]^{2}$. In fact, one can also handle the
mixed $%
(1,\rho)$-variation of $R^{+}$ on $[0,\pi]\times [\pi,2\pi]$ by
playing it back to the mixed variation of $R^{-}$ on $[0,\pi]\times
{}[
0,\pi]$, using the fact that $K$ is given by cosine series, hence is even
around $\pi$.})
\end{example}

%
\begin{example}[(Fourier fractional Brownian bridge)]\label
{exfourierFBB} Fourier fractional Brownian bridge is the Gaussian
process given by the random Fourier series
\[
W_{t}^{\alpha}=\sum_{k=1}^{\infty}
\frac{Y_{k}\sin( (k/2) t ) }{k^{\alpha}}\qquad\mbox{for }t\in{}[0,2\pi], \alpha\in\biggl(%
\frac{1}{2},1\biggr],
\]
with $Y_{k}$ as above. This process arises by replacing the covariance
operator of Brownian bridge (the Dirichlet Laplacian $-\D$) by its
fractional power $(-\D)^{\a}$. Clearly, this is a special case of the
previous example.
\end{example}

%
\begin{example}[(Stationary processes: spectral measure)]\label
{exspectralmeasure}
Let $X_{t}$ be a stationary, zero-mean process with covariance
\[
R ( s,t ) =K \bigl( \llvert t-s\rrvert\bigr)
\]
for some continuous function $K$. By a well-known theorem of Bochner,
\begin{eqnarray*}
K ( t ) &=&\int\cos( t\xi) \mu( d\xi),
\\
\sigma^{2} ( t ) &:=&  \sigma^2(0,t) =2 \bigl( K ( 0 ) -K
( t ) \bigr) = 4 \int\sin^{2} ( t\xi/2 ) \mu( d\xi),
\end{eqnarray*}
where $\mu$ is a finite positive symmetric measure on $\mathbb{R}$ (``spectral
measure''). The case of discrete $\mu$ corresponds to Example~\ref
{exrfsstationary}. Another example is given by the \textit{fractional Ornstein--Uhlenbeck process},
\[
X_{t}=\int_{-\infty}^{t}e^{-\lambda( t-u )
}\,dB_{u}^{H},\qquad
t\in\mathbb{R},
\]
which should be viewed as the stationary solution to $dX=-\lambda
X\,dt+dB^{H}$. In this case, it is known that $X$ has a \emph{spectral density}
of the form\footnote{%
This generalizes the well-known fact that the spectral density of the
classical OU process is of Cauchy type.}
\[
\frac{d\mu}{d\xi}=c_{H}\frac{\llvert \xi\rrvert
^{1-2H}}{\lambda
^{2}+\xi^{2}}.
\]
Clearly, the decay of the density is related to the regularity of $K$. More
precisely, writing
\begin{eqnarray*}
\hat{K} ( \xi) &:=& \frac{\llvert \xi\rrvert
^{1-2H}}{%
\lambda^{2}+\xi^{2}}\sim\langle\xi\rangle^{-1-2H}
\qquad\mbox{where } \langle\xi\rangle= \bigl( 1+\xi^{2} \bigr)
^{1/2},
\\
\langle\xi\rangle^{s}\hat{K} ( \xi) &\sim&\langle\xi
\rangle^{s-1-2H}\in L^{2}\qquad\mbox{iff }2 ( s-1-2H ) <-1,
\end{eqnarray*}
that is, if and only if $s < s^{\ast}:=1/2+2H$. It follows that $K\in
H^{s}$ for any $s<s^{\ast}$ and
thus by a standard Sobolev embedding,
$K$ is $\alpha$-H\"{o}lder for $\alpha<s^{\ast}-1/2=2H$. Alternatively,
and a little sharper, Theorem~7.3.1 in \cite{MR06} tells us that if
$\hat{K}$ is regularly varying at $\infty$, then
\[
\sigma^{2} ( t ) \sim C\hat{K} ( 1/t ) /t\qquad\mbox{as } t
\rightarrow0.
\]
Applied to the situation at hand we see that $\sigma^{2} (
t ) =O (
t^{2H} )$, since $\hat{K} ( \xi) \sim(
1/\xi) ^{1+2H}$. With focus on the rough case $H\leq1/2$, this gives
condition (B.iii) with $\rho=1/ ( 2H ) $, $\omega(
s,t ) =t-s$. Moreover, it can be seen that there is a $T > 0$
such that
$K$ is convex on $[0,T]$ (cf. Example~\ref{exfracOUcalc} below),
which implies (B.i) and
(B.ii) as in Example~\ref{StInII}. Hence it follows that $V_{1,\rho
} (
R; [ s,t ] ^{2} ) =O ( \llvert t-s\rrvert
^{2H} ) $ for all $[s,t] \subseteq[0,T']$.
\end{example}

%

\subsection{Proof of Theorem \texorpdfstring{\protect\ref{teomain}}{2.2}, part \textup{A}}\label{sec24}
From (A.i), the distributional mixed derivative of $R$ on $ (
0,T )
^{2}\setminus D$ is given by
%
\begin{equation}
\frac{\partial^{2}R}{\partial s\,\partial t}=\mu_{+}-\mu_{-}, \label{eqnRdecomp1}
\end{equation}
where $\mu_{-}$ (trivially extended to $ [ 0,T ] ^{2}$ whenever
convenient) has finite mass. By assumption, the distribution function
of $%
\mu_{-}$
\[
R^{-} ( s,t ):=\mu_{-} \bigl( [ 0,s ]\times [ 0,t]
\bigr),
\]
is continuous. We may then define $R^{+}\in C( [ 0,T ]) $ by
imposing the decomposition
\[
R=R^{+}-R^{-}.
\]
Clearly, the distributional mixed derivatives of $R^{\pm}$ on $ (
0,T ) ^{2}\setminus D$ are given by
%
\begin{equation}
\frac{\partial^{2}R^{\pm}}{\partial_{t}\,\partial_{s}}=\mu_{\pm}. \label
{eqnRdecomp2}
\end{equation}
Noting that all rectangular increments of $R^{-}$ are nonnegative, $%
R^{-} ( A ) = \mu_{-} ( A ) \geq0$, we
immediately have
\[
V_{1}\bigl(R^{-};A\bigr)=R^{-} ( A ) =
\mu_{-} ( A )
\]
for all $A=[s,t]\times [ u,v]\subset[ 0,T ] ^{2}$. On the
other hand, any such rectangle $A$ may be split up in finitely many ``small
squares,'' say $Q_{i}= [ t_{i},t_{i+1} ] ^{2}$ with $%
t_{i+1}-t_{i}\leq h$ for all $i$, and a (finite) number of ``off-diagonal''
rectangles $A_{j}$, whose interior does not intersect the diagonal.
Since $%
R ( Q_{i} ) =\sigma^{2} ( t_{i},t_{i+1} ) \geq
0$, by
(A.ii), and $R ( A_{j} ) \geq-R^{-} ( A_{j} )
=-\mu
_{-} ( A_{j} ) $, we have
\begin{eqnarray*}
R ( A ) &=& \sum_i R ( Q_{i} ) +\sum
_{j}R ( A_{j} )
\\
&\geq& -\sum_{j}\mu_{-} (
A_{j} ) \geq-\mu_{-} ( A ),
\end{eqnarray*}
for all rectangles $A$. This implies finite $1$-variation over every
rectangle $A=[s,t]\times [ u,v]$. Indeed, for any dissections $(t_{i})$
of $[s,t]$ and $(t_{j}')$ of $[u,v]$ we have
\begin{eqnarray*}
\sum_{t_{i},t_{j}'}\left| R\pmatrix{
t_{i},t_{i+1}
\cr
t_{j}',t_{j+1}'}
\right|& \leq&\sum_{t_{i},t_{j}'} \biggl\{ \biggl\llvert R
\pmatrix{ t_{i},t_{i+1}
\cr
t_{j}',t_{j+1}'}+\mu_{-}\bigl([t_{i},t_{i+1}]\times \bigl[
t_{j}^{\prime
},t_{j+1}'\bigr]\bigr)
\biggr\rrvert
\\
&&\hspace*{89pt} +\mu_{-}\bigl([t_{i},t_{i+1}]\times \bigl[ t_{j}',t_{j+1}'\bigr]
\bigr) \biggr\}
\\
& =&R\pmatrix{s,t
\cr
u,v} +2\mu_{-}\bigl([s,t]\times [ u,v]\bigr),
\end{eqnarray*}
and so, for all rectangles $A$,
\[
V_{1}(R;A)\leq R ( A ) +2\mu_{-}(A).
\]
%

\subsection{Proof of Theorem \texorpdfstring{\protect\ref{teomain}}{2.2}, part \textup{B}}\label{sec25}

Let us start with a few definitions.

\begin{definition}\label{defother2dvariations}
For $\gamma, \rho\geq1$ set
\begin{eqnarray*}
&& V_{\gamma, \rho}^+ \bigl(R;[s,t]\times[u,v]\bigr)
\\[-2pt]
&&\qquad := \sup_{(t'_{j})
\in
\mathcal{D}([u,v])}
\biggl( \sum_{t'_j} \sup_{(t_{i}) \in
\mathcal{D}%
([s,t])}
\biggl( \sum_{t_i} \biggl\llvert R \pmatrix{
t_{i},t_{i+1}
\cr
t_{j}',t_{j+1}'}
\biggr\rrvert^{\gamma} \biggr)^{\rho/\gamma} \biggr)^{1/\rho}
\end{eqnarray*}
and
\begin{eqnarray*}
V_{\gamma, \rho}^+ (R;U_{[s,t]}) &:=& \sup_{(t'_{j}) \in
\mathcal{D}%
([s,t])}
\biggl( \sum_{t'_j} \sup_{(t_{i}) \in\mathcal{D}%
([s,t'_j])}
\biggl( \sum_{t_i} \biggl\llvert R \pmatrix{
t_{i},t_{i+1}
\cr
t_{j}',t_{j+1}'}
\biggr\rrvert^{\gamma} \biggr)^{\rho/\gamma} \biggr)^{1/\rho},
\\[-2pt]
V_{\gamma, \rho}^+ (R;L_{[s,t]}) &:=& \sup_{(t'_{j}) \in
\mathcal{D}%
([s,t])}
\biggl( \sum_{t'_j} \sup_{(t_{i}) \in\mathcal{D}%
([t'_{j+1},t])}
\biggl( \sum_{t_i} \biggl\llvert R \pmatrix{
t_{i},t_{i+1}
\cr
t_{j}',t_{j+1}'}
\biggr\rrvert^{\gamma} \biggr)^{\rho/\gamma} \biggr)^{1/\rho},
\\[-2pt]
V_{\gamma, \rho}^+(R;D_{[s,t]}) &:=& \sup_{(t'_{j}) \in
\mathcal{D}%
([s,t])}
\biggl( \sum_{t'_j} \sup_{(t_{i}) \in\mathcal{D}%
([t'_j,t'_{j+1}])}
\biggl( \sum_{t_i} \biggl\llvert R
\pmatrix{t_{i},t_{i+1}
\cr
t_{j}',t_{j+1}'}
\biggr\rrvert^{\gamma} \biggr)^{\rho/\gamma} \biggr)^{1/\rho}.
\end{eqnarray*}
%
\end{definition}

For any rectangle $A\subseteq I^2$ it is easy to see that
\begin{eqnarray*}
&& V_{\gamma, \rho}(R;A) \leq V_{\gamma, \rho}^+ (R;A)
\end{eqnarray*}
and also (e.g., as a consequence of \cite{FV11}, Theorem 1.i)
\begin{eqnarray*}
&& V_{1}(R;A) = V_{1,1}^+(R;A).
\end{eqnarray*}

The main reason for introducing $V^+$ as above is the following lemma:

\begin{lemma}[(Concatenation Lemma 1)]
\label{lemconc1} Let $R$ be as before. Then
%
\begin{eqnarray}
V_{\gamma, \rho}^{+}\bigl(R; [s,t]^2\bigr) &\leq& C
\bigl( V_{\gamma, \rho}^+ (R;U_{[s,t]}) + V_{\gamma, \rho}^+(R;D_{[s,t]})
+ V^+_{\gamma,\rho}(R;L_{[s,t]}) \bigr)\nonumber
\\[-2pt]
\eqntext{\forall[s,t] \subseteq I,}
\end{eqnarray}
for some constant $C=C(\rho,\gamma)$.
\end{lemma}

\begin{pf}
Let\vspace*{1pt} $(t'_j)$ be a partition of $[s,t]$. Fix $[t^{\prime
}_j,t'_{j+1}]$, and let $(t_i)$ be a partition of $[s,t]$. By
subdividing rectangles which lie on the diagonal into at maximum three
parts, we see that
\begin{eqnarray*}
&& 3^{1 - \gamma} \sum_{t_i} \biggl\llvert R
\pmatrix{ t_{i},t_{i+1}
\cr
t_{j}',t_{j+1}'}
\biggr\rrvert^{\gamma}
\\[-2pt]
&&\qquad \leq\sup_{(t_{i}) \in\mathcal{D}([s,t'_j])} \sum_{t_i} \biggl
\llvert R \pmatrix{ t_{i},t_{i+1}
\cr
t_{j}',t_{j+1}'}
\biggr\rrvert^{\gamma} + \sup_{(t_{i}) \in\mathcal
{D}([t^{\prime
}_j,t'_{j+1}])} \sum
_{t_i} \biggl\llvert R \pmatrix{ t_{i},t_{i+1}
\cr
t_{j}',t_{j+1}'} \biggr
\rrvert^{\gamma}
\\[-2pt]
&&\quad\qquad{} + \sup_{(t_{i}) \in\mathcal{D}([t'_{j+1},t])}\sum_{t_i} \biggl
\llvert R \pmatrix{ t_{i},t_{i+1}
\cr
t_{j}',t_{j+1}'}
\biggr\rrvert^{\gamma}.
\end{eqnarray*}
Now we take the supremum, sum over $t'_j$ and take the
supremum again
to see that
\begin{eqnarray*}
&&\sup_{(t'_{j}) \in\mathcal{D}([s,t])} \biggl( \sum_{t'_j}
\sup_{(t_{i}) \in\mathcal{D}([s,t])} \biggl( \sum_{t_i}
\biggl\llvert R \pmatrix{t_{i},t_{i+1}
\cr
t_{j}',t_{j+1}'}\biggr\rrvert
^{\gamma} \biggr)^{\rho/\gamma} \biggr)^{1/\rho}
\\
&&\qquad \leq C \bigl( V_{\gamma, \rho}^+ (R;U_{[s,t]}) + V_{\gamma, \rho}^+
(R;D_{[s,t]}) + V^+_{\gamma, \rho}(R;L_{[s,t]}) \bigr).
\end{eqnarray*}\upqed
\end{pf}

\begin{lemma}[(Concatenation Lemma 2)]
\label{lemconc2} Assume that there is an $h>0$ such that
\begin{eqnarray*}
&& V_{\gamma, \rho} \bigl(R; [s,t] \times[u,v]\bigr) \leq\Phi(s,t;u,v)
\end{eqnarray*}

\noindent
holds for all squares $[s,t]\times[u,v]=[s,t]^2 \subseteq D_h$ and all
off-diagonal rectangles $(s,t)\times(u,v) \subseteq I^2\setminus{D}$,
where $%
\Phi\colon\Delta_I \times\Delta_I \to\R_+$ is a nondecreasing function
in $t-s$ and $v-u$. Then there is a constant $C = C(\gamma,\rho,h,T)$ such
that
\begin{eqnarray*}
&& V_{\gamma, \rho} \bigl(R; [s,t] \times[u,v]\bigr) \leq C \Phi(s,t;u,v)
\end{eqnarray*}
holds for all rectangles $[s,t]\times[u,v]$. The constant $C$ can be chosen
independently of $h$ and $T$ when considering only rectangles $[s,t]
\times
[u,v] \subset D_h$. The same is true if one replaces $V_{\gamma, \rho
}$ by $%
V_{\gamma, \rho}^+$.
\end{lemma}

\begin{pf}
\textit{Step} 1. Consider any square of the form $[s,t]^2 \subseteq I^2$.
Then\vspace*{1pt} we can subdivide this square into $N^2$ smaller squares
$(A_{i,j})_{i,j =
1}^N$ with equal side length~$\tilde{h}$, which can be chosen such
that $h/2
\leq\tilde{h} \leq h$ and $N \leq M$ where $M$ is a number depending
on $T$
and $h$; see Figure~\ref{figfigure1}. Each of these small squares
does either lie on the diagonal, or its
inner part does not intersect with the diagonal. Hence
\begin{eqnarray*}
&& V_{\gamma,\rho}\bigl(R;[s,t]^2\bigr) \leq c_1(N,
\gamma,\rho) \sum_{i,j = 1}^N
V_{\gamma,\rho}(R;A_{i,j}) \leq c_2(N,\gamma,\rho)
\Phi(s,t;u,v)
\end{eqnarray*}
by monotonicity.

%
\begin{figure}[t]

\includegraphics{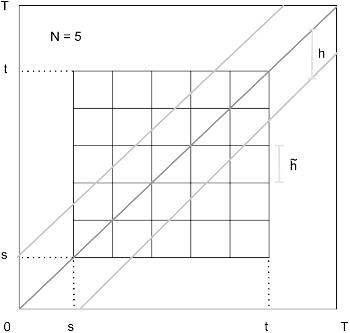}

\caption{Subdivision of square as used in step~1 of Lemma \protect\ref{lemconc2}.}\label{figfigure1}
\end{figure}

%
\begin{figure}[b]

\includegraphics{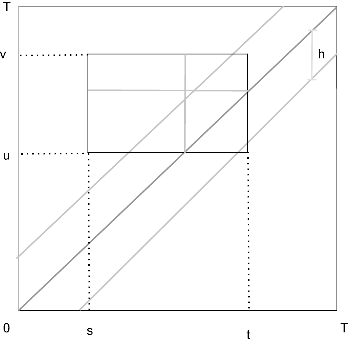}

\caption{Subdivision of square as used in step 2 of Lemma \protect\ref{lemconc2}.}\label{figfigure2}
\end{figure}

\textit{Step} 2. Let $[s,t]\times[u,v]$ be any rectangle in $I^2$.
Then we
can subdivide it into one square lying on the diagonal and three rectangles
for which the inner part does not intersect with the diagonal; see
Figure~\ref{figfigure2}. We conclude
as in step~1.
\end{pf}

\begin{lemma}
\label{lemincrcorr} Let $R$ as before and $\s$ as in \eqref{eqnsigma}.
Then the following two assertions are equivalent:
\begin{longlist}[(iii)]
\item[(i)] $\frac{\partial^2 \s^2}{\partial_t\, \partial_s}=-2\frac
{\partial^2
R}{\partial_t\, \partial_s} \ge0$ in the sense of distributions, that
is, for
every nonnegative $\phi\in C^\infty_c(I^2\setminus D)$,
\begin{eqnarray*}
&& \int_{I^2} \frac{\partial^2 \phi}{\partial_t\, \partial_s}(s,t) \s%
^2(s,t)
\,ds \,dt = -2\int_{I^2} \frac{\partial^2 \phi}{\partial_t\,
\partial_s}%
(s,t)
R(s,t) \,ds \,dt \ge0.
\end{eqnarray*}

\item[(ii)] For all off-diagonal rectangles $(s,t) \times(u,v)
\subseteq
I^2 \setminus D$, we have
\begin{eqnarray*}
&& R\pmatrix{ s,t
\cr
u,v} \leq0.
\end{eqnarray*}
\end{longlist}

In addition, if either of the above conditions is satisfied, then
\begin{eqnarray*}
&& R\pmatrix{ s,t
\cr
u,v}\le\sigma^2(u,v)\qquad\forall[u,v]
\subseteq[s,t] \subseteq I.
\end{eqnarray*}
All assertions remain true if we substitute $\leq$ by $\geq$ in the three
inequalities.
\end{lemma}

\begin{pf}
We will only consider the $\leq$-case. Let $\vp\in C_c(B_1(0))$
nonnegative with $\llVert \vp\rrVert _{L^1(\R^2)}=1$. We then define the
standard Dirac
sequence\break $\vp^\ve((s,t)):=\ve^{-2}\vp(\frac{1}{\ve}(s,t))$ and
observe $\supp%
(\vp^\ve)\subseteq B_\ve(0)$. We extend $R$ by $0$ to all of $\R^2$
and set $%
R^\ve:= R \ast\vp^\ve$. Then
\[
\frac{\partial^2 R^\ve}{\partial_t\, \partial_s}(a,b)=\int_{I^2} R(s,t)
\frac{%
\partial^2 }{\partial_t\, \partial_s}
\bigl(\vp^\ve(s-a,t-b)\bigr) \,ds \,dt.
\]
For $(a,b) \in\mathop{\D}\limits^{\circ}{}_I = \{(s,t)\mid s < t
\in
\mathop{I}\limits^{\circ}\}$, we note
\[
\supp\bigl(\vp^\ve(s-a,t-b)\bigr)\subseteq B_\ve
\bigl((a,b)\bigr) \subseteq\mathop{\D}\limits
^{\circ}{}_I
\]
for all $\ve$ small enough. Hence, $\vp^\ve(\cdot-a,\cdot-b)$ is an
admissible test-function for:
\begin{longlist}[(ii)]
\item[(i)] and thus $\frac{\partial^2 R^\ve}{%
\partial_t\, \partial_s} (a,b)\le0$. Since
\begin{eqnarray*}
&& R^\ve\pmatrix{ s,t
\cr
u,v}= \int_{I^2}
\mathbbm{1}_{[s,t]}(x) \mathbbm{1}_{[u,v]}(y)
\frac{%
\partial^2 R^\ve}{\partial_t\, \partial_s}(x,y) \,dx \,dy\qquad\forall s \le
t \le u \le v \in I,
\end{eqnarray*}
\item[(ii)] follows using continuity of $R$.
\end{longlist}

Suppose now that (ii) is satisfied. We may approximate $R$ by $R^\ve
\in
C^\infty(\D_I)$ such that
\[
\bigl\llVert R-R^\ve\bigr\rrVert_{C(\D_I)}\le
\frac{\ve}{4}.
\]
By (ii) we have
\begin{eqnarray*}
&& \int_{I^2} \mathbbm{1}_{[s,t]}(x)
\mathbbm{1}_{[u,v]}(y) \frac
{\partial^2 R^%
\ve}{\partial_t\, \partial_s}(x,y) \,dx \,dy = R^\ve
\pmatrix{ s,t
\cr
u,v} \leq\ve,
\end{eqnarray*}
for all $s \le t \le u \le v \in I$. We note that the set of
nonnegative $%
f\in L^1(\D_I)$ satisfying
\[
\int_{\D_I} f(x,y) \frac{\partial^2 R^\ve}{\partial_t\, \partial_s} 
(x,y) \,dx \,dy
\leq\ve
\]
is a monotone class. By the monotone class theorem, we thus have
\[
\int_{\D_I} f(x,y) \frac{\partial^2 R^\ve}{\partial_t\, \partial_s} 
(x,y) \,dx \,dy
\leq\ve
\]
for all nonnegative $f \in L^1(\D_I)$. Considering nonnegative $f\in
C^\infty_c(\mathop{\D}\limits^{\circ}{}_I)$, a partial integration
and letting $\ve\to0$ yields
(i).

To prove the remaining inequality we note
\begin{eqnarray*}
&& R\pmatrix{ s,t
\cr
u,v} = R\pmatrix{ s,u
\cr
u,v}+R \pmatrix{ u,v
\cr
u,v}+R
\pmatrix{ v,t
\cr
u,v} \le R\pmatrix{ u,v
\cr
u,v}.
\end{eqnarray*}\upqed
\end{pf}

We are now able to prove part \textup{B} of our main theorem.

\begin{pf*}{Proof of Theorem \protect\ref{teosuffcritmixed2dvar},
part \textup{B}}
We decompose $R$ as in \eqref{eqnRdecomp1}, \eqref{eqnRdecomp2}.
We start by proving \eqref{eqnon-diag-est}
by an application of Lemma~\ref{lemconc1}: let $(t_{j}')$
be a
partition of $[s,t]$. Fix $[t_{j}',t_{j+1}']$, and
let $%
(t_{i})$ be a partition of $[s,t_{j}']$. Apply Lemma~\ref{lemincrcorr} with $R$ equal to $-R^{-}$ and then $-R^{+}$ to get
\[
-R^{-}(A_{i,j})\leq0\leq R^{+}(A_{i,j})
= \mu_+(A_{i,j})
\]
for all $A_{i,j}=[t_{j}',t_{j+1}']\times [
t_{i},t_{i+1}]$. Hence, with condition (B.ii) we have
\begin{eqnarray*}
\sum_{t_{i}}\biggl\llvert R\pmatrix{
t_{i},t_{i+1}
\cr
t_{j}',t_{j+1}'}
\biggr\rrvert& \leq&\sum_{t_{i}}\biggl\llvert
R^{-} \pmatrix{ t_{i},t_{i+1}
\cr
t_{j}',t_{j+1}'}\biggr\rrvert+
\biggl\llvert R^{+} \pmatrix{ t_{i},t_{i+1}
\cr
t_{j}',t_{j+1}'}\biggr\rrvert
\\
& =&R^{-}\pmatrix{ s,t_{j}'
\vspace*{2pt}\cr
t_{j}',t_{j+1}'}
+R^{+}\pmatrix{ s,t_{j}'
\vspace*{2pt}\cr
t_{j}',t_{j+1}'}
\\
& =&-R\pmatrix{s,t_{j}'
\vspace*{2pt}\cr
t_{j}',t_{j+1}'}
+2R^{+}\pmatrix{s,t_{j}'
\vspace*{2pt}\cr
t_{j}',t_{j+1}'}
\\
& =&-R\pmatrix{ s,t_{j+1}'
\vspace*{2pt}\cr
t_{j}',t_{j+1}'}
+R\pmatrix{ t_{j}',t_{j+1}'
\vspace*{2pt}\cr
t_{j}',t_{j+1}'}
+2R^{+}\pmatrix{s,t_{j}'
\vspace*{2pt}\cr
t_{j}',t_{j+1}'}
\\
& \leq&\sigma^{2}\bigl(t_{j}',t_{j+1}'
\bigr)+2R^{+}\pmatrix{ s,t_{j}'
\vspace*{2pt}\cr
t_{j}',t_{j+1}'}
\\
& \leq&\omega\bigl(t_{j}',t_{j+1}'
\bigr)^{1/\rho} + 2\mu_+\bigl(\bigl[s,t'_j
\bigr] \times\bigl[t_j',t'_{j+1}
\bigr]\bigr). 
\end{eqnarray*}

Taking the supremum over $(t_i)$, then the $\rho$th power, summing
over $%
(t'_j)$ and finally taking the supremum over $(t_j')$
gives
\begin{eqnarray*}
V_{1, \rho}^+ (R;U_{[s,t]}) &\leq& C \bigl( \omega(s,t)+ \mu_+
\bigl(\bigl\{(u,v) \in[s,t]^2\mid u \leq v \bigr\}
\bigr)^{\rho} \bigr)^{1/\rho}
\\
&\leq& C \bigl( \omega(s,t)^{1/\rho} + \mu_{+}
\bigl([s,t]^2\bigr) \bigr),
\end{eqnarray*}
for some constant $C$ depending on $\rho$ only. Similarly,
\begin{eqnarray*}
V_{1, \rho}^+ (R;L_{[s,t]}) &\leq& C \bigl( \omega(s,t)^{1/\rho}
+ \mu_{+}\bigl([s,t]^2\bigr) \bigr).
\end{eqnarray*}
Now let $(t'_j)$ be a partition of $[s,t]$, fix $[t^{\prime
}_j,t'_{j+1}]$ and let $(t_i)$ be a partition of $[t^{\prime
}_j,t'_{j+1}]$. By (B.ii), $R(A_{i,j}) \geq0$ for all
$A_{i,j} =
[t'_j,t'_{j+1}]\times[t_i,t_{i+1}]$, thus
\begin{eqnarray*}
&& \sum_{t_i} \biggl\llvert R \pmatrix{
t_{i},t_{i+1}
\cr
t_{j}',t_{j+1}'}
\biggr\rrvert= \biggl\llvert R \pmatrix{ t_{j}',t_{j+1}'
\vspace*{2pt}\cr
t_{j}',t_{j+1}'} \biggr
\rrvert= \sigma^2\bigl(t'_j,t'_{j+1}
\bigr)
\end{eqnarray*}
and hence
\begin{eqnarray*}
&& V_{1, \rho}^+ (R;D_{[s,t]}) \leq\omega(s,t)^{1/\rho}.
\end{eqnarray*}
By Lemma~\ref{lemconc1} we conclude
\begin{eqnarray*}
V_{1, \rho}^{+}\bigl(R; [s,t]^2\bigr) &\leq& C
\bigl( V_{1, \rho}^+ (R;U_{[s,t]}) + V_{1, \rho}^+(R;D_{[s,t]})
+ V^+_{1, \rho}(R;L_{[s,t]}) \bigr)
\\
&\le& C \bigl( \omega(s,t)^{1/\rho} + \mu_{+}
\bigl([s,t]^2\bigr) \bigr)
\end{eqnarray*}
and \eqref{eqnon-diag-est} has been shown.\footnote{Note that in
fact we may deduce
the somewhat stronger conclusion
\[
V_{1, \rho}^{+}\bigl(R; [s,t]^2\bigr) \leq C
\bigl( \omega(s,t)^{1/\rho} + V_{1,\rho}^{+}
\bigl(R^+;[s,t]^2\bigr) \bigr)\qquad\forall[s,t]^2
\subseteq D_h.
\]\vspace*{-6pt}}

We now establish \eqref{eqnoff-diag-est-1rho}. Let $(s,t) \times(u,v)
\subseteq I^2\setminus D$, and let $(t_i)$ be a partition of $[s,t]$
and $%
(t'_j)$ be a partition of $[u,v]$. By nonnegativity of
nonoverlapping increments,
\begin{eqnarray*}
\sum_{t_i,t'_j} \biggl\llvert R \pmatrix{
t_{i},t_{i+1}
\cr
t_{j}',t_{j+1}'}
\biggr\rrvert&\le&\sum_{t_i,t'_j} \biggl\llvert R^-
\pmatrix{t_{i},t_{i+1}
\cr
t_{j}',t_{j+1}'}
\biggr\rrvert+ \biggl\llvert R^+ \pmatrix{ t_{i},t_{i+1}
\cr
t_{j}',t_{j+1}'}\biggr\rrvert
\\
&=& R^- \pmatrix{s,t
\cr
u,v} + R^+ \pmatrix{s,t
\cr
u,v}
\\
&\le&\bigg\llvert R \pmatrix{ s,t
\cr
u,v} \bigg\rrvert+ 2R^+ \pmatrix{s,t
\cr
u,v}.
\end{eqnarray*}
Taking the supremum over all partitions, the Cauchy--Schwarz inequality
\[
\bigg\llvert R \pmatrix{ s,t
\cr
u,v}\bigg\rrvert\le\bigg\llvert R\pmatrix{ s,t
\cr
s,t}
\bigg\rrvert^{1/2} \bigg\llvert R\pmatrix{ u,v
\cr
u,v}\bigg\rrvert
^{1/2}
\]
gives
\begin{eqnarray*}
V_{1,\rho}\bigl(R;[s,t] \times[u,v]\bigr) &\leq&\bigg\llvert R \pmatrix{ s,t
\cr
u,v} \bigg\rrvert+ 2 R^+ \pmatrix{ s,t
\cr
u,v}
\\
&\leq& C \bigl(\omega(s,t)^{1/(2\rho)} \omega(u,v)^{1/(2\rho)
}+
\mu_{+}\bigl([s,t]\times[u,v]\bigr) \bigr),
\end{eqnarray*}
and Lemma~\ref{lemconc2} completes the proof.
\end{pf*}

\section{Random Fourier series}\label{sec3}
\label{secrfs}

Let us now consider a (formal) random Fourier series
%
\begin{equation}
\label{eqnrfs} \Psi(t)=\frac{\a_{0}Y_{0}}{2}+\sum_{k=1}^\infty
\a_{k}Y^{k}\sin(kt)+\a%
_{-k}Y^{-k}
\cos(kt),
\end{equation}
where $Y^{k}$ are real-valued, centered random variables with $\E%
Y^{k}Y^{l}=\delta_{k,l}$ for all $k,l\in\Z$ and $\a_{k}$ are real-valued
coefficients. Since we are interested in properties of the covariance
of $%
\Psi$, we will formulate our conditions in terms of the squared
coefficients $a_k:= \a_k^2$, $k\in\Z$.

\begin{remark}
\label{remarkrFSexists} \label{remcontinuityrfs} Assume that $\a%
_k^2=O (\llvert k\rrvert ^{-(1+{1/\rho})} )$ for some $\rho> 0$.
Then %
\eqref{eqnrfs} converges uniformly almost surely, and the limit
yields a
continuous function. Moreover, if the $Y_k$ are Gaussian, $\Psi$ has
$\beta$%
-H\"older continuous trajectories\footnote{%
If $\beta= n + \tilde{\beta}$ for some $\tilde{\beta} \in(0,1]$, this
means that the trajectories are $n$-times differentiable and the $n$th
derivative is $\tilde{\beta}$-H\"older continuous.} almost surely for
all $%
\beta< \frac{1}{2\rho}$. This follows from \cite{K85}, Theorems 7.4.3~and~5.8.3.
\end{remark}

Our main theorem on random Fourier series follows:

%
\begin{theorem}
\label{teorfs} Consider the random Fourier series \eqref{eqnrfs}
with $%
(a_k)$ satisfying $\D^2(k^2 a_k)
\le0$ for all $k \in\Z$,
\[
\lim_{k \to\pm\infty} k^3\bigl\llvert\D^2
a_k\bigr\rrvert+ k^2\llvert\D a_k\rrvert=
0,
\]
$a_k = O (\llvert k\rrvert ^{-(1+{1/\rho})} )$ for some $\rho\ge
1$ for $%
k\to\pm\infty$ and $a_k$, $a_{-k}$ nonincreasing\footnote{The
monotonicity of $a_k$, $a_{-k}$ is required for the sole purpose of
using Lemma~\ref{lemmahoelder} below. In fact, it can be dropped when
we use Sobolev embeddings instead; cf. Remark~\ref{remsobolevemb}
below. However, we may only conclude finite $(1,\rho')$-variation for
any $\rho' > \rho$ in this case.} for $k\geq1$.
Then the covariance $R_\Psi$ of $\Psi$ has finite H\"older controlled
$%
(1,\rho)$-variation, and there is a constant $C>0$ such that
%
\begin{eqnarray}\label{eqn1rhovarboundcovrFS}
V_{1,\rho}\bigl(R_\Psi;[s,t]\times
[u,v]\bigr)
\le C \llvert t-s\rrvert^{1/(2\rho)}\llvert v-u\rrvert
^{1/(2\rho)}
\nonumber\\[-8pt]\\[-8pt]
\eqntext{\forall[s,t]\times[u,v] \subseteq[0,2
\pi]^2.}
\end{eqnarray}
The constant $C$ depends only on $\rho$ and $C_1$, where $C_1 \geq
\sup_{k
\in\Z} a_k \llvert k\rrvert ^{1+{1/\rho}}$.
\end{theorem}

Note that the model case $(a_k)=(\llvert k\rrvert ^{-2\a})$ for $\a\in
(\frac
{1}{2},1]$ is
contained as a special case in Theorem~\ref{teorfs}.

\begin{pf*}{Proof of Theorem \protect\ref{teorfs}}
Note first that, as already seen in Remark~\ref{remarkrFSexists},
$\Psi$~exists as a uniformly almost sure limit.\vspace*{1pt} Since $(a_k) \in l^1(\Z)$ we
have $%
(\a_k) \in l^2(\Z)$. Thus for fixed $t \in[0,2 \pi]$, $\Psi$
exists also as
a convergent sum in $L^2(\Omega)$. Set $Q_1 = [0,\pi]^2$, $Q_2 =
[0,\pi]\times[\pi,2\pi]$, $Q_3 = [\pi,2\pi]^2$
and $Q_4 = [\pi,2\pi]\times[0,\pi]$. We first show that %
\eqref{eqn1rhovarboundcovrFS} holds provided $[s,t]\times[u,v]
\subseteq Q_i$ for some $i = 1,\ldots,4$. Recall from \eqref
{eqnfourierseriesdecomp} that we can decompose the covariance as
\begin{eqnarray}\label{eqnrfsRdecomp}
R_{\Psi}(s,t) &=& K\bigl(\llvert t-s\rrvert\bigr) + K\bigl(\llvert t+s
\rrvert\bigr) + \tilde{K}\bigl(\llvert t-s\rrvert\bigr) + \tilde
{K}\bigl(
\llvert t+s\rrvert\bigr)
\nonumber\\[-8pt]\\[-8pt]\nonumber
&=:& R^{-} ( s,t ) + R^+ ( s,t ) + \td
R^{-} (
s,t ) - \td R^+ ( s,t ),
\nonumber
\end{eqnarray}
where
\begin{eqnarray*}
&& K(t) = \sum_{k = 1}^{\infty}
\alpha_k^2 \cos(kt) \quad\mbox{and}\quad\tilde{K}(t) =
\sum_{k = 1}^{\infty} \alpha_{-k}^2
\cos(kt).
\end{eqnarray*}

Using the triangle inequality it is enough to show the estimate %
\eqref{eqn1rhovarboundcovrFS} for $R^\pm,\td R^\pm$
separately. From Lemma~\ref{lemmarhobound} below we know that $K$
and $\td K$ are convex on $%
[0,2\pi]$ and nonincreasing on $[0,\pi]$. By Lemma~\ref{lemmahoelder}
below, $K$ and $\td K$ are $\frac{1}{\rho}$-H\"older continuous.
Convexity implies that
\[
\partial_{s,t} R^- = -K'' \leq0.
\]
Therefore $\mu:= -\mu_{-}:= \partial_{s,t}^2 R^-$ yields a Radon
measure on $(0,1)^2 \setminus D$ and condition (B.i) of Theorem~\ref
{teomain} is satisfied. Condition (B.ii) holds for $h = \pi$ since
$K$ is nonincreasing. $(JM)_{\rho,\omega}$ follows from H\"
older-continuity with $\omega(s,t) = C\llvert t-s\rrvert $. Since
$R^-$ is a
covariance function, it satisfies the Cauchy--Schwarz inequality. Thus
we may apply part \textup{B} in Theorem~\ref{teomain} to conclude that there
is a constant $C$ such that
\[
V_{1,\rho}\bigl(R^-;[s,t] \times[u,v]\bigr) \leq C \llvert t-s\rrvert
^{1/(2\rho)} \llvert v-u\rrvert^{1/(2\rho)}
\]
holds for all $[s,t] \times[u,v] \in Q_1$. The same reasoning works
for $\tilde{R}^-$. Using again convexity of $K$, we have $\partial
^2_{s,t} R^+ = K'' \geq0$ which shows that $\nu:= \nu_+:= \partial
^2_{s,t} R^+$ is a Radon measure on $(0,T) \setminus D$. Hence
condition (A.i) of Theorem~\ref{teomain} holds for~$R^+$. In \eqref
{eqncondAiiforrfs} we have seen that also (A.ii) is satisfied for
$R^+$ on $Q_1$, and we may conclude, using part \textup{A} of Theorem~\ref
{teomain}, that
\[
V_{1,\rho}\bigl(R^+;[s,t] \times[u,v]\bigr) \leq V_{1}
\bigl(R^+;[s,t] \times[u,v]\bigr) \leq R^+ \pmatrix{ s,t
\cr
u,v}
\]
holds for all $[s,t] \times[u,v] \in Q_1$. $R^+$ will in general not
be a covariance function, but we may use the $2\pi$-periodicity of $K$
to deduce the Cauchy--Schwarz inequality for $R^+$ as well. Indeed,
\[
R^+(s,t) = K(t+s) = K\bigl(t - (2\pi- s)\bigr) = R^-(2\pi- s,t),
\]
and using the Cauchy--Schwarz inequality for $R^-$ implies that
%
%
\begin{eqnarray*}
R^+ \pmatrix{ s,t
\cr
u,v} &\leq&\sqrt{R^+\pmatrix{ s,t
\cr
s,t}} \sqrt{R^+
\pmatrix{ u,v
\cr
u,v} }
\\
&\leq&\llVert K\rrVert_{1/\rho\mbox{-H\"ol}}\llvert t-s\rrvert
^{1/(2\rho)} \llvert u-v\rrvert^{1/(2\rho)},
\end{eqnarray*}
where the second estimate follows from H\"older continuity of $K$ as
seen in \eqref{eqnpospartrfshoelder}. The same is true for $\tilde
{%
R}^+$ which shows \eqref{eqn1rhovarboundcovrFS} for $R^+,\td
R^+$ and $%
[s,t]\times[u,v] \subseteq Q_1$. The process $t\mapsto\Psi_{t+\pi}$
has the
same covariance as $\Psi$. Thus estimate %
\eqref{eqn1rhovarboundcovrFS} also holds for $[s,t]\times[u,v]
\subseteq Q_3$. By symmetry considerations, if $E$ is any rectangle in $Q_2$
or $Q_4$, there is a rectangle $\bar{E}$ in $Q_1$ (or in $Q_3$)\vspace*{2pt} with the
same side length such that $R^+(E) = R^-(\bar{E})$, $\tilde{R}^+(E) =
\tilde{%
R}^-(\bar{E})$ and vice versa for $R^-$, $\tilde{R}^-$. Thus~\eqref{eqn1rhovarboundcovrFS} also holds for $[s,t]\times[u,v]
\subset Q_i$ for $i=2,4$. The general case just follows by subdividing a
given rectangle $[s,t]\times[u,v]$ in at maximum four rectangles lying
in $%
Q_1,\ldots,Q_4$ and using the estimates above (which only leads to a larger
constant). This proves the theorem.
\end{pf*}

\subsection{Convexity, monotonicity and H\"older regularity of cosine series}\label{sec31}

We start by deriving conditions for convexity and monotonicity of cosine
series
%
\begin{equation}
\label{eqncos-series} K(t)=\frac{a_{0}}{2}+\sum_{k=1}^{\infty}a_{k}
\cos(kt).
\end{equation}
In the following let $\D$, $\D^2$ be the first and second forward-difference
operators, that is, for a sequence $\{a_k\}_{k\in\N}$
\[
\D a_k:= a_{k+1}-a_k
\]
and $\D^2:= \D\circ\D$. Moreover, let
\[
D_n(t):= 1 + 2\sum_{k=1}^n
\cos(kt), \qquad t \in\R,
\]
be the Dirichlet kernel and
\[
F_n(t):= \sum_{k=0}^{n}D_k(t),\qquad t \in\R,
\]
be the unnormalized Fej\'er kernel.

\begin{lemma}
\label{lemmarhobound} Let $\{a_k\}_{k \in\N}$ be such that
%
%
\begin{equation}
\Delta^2\bigl(k^2a_k\bigr)\le0,\qquad k\in
\mathbb{N} \label{E1}
\end{equation}
and
%
\begin{equation}
\label{eqnconvass} \lim_{k \to\infty} k^3\bigl\llvert
\D^2 a_k\bigr\rrvert+ k^2\llvert\D
a_k\rrvert+ k\llvert a_k\rrvert= 0.
\end{equation}
Then the cosine series \eqref{eqncos-series}
exists locally uniformly in $(0,2\pi)$, is convex on $[0,2\pi]$ and
decreasing on $[0,\pi]$.
\end{lemma}

\begin{pf}
The proof follows ideas from \cite{K11}; we include it for the reader's
convenience. We first note that since
\[
\D\bigl(k^2 a_k\bigr) = k^2 \D
a_k + (2k+1)a_{k+1}
\]
and
\[
\D^2\bigl(k^2a_k \bigr) = k^2
\D^2 a_k + 2(2k+1) \D a_{k+1} + 2
a_{k+2},
\]
assumption \eqref{eqnconvass} is equivalent to
%
\begin{equation}
\label{eqnaltconvexass} \lim_{k \to\infty} \bigl\llvert k\D^2
\bigl(k^2 a_k\bigr)\bigr\rrvert+ \bigl\llvert\D
\bigl(k^2 a_k\bigr)\bigr\rrvert+ k\llvert a_k
\rrvert= 0.
\end{equation}
Using the Abel transformation, we observe
\begin{eqnarray*}
S_n(t) = \frac{a_0}{2} + \sum_{k=1}^n
a_k \cos(kt) &=& \frac{1}{2}%
\sum
_{k=0}^n \D a_{k+1} D_k(t) +
\frac{1}{2} a_{n+1} D_n(t).
\end{eqnarray*}
By the assumptions and \eqref{eqnaltconvexass} we have $\sum
_{k=1}^\infty
\llvert \D a_k\rrvert < \infty$. Since $\sup_{n \in\N} D_n(t)$ is
bounded locally
uniformly on $(0,2\pi)$ and $a_n \to0$, we observe that
\[
K(t):= \frac{a_0}{2} + \sum_{k=1}^\infty
a_k \cos(kt) = \frac{1}{2} 
\sum
_{k=0}^\infty\D a_k D_k(t)
\]
exists locally uniformly and is continuous in $(0,2\pi)$.

The Ces\`aro means of the sequence $S_n(t)$ are given by
\[
\s_n(t) = \frac{a_0}{2} + \sum_{k=1}^n
\biggl(1-\frac{k}{n+1} \biggr) a_k \cos(kt).
\]
By Fej\'er's theorem (\cite{Z59}, Theorem III.3.4) and continuity of
$K$, $\s%
_n \to K$ locally uniformly in $(0,2\pi)$. Hence, $\s_n^{\prime
\prime}\to
K^{\prime\prime}$ in the space of distributions on $(0,2\pi)$. Clearly,
\[
\s_n^{\prime\prime}(t) = -\sum_{k=0}^n
\biggl(1-\frac{k}{n+1} \biggr) k^2a_k \cos(kt).
\]
Let $\b_k:= (1-\frac{k}{n+1} ) k^2a_k$. Using summation
by parts
twice we obtain
\begin{eqnarray*}
2\s_n^{\prime\prime}(t) &=& \sum_{k=0}^{n}
\D\b_k D_k(t)
\\
&=& \D\b_n F_n(t) - \sum_{k=0}^{n-1}
\D^2\b_k F_k(t)
\\
&=& - \sum_{k=0}^{n-1} \D^2
\bigl(k^2a_k\bigr) F_k(t)
\\
&&{}- \sum
_{k=0}^{n-1} \biggl(\frac{k \D%
^2(k^2a_k)}{n+1} -
\frac{2\D((k+1)^2 a_{k+1}
)}{n+1} \biggr) F_k(t) + \frac{n^2}{n+1}a_n
F_n(t).
\end{eqnarray*}


We have $0\le F_n(t)\le\frac{C}{t^2}+\frac{C}{(2\pi-t)^2}$, where
$C> 0$ is an absolute constant. Therefore, for every
$\varepsilon$ with $0<\varepsilon<2\pi$,
%
\begin{equation}
\sup_{n\ge0;t\in[\varepsilon,2\pi-\varepsilon
]}F_n(t)=C_{\varepsilon}<\infty.
\label{E2}
\end{equation}

It follows from (\ref{E1}) that for all $t\in[0,2\pi]$ and $n\ge1$,
\[
-\sum_{k=0}^{n-1}\Delta^2
\bigl(k^2a_k\bigr)F_k(t)-\frac{1}{n+1}
\sum_{k=0}^{n-1}k\Delta^2
\bigl(k^2a_k\bigr)F_k(t)\ge0. %
\]
Moreover, since $k\llvert a_k\rrvert \rightarrow0$ as $k\rightarrow
\infty$ [see
(\ref{eqnconvass})], and (\ref{E2}) holds, we have
\[
\sup_{t\in[\varepsilon,2\pi-\varepsilon]}\frac
{n^2}{n+1}\llvert a_n\rrvert
F_n(t)\rightarrow0 %
\]
as $n\rightarrow\infty$.
Finally, set
%
\begin{equation}
S_n(t)=\frac{2}{n+1}\sum_{k=0}^{n-1}
\Delta\bigl((k+1)^2a_{k+1}\bigr)F_k(t).
\label{E3}
\end{equation}
It is easy to see, using (\ref{E2}), that for all $n\ge1$,
\[
\sup_{t\in[\varepsilon,2\pi-\varepsilon]}\bigl\llvert S_n(t)\bigr
\rrvert\le
\frac
{2C_{\varepsilon}}{n+1}\sum_{k=0}^{n-1}\bigl
\llvert\Delta\bigl((k+1)^2a_{k+1}\bigr)\bigr\rrvert.
\]
Next, taking into account (\ref{eqnconvass}) and the Ces\`{a}ro summability theorem
for convergent sequences, we obtain
\[
\sup_{t\in[\varepsilon,2\pi-\varepsilon]}\bigl\llvert S_n(t)\bigr\rrvert
\rightarrow0 %
\]
as $n\rightarrow\infty$, for all $t\in[\varepsilon,2\pi
-\varepsilon]$. Summarizing what was said above, we see that for every
$0<\varepsilon< 2\pi$,
\[
\liminf_{n\rightarrow\infty}\inf_{t\in[\varepsilon,2\pi
-\varepsilon]}
\sigma^{\prime\prime}_n(t)\ge0. %
\]

For any nonnegative test-function $\vp\in C_c^\infty(0,2\pi)$, Fatou's
lemma implies
\begin{eqnarray*}
K^{\prime\prime}(\vp) = \lim_{n\to\infty} \s_n^{\prime\prime
}(
\vp) \ge\int_0^{2\pi} \liminf
_{n\to\infty} \s_n^{\prime\prime}(t)\vp(t)\,dt \ge0;
\end{eqnarray*}
that is, $K^{\prime\prime}$ is a nonnegative distribution on
$(0,2\pi)$.
Thus $K$ is convex on $[0,2\pi]$.

Assume now that $K$ is not decreasing on $[0,\pi]$; that is, there are
$s<t \in
[0,\pi]$ such that $K(s) < K(t)$. Since $K$ is given as a cosine
series, we
have $K(s) = K(s')$ and $K(t) = K(t')$ for $s'=
2\pi- s$ and $t'= 2\pi- t$. Choose $\l\in(0,1)$ such that
$\l s
+ (1-\l)s'= t$. Then
\[
K\bigl(\l s + (1-\l)s'\bigr) = K(t) > K(s) = \l K(s) + (1-\l) K
\bigl(s'\bigr)
\]
which is a contradiction to the convexity of $K$.
\end{pf}

Concerning H\"older regularity of cosine series we recall the following:

\begin{lemma}[(\protect\cite{L48}, Satz 8)]
\label{lemmahoelder} A cosine\vspace*{2pt} series \eqref{eqncos-series} with
nonincreasing coefficients $a_k \downarrow0$ for $k \to\infty$ is $%
\frac{1}{\rho}$-H\"older continuous if and only if $a_k = O
(k^{-(1+(1/\rho))} )$ for $k \to\infty$.
\end{lemma}

\begin{remark}\label{remsobolevemb}
The above lemma gives a sharper result than what is obtained by usual
Sobolev embeddings. Indeed: recall that an $L^2$ function on the torus
with Fourier coefficients $(a_k)$ is in the Sobolev space $H^s$ if and
only if \mbox{$((1+\llvert k\rrvert ^s) a_k) \in l^2$}. By a standard
Sobolev embedding (here in dimension 1), such functions are
$(s-1/2)$-H\"older, provided $s>1/2$. Hence, a cosine series \eqref
{eqncos-series} with coefficients $a_k = O (k^{-(1+(1/\rho))} )$ for $k \to\infty$ is $\a$-H\"older for all $\a< 1/\rho$.
\end{remark}

\subsection{Stability under approximation}\label{sec32}

We now aim to prove stability of the estimates provided in Theorem~\ref{teorfs} under approximations of $\Psi$. These stability properties
will be
used in Section~\ref{secrandomFseriesasRP} to prove the
convergence (in
rough path topology) of Galerkin and hyper-viscosity approximations of
random Fourier series. Let us consider
%
\begin{equation}
\label{eqnrfs2} \td\Psi(t)=\frac{\a_{0}\b_0 Y^{0}}{2}+\sum_{k=1}^\infty
\a_{k}\b_k Y^{k}\sin(kt)+\a_{-k}
\b_{-k} Y^{-k}\cos(kt),
\end{equation}
with $Y^k$ as above and $(\a_k),(\b_k)$ real-valued sequences. In the
applications, the multiplication of the coefficients by $\b_k$ will
correspond to a smoothing of $\Psi$. We thus aim to prove that the estimates
given in Theorem~\ref{teorfs} remain true uniformly for $(b_k)=(\b
_k^2)$ in
an appropriate class of sequences. This will naturally lead to the following:

\begin{definition}\label{defnegli}
(1)~A sequence $(b_k)_{k\in\Z}$ is negligible if there are finite, signed,
real Borel measures $\mu_1,\mu_2$ on $\mcS^1:=\R/2\pi\Z$ such that
\[
b_k = \int_{0}^{2\pi} \cos(kr)
\mu_1(dr),\qquad b_{-k} = \int_{0}^{2\pi}
\cos(kr) \mu_2(dr)\qquad\forall k\in\N.
\]

(2)~A family of sequences $(b^\tau_k)$ is uniformly negligible if
each $%
(b^\tau_k)$ is negligible with associated measures $\mu_1^\tau,\mu
_2^\tau$
being uniformly bounded in total variation norm.

(3)~For two bounded sequences $(a_k)$, $(c_k)$ we write $(a_k) \preceq
(c_k) $ if there is a negligible sequence $(b_k)$ such that $a_k = c_k b_k$
for every $k\in\Z$.
\end{definition}

%
\begin{example}
Some (simple) examples of
negligible sequences are:
\begin{longlist}[(2)]
\item[(1)] $(b_k)\equiv C$, with $\mu_1=\mu_2=C\delta_0$,

\item[(2)] $(b_k) \in l^1(\Z)$, with $\mu_1 = \sum_{k=1}^\infty b_k \cos(kt)\,dt$
and $\mu_2 = \sum_{k=1}^\infty b_{-k} \cos(kt)\,dt$.
\end{longlist}
\end{example}

In the forthcoming Lemmas~\ref{lemmakolm} and~\ref
{lemmauniformseriesbound}, we will give sufficient conditions for
(uniform) negligibility.

As will be seen below, our results are uniform relative to ``negligible''
perturbations as in \eqref{eqnrfs2}.


\begin{proposition}
\label{propmollRFS} Consider the random Fourier series \eqref{eqnrfs2}
with $(a_k)$ satisfying the assumptions of Theorem~\ref{teorfs}. Let $(b_k)$
be negligible. Then
\[
V_{1,\rho}\bigl(R_{\td\Psi};[s,t]^2\bigr) \le C \llvert
t-s\rrvert^{1/\rho}\qquad\forall[s,t]^2
\subseteq[0,2\pi]^2.
\]
The constant $C$ depends only on $\rho$, the constant $C_1=\sup_{k
\in\Z}
a_k \llvert k\rrvert ^{1+{1/\rho}}$ and a constant $C_2$ which
bounds $%
\llVert \mu_1\rrVert _{\mathrm{TV}} $ and $\llVert \mu_2\rrVert _{\mathrm{TV}}$ with
$\mu_1,\mu_2$
corresponding to $%
(b_k)$; cf. Definition~\ref{defnegli}.
\end{proposition}

This proposition is a special case of Proposition~\ref{proppairRFS} below.
Consider another random Fourier series
\[
\Phi(t)=\frac{\gamma_{0} Z^{0}}{2}+\sum_{k=1}^\infty
\gamma_{k} Z^{k}\sin(kt)+\gamma_{-k}
Z^{-k}\cos(kt),
\]
and assume that the $Z^k$ fulfill the same conditions as the $Y^k$.
Furthermore, assume that $\{Y^k,Z^k\}_{k\in\Z}$ are uncorrelated random
variables, and set $c_k:= \gamma_k^2$, $\varrho_k:= \E Y^k Z^k$ and
\[
R_{\Psi,\Phi}(s,t):= \E\Psi(s) \Phi(t).
\]

Then the following holds:

\begin{proposition}
\label{proppairRFS} Assume that there is a sequence $(d_k)$
satisfying the
assumptions of Theorem~\ref{teorfs} such that
\[
(b_k):= \biggl(\frac{\alpha_k \gamma_k \varrho_k}{d_k} \biggr)
\]
is negligible with associated measures $\mu_1$, $\mu_2$. Then
\[
V_\rho\bigl(R_{\Psi,\Phi};[s,t]^2\bigr) \le
V_{1,\rho}\bigl(R_{\Psi,\Phi
};[s,t]^2\bigr) \le C \llvert
t-s\rrvert^{1/\rho}\qquad\forall[s,t]^2
\subseteq[0,2\pi]^2.
\]
The constant $C$ depends only on $\rho$, the constant $C_1=\sup_{k
\in\Z}d_k \llvert k\rrvert ^{1+{1/\rho}}$ and a constant $C_2$ which
bounds $%
\llVert \mu_1\rrVert _{\mathrm{TV}} $ and $\llVert \mu_2\rrVert _{\mathrm{TV}}$.
\end{proposition}

\begin{pf}
Arguing as for Theorem~\ref{teorfs} we observe
\begin{eqnarray*}
&& V_{1,\rho}\bigl(R_{\Psi,\Phi};[s,t]\times[u,v]\bigr)
\\
&&\qquad  \lesssim
V_{1,\rho
}\bigl(R^{-};[s,t]%
\times[u,v]\bigr) +
V_{1,\rho}\bigl(R^{+};[s,t]\times[u,v]\bigr)
\\
&&\quad\qquad{} + V_{1,\rho}\bigl(\tilde{R}^{-};[s,t]\times[u,v]\bigr) +
V_{1,\rho}\bigl(\tilde{R}%
^{+};[s,t]\times[u,v]
\bigr),
\end{eqnarray*}
with $R^-(s,t)=K(t-s)$, $R^+(s,t)=K(t+s)$, $\tilde{R}^-(s,t)=\tilde{K}(t-s)$
and $\tilde{R}^+(s,t)=\tilde{K}(t+s)$
\begin{eqnarray*}
K(t) &:=& \frac{1}{2}\sum_{k=1}^{\infty}d_{-k} b_{-k} \cos(kt),
\\
\tilde{K}(t)&:=&
\frac{1}{2}\sum_{k=1}^{\infty}
d_{k} b_k \cos(kt).
\end{eqnarray*}
We thus need to estimate the mixed $(1,\rho)$-variation of cosine series
under multiplication with negligible sequences. In the following we
consider $R^\pm$, $\tilde R^\pm$ can be treated analogously. Let
\[
R_0^{\pm} (t,s):= \frac{d_0}{2}+\sum
_{k=1}^\infty d_{-k} \cos\bigl(k(t \pm s)
\bigr).
\]
We then apply Proposition~\ref{proprhovariationconv} below with
$R^\pm_{\mu}=R^\pm$, $R^\pm=R_0^\pm$, $a_k=d_{-k}$, $b_k=b_{-k}$
to obtain
\begin{eqnarray*}
V_{1,\rho}\bigl(R^{\pm};[s,t]^2\bigr) &\le&\llVert
\mu\rrVert_{\mathrm{TV}} \sup_{0\leq z
\leq2\pi} V_{1,\rho}
\bigl(R_0^{\pm}; [s-z,t-z]\times[s,t]\bigr)
\end{eqnarray*}
for every $[s,t] \subseteq[0,2\pi]$. By Theorem~\ref{teorfs}
applied to $R_0^{\pm}$, we have
\[
\sup_{0\leq z \leq2\pi} V_{1,\rho}\bigl(R_0^{\pm}; [s-z,t-z]\times
[s,t]\bigr) \le C \llvert t-s\rrvert^{1/\rho},
\]
which completes the proof.
\end{pf}

In the following let $\mcM(\mcS^1)$ be the space of signed, real
Borel-measures on the circle $\mcS^1$ with finite total variation $%
\llVert \cdot\rrVert _{\mathrm{TV}}$. Define $\mcM^w(\mcS^1)$ to be $\mcM(\mcS^1)$
endowed with
the topology of weak convergence. For $B \in L^1(\mcS^1)$ we set $\mu
_B:=
B \,dt \in\mcM(\mcS^1)$ to be the associated measure with density $B$.

\begin{lemma}
\label{lemmarhovarunifbddunderconv} Let $\mu\in\mcM(\mcS
^1)$, $%
R\colon\mcS^1 \times I \to\R$ and set $R_{\mu}(s,t):= (R(\cdot,t) \ast
\mu)(s)$. Then
\[
V_{\gamma,\rho}\bigl(R_{\mu};[s,t]\times[u,v]\bigr) \leq\llVert\mu
\rrVert_{\mathrm{TV}} \sup_{x \in\mcS^1} V_{\gamma,\rho}\bigl
(R;[s-x,t-x]\times[u,v]\bigr)
\]
for all $[s,t] \times[u,v] \subseteq\mcS^1 \times I$ and $1 \leq
\gamma
\leq\rho$.
\end{lemma}

\begin{pf}
Let $(t_i)$, $(t'_j)$ be partitions of $[s,t]$, respectively,
$[u,v]$. From
Jensen's inequality,
\begin{eqnarray*}
\biggl\llvert R_{\mu} \pmatrix{ t_{i},t_{i+1}
\cr
t_{j}',t_{j+1}'} \biggr\rrvert
^{\gamma} & \le&\biggl( \int_{\mcS^1} \biggl\llvert R
\pmatrix{ t_i-x,t_{i+1}-x
\cr
t_{j}',t_{j+1}'}
\biggr\rrvert \,d\llvert\mu\rrvert(x) \biggr)^\g
\\
& \le&\llVert\mu\rrVert_{\mathrm{TV}}^\g\int
_{\mcS^1} \biggl\llvert R \pmatrix{ t_i-x,t_{i+1}-x
\cr
t_{j}',t_{j+1}'}\biggr
\rrvert^\g \,d\frac{\llvert \mu\rrvert (x)}{\llVert \mu\rrVert _{\mathrm{TV}}}.
\end{eqnarray*}
Summing over $t_i$ and using again Jensen's inequality for $\frac{\rho
}{\g}$
yields
\begin{eqnarray*}
&& \sum_{t'_j} \biggl( \sum
_{t_i} \biggl\llvert R_{\mu} \pmatrix{
t_{i},t_{i+1}
\cr
t_{j}',t_{j+1}'}
\biggr\rrvert^{\gamma} \biggr)^{\rho/\gamma}
\\
&&\qquad \le \llVert\mu
\rrVert
_{\mathrm{TV}}^\rho\int_{\mcS^1} \sum
_{t'_j} \biggl( \sum_{t_i} \biggl
\llvert R \pmatrix{ t_i-x,t_{i+1}-x
\cr
t_{j}',t_{j+1}'}
\biggr\rrvert^\g\biggr)^{\rho/\g} \,d
\frac{\llvert \mu\rrvert (x)}{\llVert
\mu\rrVert _{\mathrm{TV}}}
\\
&&\qquad \leq \llVert\mu\rrVert_{\mathrm{TV}}^\rho\int
_{\mcS^1} V_{\g,\rho}^\rho\bigl(R;[s-x,t-x]
\times[%
u,v]\bigr) \,d\frac{\llvert \mu\rrvert (x)}{\llVert \mu\rrVert _{\mathrm{TV}}}
\\
&&\qquad \leq\llVert\mu\rrVert_{\mathrm{TV}}^\rho\sup
_{x \in\mcS^1} V_{\g,\rho}^\rho\bigl(R;[s-x,t-x]%
\times[u,v]\bigr).
\end{eqnarray*}
Taking the supremum over all partitions yields the inequality.
\end{pf}

\begin{remark}
\label{remyoungineq} In many cases, $x \mapsto V_{\gamma,\rho}(R;
[s-x,t-x]\times[s,t])$ attains its maximum at $x=0$. In this case our
inequality above reads
\[
V_{\gamma,\rho}\bigl(R \ast\mu;[s,t]^2\bigr) \leq\llVert\mu
\rrVert_{\mathrm{TV}} V_{\gamma,\rho}\bigl(R;[s,t]^2\bigr)
\]
for all squares $[s,t]^2 \subseteq[0,2\pi]^2$. Lemma~\ref{lemmarhovarunifbddunderconv} can thus be interpreted as a
Young-inequality for the mixed $(\gamma,\rho)$-variation of a
function with
two arguments. If $\mu= \delta_0$, we have $b_k = 1 $ for every $k$
and the
estimate is thus sharp.
\end{remark}

\begin{proposition}
\label{proprhovariationconv} Let $R^+_{\mu}, R^-_{\mu}\dvtx [0,2\pi
]^2 \to
\R
$ be continuous functions of the form
\begin{eqnarray*}
\label{eqnproductfourier} R_{\mu}^{\pm} (s,t) &=& \frac{a_0b_0}{2}+\sum
_{k=1}^\infty a_k b_k
\cos\bigl(k(s \pm t)\bigr)
\end{eqnarray*}
with $a_k,b_k$ being real-valued\vspace*{1.5pt} coefficients such that $\sum
_{k=1}^\infty
\llvert a_k\rrvert < \infty$, and assume that there is a measure $\mu
\in\mcM
(\mcS^1)$
such that
\begin{eqnarray*}
&&b_k = \int_{0}^{2\pi} \cos(kr)
\mu(dr).
\end{eqnarray*}
Set
\begin{eqnarray*}
&&R^{\pm} (t,s) = \frac{a_0}{2}+\sum_{k=1}^\infty
a_k \cos\bigl(k(t \pm s)\bigr).
\end{eqnarray*}
Then for every $1\leq\gamma\leq\rho$,
\begin{eqnarray*}
V_{\g,\rho}\bigl(R^{\pm}_{\mu};[s,t]\times[u,v]\bigr) & \le&\llVert\mu\rrVert_{\mathrm{TV}} \sup_{0\leq z
\leq2\pi}
V_{\g,\rho}\bigl(R^{\pm}; [s-z,t-z]\times[u,v]\bigr)
\end{eqnarray*}
for every $[s,t]\times[u,v] \subseteq[0,2\pi]^2$.
\end{proposition}

\begin{pf}
Let $a_{-k}:=a_k$, $b_{-k}:=b_k$ for $k\in\N$. Since $\sum_{k=1}^\infty
\llvert a_k\rrvert < \infty$, we observe
\begin{eqnarray*}
R_{\mu}^{\pm}(s,t) =\frac{1}{2} \sum
_{k \in\Z} a_k b_k e^{ik(t
\pm s)} =
\bigl(R^{\pm} (\cdot,t) \ast\mu\bigr) (s)
\end{eqnarray*}
and the estimate is thus a direct consequence from Lemma~\ref{lemmarhovarunifbddunderconv}.
\end{pf}

\subsection{(Uniform) negligibility}\label{sec33}\label{secnegligibility}

In order to use Proposition~\ref{proprhovariationconv} to control
the $%
(\gamma,\rho)$-variation of $R(s,t)$, we need to control $\llVert \mu
\rrVert _{\mathrm{TV}}$.
We recall the following:

\begin{lemma}
\label{lemmakolm} Let $\{b_k\}_{k \in\N}$ be a sequence satisfying $b_k
\to b \in\R$ for $k \to\infty$, and let $S_n(t):= \frac{b_0}{2} +
\sum_{k=1}^n b_k \cos(kt) $. Assume one of the following conditions:
\begin{longlist}[(2)]
\item[(1)] $\sum_{k=1}^\infty\llvert b_k-b\rrvert < \infty$;

\item[(2)] there exists a nonincreasing sequence $A_k$ such that $%
\sum_{k=0}^\infty A_k < \infty$ and $\llvert \D b_k\rrvert \le A_k$
for all $k \ge0$;

\item[(3)] $b_k$ is quasi-convex, that is,
\[
\sum_{k=0}^\infty(k+1) \bigl\llvert
\D^2 b_{k}\bigr\rrvert< \infty.
\]
\end{longlist}

Then, $B(t) = \frac{b_0}{2} + \sum_{k=1}^\infty b_k \cos(kt) $ exists
locally uniformly on $(0,2\pi)$, and the right-hand side is the Fourier
series of $B$. Moreover,
\[
\mu_{S_n} \rightharpoonup\mu_B + b \delta_0
=: \mu\qquad\mbox{weakly in } \mcM\bigl(\mcS^1\bigr)
\]
and $b_k = \int_{0}^{2\pi} \cos(kr) \mu(dr)$. Moreover, there is a
numerical constant $C > 0$ such that
%
\begin{equation}
\label{eqnL1-bound} \llVert\mu\rrVert_{\mathrm{TV}} \le\llvert b\rrvert+ C
\cases{
\displaystyle\sum_{k=0}^\infty\llvert b_k - b\rrvert, &\quad in case (1),
\vspace*{3pt}\cr
\displaystyle\sum_{k=0}^\infty A_k, &\quad in case (2),
\vspace*{3pt}\cr
\displaystyle\sum_{k=0}^\infty(k+1) \bigl\llvert\D^2 b_{k}\bigr\rrvert, &\quad in case (3).}
\end{equation}
\end{lemma}

\begin{pf}
The case $b = 0$ is classical [(1) is trivial; cf. \cite{T73} for (2)
and \cite{K23} for~(3)].
The case $b\ne0$ may be reduced to $b=0$ by noting that $bD_n(t) \to
2\pi
b\delta_0$ in $\mcM^w(\mcS^1)$, where $D_n$ is the Dirichlet kernel.
\end{pf}

Lemma~\ref{lemmakolm} in combination with Proposition~\ref{proprhovariationconv} allows us to derive bounds on the $\rho
$-variation of
covariance functions of the type discussed here, depending on $\mu$
only via its
total variation norm. Since we will use this to prove uniform
estimates, we
will need the following uniform estimates on the $L^1$-norm of cosine series.

\begin{lemma}
\label{lemmauniformseriesbound} Let $b \in C^1(0,\infty)$ with
$b(r) \to
0 $ for $r \to\infty$ and $b^\tau_k:= b(\tau^m k)$ for some $\tau, m > 0$.
If:
\begin{longlist}[(2)]
\item[(1)] $b$ is convex, nonincreasing, then $b^\tau_k$ satisfies the
assumptions of Lemma \ref{lemmakolm};

\item[(2)] $B^\tau(t) = \frac
{b^\tau_0}{2}
+ \sum_{k=1}^\infty b^\tau_k \cos(kt)$ exists locally uniformly in
$(0,2\pi)$
and
\[
\bigl\llVert B^\tau\bigr\rrVert_{L^1([0,2\pi])} \le C
b_0,
\]
for some $C > 0$;

\item[(3)] $b \in C^2(0,\infty)$ with $r\mapsto r\llvert b^{\prime\prime
}(r)\rrvert $ being
integrable, then $b^\tau_k$ satisfies the assumptions of Lemma~\ref{lemmakolm}, (3)~and
\begin{eqnarray*}
&&\bigl\llVert B^\tau\bigr\rrVert_{L^1([0,2\pi])} \le C \int
_{0}^{\infty} r \bigl\llvert b^{\prime
\prime}(r)\bigr
\rrvert \,dr,
\end{eqnarray*}
for some $C > 0$ with $B^\tau$ as in (1).
\end{longlist}
\end{lemma}

\begin{pf}
(1)~Since $b$, $\llvert b'\rrvert $ are nonincreasing $\D b^\tau_k
\le0$
and $-%
\D b_k$ is nonincreasing. We set $A_k:= -\D b_k$. Clearly, $%
\sum_{k=0}^\infty A_k = 2b_0$, and the claim follows from Lemma~\ref{lemmakolm}.

(2)~Let $b^\tau(r):= b(\tau^m r)$, and observe
\[
\D^2 b^\tau_k = \int_{k+1}^{k+2}
\int_{s-1}^{s} \bigl(b^\tau
\bigr)^{\prime
\prime}(r) \,dr\,ds.
\]
Since $ (b^\tau)^{\prime\prime}\,dr=\tau^mb^{\prime
\prime}(\tau^m
r) \,d(\tau^m r)$, elementary calculations show
\begin{eqnarray*}
\sum_{k=0}^\infty(k+1) \bigl\llvert
\D^2 b^\tau_{k}\bigr\rrvert&\le& 2 \int
_{0}^{\infty} r \bigl\llvert b^{\prime\prime}(r)\bigr
\rrvert \,dr,
\end{eqnarray*}
and Lemma~\ref{lemmakolm} completes the proof.
\end{pf}

%
\begin{example}
As an application of Lemma~\ref{lemmakolm}, we see that the sequence
$(b_k) = (k^{-\alpha})$, $\alpha> 0$, is negligible. Furthermore, the
sequence $(b_k) = (e^{- \tau k^{\alpha}})$, $\alpha, \tau> 0$, is
uniformly negligible in $\tau$ which follows from Lemma~\ref
{lemmauniformseriesbound}.
\end{example}

\subsection{Random Fourier series as rough paths}\label{sec34}
\label{secrandomFseriesasRP}

We now return to the initial problem of showing the existence of a lift to
vector-valued versions of \eqref{eqnrfs} to a process with values in a
rough paths space.

Recall that we write $(a_k) \preceq(b_k)$ for two sequences $(a_k)$
and $(b_k)$ if there is a negligible sequence $(c_k)$ such that $a_k =
c_k b_k$; cf. Definition~\ref{defnegli}. We will extend this notation
as follows: if $(A_k) = (a_k^{i,j})$ is a sequence
of matrices, and $(b_k)$ is a sequence of real numbers, $(A_k) \preceq(b_k)$
means that $(a^{i,j}_k) \preceq(b_k)$ for every $i,j$. If $(A_k) =
(A^1_k,\ldots,A^m_k)$ is a sequence of vectors whose entries are
matrices or
real numbers, we will write $(A_k) \preceq(b_k)$ if $(A^i_k) \preceq(b_k)$
for all $i=1,\ldots,m$.

Let $\Psi= (\Psi^1, \ldots, \Psi^d)$ where the $\Psi^i$ are given
as random
Fourier series
%
\begin{equation}
\label{eqnrfscoord} \Psi^i(t)=\frac{\a^i_{0}Y^{0,i}}{2}+\sum
_{k=1}^\infty\a^i_{k}Y^{k,i}\sin(kt)+\a^i_{-k}Y^{-k,i}\cos(kt),
\end{equation}
with\vspace*{1pt} $(Y^{k,i})_{k\in\Z,i=1,\ldots,d}$ being independent, $\mathcal{N}(0,1)$
distributed random variables. As before, set $a^i_k:= (\alpha_k^i)^2$
and $%
(a_k):= (a^1_k,\ldots,a^d_k)$. Our main existence result is the following:

%
\begin{theorem}
\label{teoexistencerfsasrp} Assume $(a_k) \preceq(\llvert k\rrvert
^{-(1+1/\rho)})$
for some $\rho\in[1,2)$ with associated measures $\mu^i_1, \mu
^i_2$, $%
i=1,\ldots,d$, as in Definition~\ref{defnegli}, and let $K\geq
\max_{i=1,\ldots,d}\{ \llVert \mu^i_1\rrVert _{\mathrm{TV}}, \llVert
\mu^i_2\rrVert
_{\mathrm{TV}}\}$.
Then for every $\beta<\frac{1}{2\rho}$, there exists a continuous $%
G^{[1/\beta]}(\R^d)$-valued process $\bolds{\Psi}$ such that:
\begin{longlist}[(3)]
\item[(1)] $\bolds{\Psi}$ has geometric $\beta$-H\"older rough sample paths,
that is,
\[
\bolds{\Psi} \in C_0^{0, \beta\mbox{-H\"ol}}\bigl([0,2\pi
],G^{[1/\beta]}\bigl(\R^d\bigr)\bigr)
\]
almost surely,

\item[(2)] $\bolds{\Psi}$ lifts $\Psi$ in the sense that $\pi_1(\bolds
{\Psi}_t)
= \Psi_t-\Psi_0 $,

\item[(3)] there is a $C = C(\rho,K)$ such that for all $s < t$ in $[0,2\pi
]$ and
$q \in[1,\infty)$,
\[
\bigl\llvert d(\bolds{\Psi}_s,\bolds{\Psi}_t)\bigr
\rrvert_{L^q} \le C \sqrt{q} \llvert t-s\rrvert^{1/(2\rho)},
\]

\item[(4)] there exists $\eta= \eta(\rho,K,\beta) > 0$, such that
\[
\E e^{\eta\llVert \bolds{\Psi}\rrVert _{\beta\mbox{-H\"ol};[0,2\pi
]^2}} < \infty.
\]
\end{longlist}
\end{theorem}

\begin{pf}
By assumption,
\begin{eqnarray*}
\Psi^i(t)&=&\frac{\gamma^i_{0}Y^{0,i}}{2}+\sum_{k=1}^\infty
\gamma^i_{k} \llvert k\rrvert^{- ({1/2} + {1/(2\rho)} )}
Y^{k,i}\sin(kt)
\\
&&{}+\gamma^i_{-k} \llvert k\rrvert
^{- (1/2 + {1/(2\rho)})} Y^{-k,i}\cos(kt)
\end{eqnarray*}
for every $i=1,\ldots,d$ where $(c^i_k)=((\gamma_k^i)^2)$ is a negligible
sequence. Hence, we may apply Proposition~\ref{propmollRFS} to see that
the covariance of $\Psi^i$ has finite H\"older dominated $\rho$-variation
for every $i$; thus \cite{FV10-2}, Theorem 35, applies.
\end{pf}

We\vspace*{2pt} will now compare the lifts of two random Fourier series $\Psi=
(\Psi^1,\ldots,\Psi^d)$ and $\tilde{\Psi} = (\tilde{\Psi
}^1,\ldots,\tilde{%
\Psi}^d)$ with
\begin{eqnarray*}
\Psi^i(t) &=& \frac{\a^i_{0}Y^{0,i}}{2}+\sum_{k=1}^\infty
\a%
^i_{k}Y^{k,i}\sin(kt)+
\a^i_{-k}Y^{-k,i}\cos(kt),
\\
\tilde{\Psi}^i(t) &=& \frac{\tilde{\a}^i_{0}\tilde{Y}^{0,i}}{2}%
+\sum
_{k=1}^\infty\tilde{\a}^i_{k}
\tilde{Y}^{k,i}\sin(kt) + \tilde{\a}%
^i_{-k}
\tilde{Y}^{-k,i}\cos(kt).
\end{eqnarray*}
We make the following assumption:
\begin{eqnarray*}
&&\bigl\{\bigl(Y^{k,i},\tilde{Y}^{k,i}\bigr)\dvtx k\in\Z, i=1,
\ldots,d\bigr\}
\end{eqnarray*}
are independent, normally distributed random vectors with $Y^{k,i},
\tilde{Y}%
^{k,i} \sim\break \mathcal{N}(0,1)$ for all $k\in\Z$ and $i=1,\ldots,d$. It
follows\vspace*{1pt} that $\E Y^{k,i} \tilde{Y}^{l,j} = 0$ for $k\neq l$ or $i \neq j$,
and we set $\varrho_k^i:= \E Y^{k,i} \tilde{Y}^{k,i}$. As before,
let $%
a^i_k:= (\alpha_k^i)^2$ and $\tilde{a}^i_k:= (\tilde{\alpha}_k^i)^2$.
Define the matrix
\begin{eqnarray*}
A_k^i:= \pmatrix{ a_k^i &
\alpha_k^i \tilde{\alpha}_k^i
\varrho_k^i
\vspace*{3pt}\cr
\alpha_k^i \tilde{
\alpha}_k^i \varrho_k^i &
\tilde{a}_k^i},
\end{eqnarray*}
and set $A_k:= (A_k^1,\ldots,A_k^d)$.

%
\begin{theorem}
\label{teodistliftRFSasRP} Assume that $(A_k) \preceq
(\llvert k\rrvert ^{-(1+1/\rho)})$ for some $\rho\in[1,2)$ and that
the total variation
of all associated measures is bounded by a constant $K$. Then we can
lift $%
\Psi$ and $\tilde{\Psi}$ to processes with values in a rough paths
space as
in Theorem~\ref{teoexistencerfsasrp}, and for all $\gamma< 1 -
\frac{\rho%
}{2}$ and $\beta< \frac{1}{\rho} (\frac{1}{2} - \gamma
)$ there is
a constant $C = C(\rho,K,\beta,\gamma)$ such that
%
\begin{eqnarray}
\label{eqnestimaterpdistancelq} &&\bigl\llvert\rho_{\beta\mbox{-H\"
ol}}(\bolds{\Psi},\bolds{\tilde{
\Psi}})\bigr\rrvert_{L^q} \leq C q^{(1/2) \lfloor {1/\beta
} \rfloor} \Bigl(\sup
_{t\in[0,2\pi]} %
\E\bigl\llvert\Psi(t) - \tilde{\Psi}(t)
\bigr\rrvert^2 \Bigr)^{\gamma}
\end{eqnarray}
for all $q\in[1,\infty)$.
\end{theorem}

\begin{pf}
The existence of the lifted processes $\bolds{\Psi}$ and $\bolds
{\tilde{%
\Psi}}$ follows from Theorem~\ref{teoexistencerfsasrp}. The $L^q$
norm of the difference of two such processes in rough paths metric can
be estimated by the $\rho$-variation of the covariance of the
difference of the two processes, and an interpolation argument shows
that this quantity can actually be bounded by the right-hand side of
\eqref{eqnestimaterpdistancelq} times the $\rho$-variation of the
covariance of the two processes\vspace*{1pt} and their joint covariance function. We
aim to apply
\cite{FR12a}, Theorem 5,\footnote{%
Strictly speaking, \cite{FR12a}, Theorem 5, assumes that $\tilde{\Psi
}$ is a
certain approximation of $\Psi$. However, it is shown in \cite{RX12} that
this is not necessary, and (\cite{FR12a}, Theorem 5) can be used more generally
to give an upper bound for the distance between $\bolds{\Psi}$ and
$\tilde{\bolds{\Psi}}$ as we need it here.} where the estimate \eqref
{eqnestimaterpdistancelq} was given for the optimal parameter
$\gamma$. To obtain a uniform estimate, we need to show that the
joint covariance function of the process $(\Psi,\tilde{\Psi})$ has finite,
H\"older dominated $\rho$-variation, bounded by a constant depending
only on
the parameters above. From independence of the components, it suffices to
estimate the $\rho$-variation of $R_{\Psi^i,\tilde{\Psi}^i}(s,t) =
\E%
\Psi^i(s) \tilde{\Psi}^i(t)$ for every $i=1,\ldots,d$. This can be done
using Proposition~\ref{proppairRFS}.
\end{pf}

As an application, we consider the \emph{truncated random Fourier series},
that is, we define $\Psi^N = (\Psi^{1,N},\ldots,\Psi^{d,N})$ by
%
\begin{eqnarray}
&&\Psi^{i,N} (t) = \frac{\a^i_{0}Y^{0,i}}{2}+\sum_{k =1 }^N
\a^i_{k}Y^{k,i}\sin(kt)+
\a^i_{-k}Y^{-k,i}\cos(kt)
\nonumber\\[-12pt]\\[-8pt]
\eqntext{\mbox{for } i=1, \ldots,d.}
\end{eqnarray}
It is then easy to show that convergence also holds for the corresponding
rough paths lifts, and we can even give an upper bound for the order of
convergence.

%
\begin{corollary}
\label{corconvtruncFR} Under the assumptions of Theorem~\ref{teoexistencerfsasrp}, choose some $\eta< \frac{1}{\rho} - \frac{1}{2}$
and $\beta< \frac{1}{2\rho} - \eta$. Then there is a constant $C =
C(\rho,K,\beta,\eta)$ such that
\begin{eqnarray*}
\bigl\llvert\rho_{\beta\mbox{-H\"ol}}\bigl(\bolds{\Psi},\bolds{\Psi}^N
\bigr)\bigr\rrvert_{L^q} \leq C q^{(1/2) \lfloor{1/\beta} \rfloor} \biggl(
\frac{1}{N} \biggr)^{\eta}
\end{eqnarray*}
for every $N \in\N$, $q \in[1,\infty)$. In particular, $\rho
\mbox{\tsub{$\beta$-\textup{H\"ol}}}(\bolds{\Psi},\bolds{\Psi}^N) \to0$ for
$N \to\infty$ almost
surely and in $L^q$ for any $q\in[1,\infty)$ with rate $\eta$.
\end{corollary}

\begin{remark}
We emphasize that $\Psi,\Psi^N$ are lifted to level $\lfloor1/\beta
\rfloor$
above. In particular, a ``good'' rate $\eta$ forces $\beta$ to be
small so
that, in general, it is not enough to work with 3 levels, as is the usual
setting in Gaussian rough paths theory.
\end{remark}

\begin{pf*}{Proof of Corollary \protect\ref{corconvtruncFR}}
We aim to apply Theorem~\ref{teodistliftRFSasRP} with $\tilde
{\alpha}%
^i_k = \mathbbm{1}_{\llvert k\rrvert \leq N} \alpha^i_k$ and $\varrho
^i_k \equiv
1$. We
will first show that $(a^i_k \mathbbm{1}_{\llvert k\rrvert \leq N})
\preceq
(\llvert k\rrvert ^{-(1 +
1/\rho')})$ for every $\rho'> \rho$, uniformly
over $i$
and $N$. Indeed, we have
\begin{eqnarray*}
&& a_k^i \mathbbm{1}_{\llvert k\rrvert \leq N} =
\bigl(a^i_k \llvert k\rrvert^{1+1/\rho}\bigr) \bigl(
\llvert k\rrvert^{1/\rho
^{\prime
}- 1/\rho}\mathbbm{1}_{\llvert k\rrvert \leq N}\bigr) \llvert k
\rrvert^{-(1+1/\rho')},
\end{eqnarray*}
and since $(a^i_k) \preceq(\llvert k\rrvert ^{-(1+1/\rho)})$ for all
$i=1,\ldots,d$, it
suffices to show that $(\llvert k\rrvert ^{-\varepsilon}\mathbbm
{1}_{\llvert k\rrvert \leq N})$ is
uniformly negligible for every $\varepsilon> 0$. Therefore, we need to show
that the cosine series
\begin{eqnarray*}
B^N(x) = \sum_{k=1}^{\infty}
\llvert k\rrvert^{-\varepsilon}\mathbbm{1}_{\llvert k\rrvert \leq N}
\cos(kx) = \sum
_{k=1}^{N} \llvert k\rrvert
^{-\varepsilon} \cos(kx)
\end{eqnarray*}
is uniformly bounded in $L^1([0,2\pi])$. Since $\D k^{-\varepsilon} =
O(k^{-\varepsilon-1})$ and\break $\lim_{k \to\infty} \log
(k)k^{-\varepsilon} = 0$,
we can apply the Sidon--Telyakovskii theorem (cf. \cite{T73}, Theorem~4) to
obtain $B^N \to B$ for $N\to\infty$ in $L^1([0,2\pi])$ which proves the
uniform negligibility, and we may apply Theorem~\ref{teodistliftRFSasRP}
for every $\rho'> \rho$. Furthermore,
\begin{eqnarray*}
\E\bigl\llvert\Psi(t) - \Psi^{N}(t)\bigr\rrvert^2 &=&
\sum_{k= N + 1}^\infty a_k \sin
^2(kt) + a_{-k} \cos^2(kt)
\\
&\leq&2 \sum_{\llvert k\rrvert \geq N + 1} a_k \lesssim4 \sum
_{k = N + 1}^{\infty} k^{-(1+1/\rho)} \lesssim
\biggl(\frac{1}{N} \biggr)^{1/\rho}.
\end{eqnarray*}
For given $\eta$, we choose $\rho'$ such that $\eta< \frac
{1}{%
\rho'} - \frac{1}{2} < \frac{1}{\rho} - \frac{1}{2}$ and apply
Theorem~\ref{teodistliftRFSasRP} to complete the $L^q$
convergence. The
almost sure convergence follows by a standard Borel--Cantelli argument;
cf. \cite{FR12a}, Theorem 6, page 41.
\end{pf*}

\section{Applications to SPDE}\label{sec4}\label{secapplication}

In this section we will apply our results on random Fourier series to
construct spatial rough path\vspace*{1pt} lifts of stationary Ornstein--Uhlenbeck
processes corresponding to the $\R^{d}$-valued (generalized) fractional
stochastic heat equation with Dirichlet, Neumann or periodic boundary
conditions
%
\begin{eqnarray}
\label{eqnfSHE} \,d\Psi_{t}=-(-\Delta)^{\alpha}\Psi_{t}
\,dt+dW_{t}\qquad\mbox{on } [0,T]\times%
[0,2\pi],
\end{eqnarray}
where the fractional Laplacian $(-\Delta)^{\alpha}$ acts on each component
of $\Psi_{t}$ and $\a\in(0,1]$. We will start by first considering the
fractional stochastic heat equation with Dirichlet boundary conditions,
proving the existence of (continuous) spatial rough paths lifts and
stability under approximations. Then we comment on Neumann boundary
conditions and on more general equations for periodic boundary conditions.

If a (spatial) rough path lift of \eqref{eqnfSHE} has been
constructed, one
can view \eqref{eqnfSHE} as an evolution in a rough path space, a
point of
view which has proven extremely fruitful in solving new classes of,
until now,
ill-posed stochastic PDE \cite{H11,HW11,H12}, arising, for example, in path
sampling problems for $\R^d$-valued SDE \cite{HSVW05,HSV07,H11}.

For a variant of \eqref{eqnfSHE} with $\a=1$, Hairer proved in \cite{H11}
finite $1$-variation of the covariance of the stationary solution to %
\eqref{eqnfSHE}, that is, of $ ( x,y ) \mapsto\E\Psi
(
t,x ) \Psi( t,y )$. This general theory then gives
the existence
of a ``canonical, level 2'' rough path $\bolds{\Psi}$ lifting $\Psi
$; cf. Theorem~\ref{teoexistencerfsasrp}; see also \cite{GIP12}.
It is
clear that in the case $\a=1$ the Brownian-like regularity of $x
\mapsto
\Psi( t,x;\omega) $ is due to the competition between the
smoothing effects of the Laplacian and the roughness of space--time white
noise. Truncation of the higher noise modes (or suitable ``coloring'') leads
to better spatial regularity; on the other hand, replacing $\Delta$ by a
fractional Laplacian, that is, considering \eqref{eqnfSHE} for some
$\alpha
\in(0,1)$, dampens the smoothing effect, and $x \mapsto\Psi(
t,x;\omega) $ will have ``rougher'' regularity properties than a standard
Brownian motion. One thus expects $\rho$-variation regularity for the
spatial covariance of $x \mapsto\Psi( t,x;\omega) $
for %
\eqref{eqnfSHE} only for some $\rho>1$ and subsequently only the existence
of a ``rougher'' rough path, that is, necessarily with higher $p$ than before.

As we shall see below, \eqref{eqnfSHE} is handled as a spatial rough path
with a number of precise estimates, provided
\[
\alpha>\alpha^{\ast}= \tfrac{3}{4}.
\]
More precisely, the resulting (geometric rough) path enjoys $\frac
{1}{p}$-H%
\"{o}lder regularity for any $p>2\rho= \frac{2}{2\alpha-1}$. When
$\alpha>
\frac{5}{6}$ we have $\rho= \frac{1}{ 2\alpha-1} < \frac{3}{2}$
and can\vspace*{1pt}
pick $p<3$. The resulting rough path can then be realized as a ``level 2''
rough path. In the general case (similar to $H\in(\frac{1}{4},\frac{1}{3}]$
in the fBm setting) one must go beyond the stochastic area and control the
third level iterated integrals. Our approach, which crucially passes through
$\rho$-variation, combined with existing theory, has many advantages: the
notoriously difficult third-level computation need not be repeated in the
present context; leave alone the higher level computations needed for rates.
A~satisfactory approximation theory is also available, based on \textit
{%
uniform} $\rho$-variation estimates; cf. Section~\ref{subsecstability}
below.

\subsection{Fractional stochastic heat equation with Dirichlet boundary
conditions}\label{sec41}
\label{secFSHEdirichlet}

We consider
%
\begin{eqnarray}
\label{eqnfSHEagaub} d\Psi_{t}=-(-\Delta)^{\alpha}\Psi_{t}
\,dt+dW_{t}\qquad\mbox{on } [0,T]\times%
[0,2\pi]
\end{eqnarray}
on $[0,2\pi]$ endowed with Dirichlet boundary conditions. Neumann and
periodic boundary conditions may be treated analogously; cf. Section
\ref{secperiodicSFHE} below. We have the following orthogonal basis of
eigenvectors with corresponding eigenvalues of $-\D$ on $L^{2}([0,2\pi])$:
\begin{eqnarray*}
e_{k}(x)=\sin\biggl(\frac{k}{2}x \biggr),\qquad
\tau_{k}= \biggl(\frac
{k}{2}%
\biggr)^{2},\qquad
k\in\N,
\end{eqnarray*}
and take $W_t = \sum_{k\in\N} \b^k_t e_k(x)$. The fractional
Laplacian has
eigenvalues $\l_k:= \tau_k^\a$ for $k\in\N$ and (informal) Fourier
expansion of the stationary solution $\Psi$ to \eqref{eqnfSHE} leads
to the
random Fourier series
%
\begin{eqnarray}
\label{eqnfourierdecomp1} \Psi(t,x) = \sum_{k =1}^\infty
\a_k Y_t^{k}\sin\biggl(k\frac
{x}{2}
\biggr),
\end{eqnarray}
with $\a_k = \frac{1}{\sqrt{2\l_k}}$ and $Y_t^k$ being a decoupled,
infinite system of $d$-dimensional, stationary, normalized
Ornstein--Uhlenbeck processes satisfying
%
\begin{equation}
\label{eqnOU} d Y^{k}_t = - \l_k
Y^{k}_t \,dt + \sqrt{2\l_k} \,d
\b_t^k.
\end{equation}
Clearly (\ref{eqnfourierdecomp1}) gives a well-defined and continuous
random field and solves (\ref{eqnfSHEagaub}) in the sense of standard SPDE
theory; cf., for example, \cite{W86,DPZ92}. Note $\E Y^{k}_t \otimes
Y^{l}_s =
e^{-\l_k\llvert t-s\rrvert } \delta_{k,l}\operatorname{Id}$, and set
\begin{eqnarray*}
a_k &=& \a_k^2 = \frac{1}{2\l_k} =
2^{2\a-1}\frac{1}{k^{2\a}},\qquad k\in\N.
\end{eqnarray*}
As an immediate consequence of our results on random Fourier series, we
get the following:

\begin{proposition}
\label{propspatiallift0} Suppose $\a\in(\frac{1}{2},1]$. Then:
\begin{longlist}[(1)]
\item[(1)] For every $t \geq0$, the spatial process $x\mapsto\Psi(t,x)$
is a
centered Gaussian process which admits a continuous modification (which we
denote by the same symbol) with covariance $R_\Psi$ of finite mixed $%
(1,\rho) $-variation for all $\rho\geq\frac{1}{2\alpha- 1}$, and all
conclusions of Theorem~\ref{teorfs} hold.

\item[(2)] If $\alpha> \frac{3}{4}$, the process $x \mapsto\Psi(t,x)$
lifts to
a process with geometric $\beta$-H\"older rough paths
\[
\bolds{\Psi}(t) \in C_0^{0,\beta\mbox{-H\"ol}}\bigl([0,2\pi
],G^{\lfloor
1/\beta\rfloor}\bigl(\R^d\bigr)\bigr)
\]
almost surely for every $\beta< \alpha- \frac{1}{2}$.

\item[(3)] Choose $\gamma$ and $\beta$ such that
\begin{eqnarray*}
\gamma< 1 - \frac{3}{4\alpha},\qquad\beta< \alpha- \frac{1}{2} -
\frac{%
2\alpha\gamma}{2\alpha- 1}.
\end{eqnarray*}
Then there is a $\gamma$-H\"older continuous modification of the map
%
\begin{eqnarray}
\bolds{\Psi}\dvtx [0,T] &\to& C_0^{0,\beta\mbox{-H\"ol}}\bigl([0,2\pi
],G^{\lfloor
1/\beta\rfloor}\bigl(\R^d\bigr)\bigr),
\nonumber\\[-8pt]\\[-8pt]\nonumber
t &\mapsto& \bolds{\Psi}(t).
\end{eqnarray}
\end{longlist}
\end{proposition}

\begin{remark}
In (3), we observe a ``trade-off'' between the parameters $\beta$~and~$
\gamma$: If we want a ``good'' time regularity (i.e., large $\gamma
$), we
have to take $\beta$ small which is tantamount to working in a rough paths
space with many ``levels'' of formal iterated integrals. For instance,
when $%
\alpha=1$, we can get arbitrarily close to $\frac{1}{4}$ in time, at the
price of working with many arbitrary levels. On the other hand, if we insist
to work with the first 3 levels only (or 2 levels in case $\alpha>5/6$),
which is the standard setting in Gaussian rough path theory, we only get
poor time regularity of the evolution in rough path space.
\end{remark}

\begin{pf*}{Proof of Proposition \protect\ref{propspatiallift0}}
Since $\Psi$ is a rescaling of
\[
\td\Psi(t,x) = \sum_{k=1}^\infty
\a_k Y_t^{k}\sin(k x )=\Psi(t,2x),
\]
it is enough to consider $\td\Psi$:

\begin{longlist}[(1)]
\item[(1)] Clearly $x \mapsto\td\Psi(t,x)$ is centered and Gaussian. Due to
\eqref{eqnrfsRdecomp} and Lemma \ref{lemmahoelder}, we have
\begin{eqnarray*}
\s^2_t(x,y)=\E\bigl\llvert\td\Psi(t,x) - \td\Psi(t,y)
\bigr\rrvert^2 \lesssim\llvert x-y\rrvert^{2\alpha- 1},
\end{eqnarray*}
which implies that there is a continuous modification of $\td\Psi$.
Theorem~\ref{teorfs} implies the claim.

\item[(2)] Follows from Theorem~\ref{teoexistencerfsasrp}.

\item[(3)] We will derive the existence of a continuous modification by
application of Kolmogorov's continuity theorem. Therefore, we need an
estimate on a $q$th moment of the distance in the $\rho_{\beta\mbox{-}\mathrm{H\ddot{o}l}}$
metric of the rough paths $\bolds{\td\Psi}(t)$, $\bolds{\td\Psi
}(s)$ at
different times $0 \le s < t \le T$. Such an estimate can be obtained by
applying Theorem~\ref{teodistliftRFSasRP}. Let $0\le s\le t \le
T$, $%
\tau:=\llvert t-s\rrvert $, and set $A_k = (A^1_k,\ldots,A^d_k)$ where
\begin{eqnarray*}
A_k^i:= \pmatrix{ a_k & a_k
e^{-\lambda_k \tau}
\cr
a_k e^{-\lambda_k \tau} & a_k}
\end{eqnarray*}
for $i=1,\ldots,d$. We claim that $(A_k) \preceq(\llvert k\rrvert
^{-2\alpha})$ uniformly
in $\tau$. Defining $b(r):= e^{- ({r/2} )^{2\a}}$
we note $%
b^\tau_k=e^{-\lambda_k\tau}=b(k\tau^{1/(2\a)})$, and $b$ is convex,
nonincreasing. Lemma~\ref{lemmauniformseriesbound} then implies
that $%
(e^{-\lambda_k\tau})$ is uniformly negligible which shows the claim. Hence,
we can apply Theorem~\ref{teodistliftRFSasRP} and obtain
\[
\bigl\llvert\rho\mbox{\tsub{$\beta$-H\"ol}}\bigl(\bolds{\tilde{\Psi
}}(t),\bolds{\tilde{
\Psi}}%
(s)\bigr)\bigr\rrvert_{L^q} \leq C q^{(1/2) \lfloor {1/\beta
} \rfloor}
\Bigl(\sup_{x \in[0,2\pi]} \E\bigl\llvert\tilde{\Psi}(t,x) - \tilde{
\Psi}(s,x) \bigr\rrvert^2 \Bigr)^{\theta}
\]
for all $\theta< \frac{4\alpha- 3}{4\alpha- 2}$, $\beta< \alpha-
\frac{1%
}{2} - \theta$ and all $q \in[1,\infty)$. In order to estimate the
right-hand side, we note
\begin{eqnarray*}
&& \E\bigl\llvert\tilde{\Psi}^1(t,x) - \tilde{\Psi}^1(s,x)
\bigr\rrvert^2
\\
&&\qquad = \E\bigl\llvert\tilde{\Psi}%
^1(t,x)
\bigr\rrvert^2 + \E\bigl\llvert\tilde{\Psi}^1(s,x)\bigr
\rrvert^2 - 2\E\tilde{\Psi}^1(t,x) \tilde{
\Psi}%
^1(s,x)
\\
&&\qquad \leq 2\sum_{k=1}^\infty a_k
\bigl(1-e^{-\lambda_k\tau} \bigr) \sin^2(kx)
\le 2\sum_{k=1}^\infty a_k
\bigl\llvert1-e^{-\lambda_k\tau}\bigr\rrvert
\\
&&\qquad \le C\sum_{k \le N} \llvert t-s\rrvert+
CN^{1 - 2\alpha'} \sum_{k > N} a_k
k^{2\alpha'- 1}
\le C \bigl(N \llvert t-s\rrvert+ N^{1 - 2\alpha'}\bigr)
\end{eqnarray*}
for all $\alpha' < \alpha$. We then choose $N \sim
\llvert t-s\rrvert ^{-1/(2\alpha')}$ to obtain
\begin{eqnarray*}
\E\bigl\llvert\tilde{\Psi}(t,x) - \tilde{\Psi}(s,x) \bigr\rrvert^2
\le C\llvert t-s\rrvert^{1 -
1/(2\alpha')}.
\end{eqnarray*}
Thus we can choose $\gamma< 1 - \frac{3}{4\alpha}$ and $\beta<
\alpha-
\frac{1}{2} - \frac{2\alpha\gamma}{2\alpha- 1}$ to obtain
\[
\bigl\llvert\rho\mbox{\tsub{$\beta$-H\"ol}}\bigl(\bolds{\tilde{\Psi
}}(t),\bolds{\tilde{
\Psi}}%
(s)\bigr)\bigr\rrvert_{L^q} \le Cq^{(1/2) \lfloor{1/\beta}
\rfloor}
\llvert t-s\rrvert^{\gamma},
\]
for all $q \in[1,\infty)$. Kolmogorov's continuity theorem gives the result.\quad\qed
\end{longlist}\noqed
\end{pf*}

\subsection{Stability and approximations}\label{sec42}\label{subsecstability}

Due to the ``contraction principle'' in the form of Proposition~\ref{proprhovariationconv}, the estimates on the $\rho$-variation of the
covariance of random Fourier series derived in Section~\ref{secrfs} are
robust with respect to approximations. In order to emphasize this
point, in
this section we consider Galerkin and hyper-viscosity approximations to
$%
\Psi$ with $\Psi$ as in Section~\ref{secFSHEdirichlet} and prove strong
convergence of the corresponding rough paths lifts. Recall that by the
general theory of rough paths, this immediately implies the strong
convergence of the corresponding stochastic integrals as well as of
solutions to rough differential equations; cf., for example, \cite{BFRS13,H11}.

\subsubsection{Galerkin approximations}\label{sec421}

The Galerkin approximation $\Psi^{N}_t$ of $\Psi_t$ is defined to be the
projection of $\Psi$ onto the $N$-dimensional subspace spanned by
$ \{
e_k \}_{k = 1,\ldots,N}$. This process solves the SPDE
%
\begin{equation}
\label{eqngenheatgal} d\Psi^{N}_t = -\bigl(P_N(-
\Delta)^{\alpha}\bigr) \Psi^{N}_t \,dt +
dP_NW_t,
\end{equation}
where $P_N(-\Delta)^{\alpha}$ has the eigenvalues $ (\frac{k}{2}
)^{2\alpha} \mathbbm{1}_{k\leq N}$, and $P_N W_t$ has the covariance
operator $Q^N$ given by $Q^N e_k = \mathbbm{1}_{k \leq N} e_k$. The
process $%
\Psi^{N}$ can be written as the truncated Fourier series
\begin{eqnarray*}
\Psi^{N}(t,x) = \sum_{k = 1}^N
\alpha_k Y^{k}_t \sin\biggl(k
\frac
{x}{2}%
\biggr),
\end{eqnarray*}
with $\alpha_k = 2^{\alpha- 1/2} k^{-\alpha}$ and $Y^{k}$
Ornstein--Uhlenbeck processes as in Section~\ref{secFSHEdirichlet}.

One easily checks that we can lift the spatial sample paths of $\Psi^{N}_t$
to Gaussian rough paths and find continuous modifications of $t \mapsto
\bolds{\Psi}^{N}_t$. Moreover, we can prove the following strong
convergence result:

\begin{proposition}
\label{propgalerkinapprox} Assume $\alpha> \frac{3}{4}$, and
choose $\eta
< 2\alpha- \frac{3}{2}$ and $\beta< \alpha- \frac{1}{2} - \eta$. Then
there is a constant $C=C(\alpha,\beta,\eta)$ such that
\begin{eqnarray*}
\bigl\llvert\rho\mbox{\tsub{$\beta$-\textup{H\"ol}}}\bigl(\bolds{\Psi}(t),\bolds{\Psi
}^N(t)\bigr)\bigr\rrvert_{L^q} \leq C q^{(1/2) \lfloor
{1/\beta} \rfloor}
\biggl(\frac{1}{N} \biggr)^{\eta}
\end{eqnarray*}
for all $t\in[0,T]$, $N\in\N$, $q \in[1,\infty)$. In particular, for
every $t\in[0,T]$,\break $\rho\mbox{\tsub{$\beta$-\textup{H\"ol}}}(\bolds{\Psi
}(t),\bolds{\Psi}%
^N(t)) \to0$ for $N\to\infty$ almost surely and in $L^q$ for any $q
\in
[1,\infty)$ with rate $\eta$.
\end{proposition}

\begin{pf}
The proof follows from Corollary~\ref{corconvtruncFR}.
\end{pf}

\subsubsection{Hyper-viscosity approximations}\label{sec422}

The hyper-viscosity approximation $\Psi^{\ve} = (\Psi^{\ve,1},\ldots
,\Psi^{%
\ve,d})$ is the solution to
%
\begin{equation}
\label{eqngenheathypervisc} d\Psi^{\ve}_{t}=-\bigl((-
\Delta)^{\alpha} + \ve(-\Delta)^{\theta}\bigr) \Psi^{\ve%
}_{t}
\,dt+dW_{t},
\end{equation}
for some (large) $\theta\geq1$ and $\ve>0$. Again, it is easy to see that
we can lift the spatial sample paths of $\Psi^{\ve}_t$ to Gaussian rough
paths and find continuous modifications of $t \mapsto\bolds{\Psi
}^{\ve}_t$.

\begin{proposition}
\label{prophyperviscapprox} Assume $\alpha> \frac{3}{4}$ and
$\theta>
\alpha$. Choose $\beta< \alpha- \frac{1}{2}$. Then there is a
function $%
r_{\alpha,\beta,\theta} \colon\R\to\R_+$ such that $r(\ve) \to
0$ for $%
\ve\to0$ and a constant $C=C(\alpha,\beta,\theta)$ such that
\begin{eqnarray*}
\bigl\llvert\rho\mbox{\tsub{$\beta$-\textup{H\"ol}}}\bigl(\bolds{\Psi}(t),\bolds{
\Psi}^{\ve
}(t)\bigr)\bigr\rrvert_{L^q} \leq C
q^{(1/2) \lfloor{1/\beta} \rfloor} r(\ve)
\end{eqnarray*}
for every $t\in[0,T]$, $\ve> 0$ and $q \in[1,\infty)$.
\end{proposition}

\begin{pf}
As before, $\Psi^{\ve}_t$ has the form of a random Fourier series
where the $%
k$th Fourier coefficients are given by $\alpha_k^{\ve} Y^{\ve,k}_t$
with $%
\alpha^{\ve}_k = \frac{1}{\sqrt{2 \lambda_k^{\ve}}}$,
\begin{eqnarray*}
\lambda^{\ve}_k = \biggl(\frac{k}{2}
\biggr)^{2\alpha} + \ve\biggl(\frac{k}{2}%
\biggr)^{2\theta},
\end{eqnarray*}
and $t\mapsto Y^{\ve,k}_t$ are $d$-dimensional, stationary
Ornstein--Uhlenbeck processes with independent components, each component
being centered with variance 1 and correlation
\begin{eqnarray*}
\E Y^{\ve,k}_t \otimes Y^{l}_t = 2
\frac{\sqrt{\lambda_k \lambda
_k^{\ve}}}{%
\lambda_k + \lambda_k^{\ve}} \delta_{k,l}\mathrm{Id}.
\end{eqnarray*}

From Theorem~\ref{teodistliftRFSasRP}, we know that it is
sufficient to
show that $(A^{\ve}_k) \preceq(\llvert k\rrvert ^{-2\alpha})$
uniformly over $\ve
> 0$
where
\begin{eqnarray*}
A_k^{\ve}:= \pmatrix{ \alpha_k^2
& \alpha_k \alpha_k^{\ve} \varrho_k^{\ve}
\cr
\alpha_k \alpha_k^{\ve}
\varrho_k^{\ve} & \bigl(\alpha_k^{\ve}
\bigr)^2}, \qquad\varrho_k^{\ve}:= 2
\frac{\sqrt{\lambda_k \lambda_k^{\ve
}}}{%
\lambda_k + \lambda_k^{\ve}}
\end{eqnarray*}
and that
%
\begin{eqnarray}
\label{eqnhverrorbound} \sup_{t\in[0,T]} \sup_{x\in[0,2\pi]} \E
\bigl\llvert\Psi(t,x) - \Psi^{\ve
}(t,x)\bigr\rrvert^2 \to
0\qquad\mbox{for } \ve\to0.
\end{eqnarray}
For the first claim, we have to show that the series
\begin{eqnarray*}
\sum_{k=1}^{\infty} \bigl(
\alpha_k^{\ve}\bigr)^2 k^{2\alpha} \cos
(kx),\qquad\sum_{k=1}^{\infty}
\alpha_k \alpha_k^{\ve} \varrho_k^{\ve}
k^{2\alpha} \cos(kx)
\end{eqnarray*}
are uniformly bounded in $L^1$ which can be done using Lemma~\ref{lemmauniformseriesbound}(2). Showing~\eqref{eqnhverrorbound} follows
by writing down the left-hand side as a Fourier series and bounding it
uniformly in $x$ and $t$ by an infinite series. Then we send $\ve\to0$,
using the dominated convergence theorem.
\end{pf}

\subsection{Various generalizations}\label{sec43}\label{secperiodicSFHE}

\subsubsection{Generalized fractional stochastic heat equation on periodic
domains}\label{sec431}

Based on the stability results for the mixed $(1,\rho)$-variation of the
covariance of random Fourier series developed in Section~\ref
{secrfs}, one
may consider more general fractional stochastic heat type equations and
different types of boundary conditions. As an example let us consider
generalized fractional stochastic heat equations on $[0,2\pi]$ with periodic
boundary conditions. An orthogonal basis of eigenvectors and corresponding
eigenvalues of $-\D$ on $L^2([0,2\pi])$ is given by
\[
\tau_k = k^2, \qquad e_k(x):= \cases{
\sin(kx), &\quad $k>0$,
\cr
\frac{1}{2}, &\quad $k=0$,
\cr
\cos(kx), &\quad $k<0$.}
\]
Via the spectral theorem we may define $A=f(-\D)$ for each Borel measurable
function $f\dvtx\R_+\to\R_+$, still having $e_k$ as a basis of
eigenvectors and
eigenvalues $f(\tau_k)$. %

In order to be able to consider stationary Ornstein--Uhlenbeck
processes, we
need to shift the spectrum of $A$ to be strictly negative. Hence, we
consider \mbox{$\R^d$-}valued SPDE of the form
%
\begin{equation}
\label{eqngenheat} d\Psi_t = (-A - \lambda)\Psi_t \,dt +
dW_t\qquad\mbox{with } \l> 0,
\end{equation}
where $W_t$ is a (possibly) colored Wiener process with covariance operator
having $e_k$ as basis of eigenvectors and $\s_k$ as eigenvalues. An
(informal) Fourier expansion of the stationary solution $\Psi$ to %
\eqref{eqngenheat} leads to the random Fourier series
%
\begin{eqnarray}
\label{eqnfourierdecomp2} \Psi(t,x) = \frac{\a_0Y_0}{2}+\sum
_{k=1}^\infty
\a_k Y_t^{k}\sin(kx)+\a%
_{-k}Y_t^{-k}
\cos(kx),
\end{eqnarray}
with $\a_k=\a_{-k}=\sqrt{\frac{\s_k}{2(\l+f(\tau_k))}}$ and
$Y_t^k$ as in %
\eqref{eqnOU}. Suppose $ ( a_k )$ to be eventually
nonincreasing and $ (
a_k ) \preceq(k^{-2\alpha})$ for some $\alpha\in(\frac{1}{2},1]$.
Then analogous results to Proposition~\ref{propspatiallift0} may be
established under various assumptions on $\s_k$ and $f(\tau_k)$, by
means of
the stability results given in Section~\ref{secrfs}.

%
\begin{example}
We consider the stochastic fractional heat equation with (possibly) colored
noise on the $1$-dimensional torus, that is,
%
\begin{equation}
\label{eqnfracheat} \qquad d\Psi^i_t = -\bigl((-\D)^\a
\Psi^i_t+\l\bigr) \,dt + d(-\D)^{-{\g/2}}
W^i_t,\qquad i=1,\ldots,d,
\end{equation}
where $\a\in(0,1]$, $\gamma\ge0$, $\lambda> 0$ and $W_t$ is a
cylindrical Wiener process. Hence, $f(\tau_k) = \llvert k\rrvert
^{2\alpha}$ and
$\s_k =
\llvert k\rrvert ^{-2\g}$. By elementary calculations we see $ ( \frac
{\sigma
_k}{\l%
+f(\tau_k)} ) \preceq(k^{-(2\gamma+ 2\alpha)})$ and thus the
conclusions of Proposition~\ref{propspatiallift0} hold if $2\g+2\a>
\frac{3}{2}$.
\end{example}

\subsubsection{Neumann boundary conditions}\label{sec432}

In the case of homogeneous Neumann boundary conditions, an orthogonal basis
of eigenvectors of $-\D$ on $L^2([0,2\pi])$ is given by
\begin{eqnarray*}
e_{k}(x)=\cos\biggl(\frac{k}{2}x \biggr),\qquad
\tau_{k}= \biggl(\frac
{k}{2}%
\biggr)^{2},\qquad
k\in\N\cup\{0\}.
\end{eqnarray*}
In order to be able to consider stationary Ornstein--Uhlenbeck
processes, we
need to shift the spectrum; that is, we consider
\[
d\Psi_{t}=- \bigl((-\Delta)^{\alpha}+1 \bigr)
\Psi_{t} \,dt+dW_{t}.
\]
We may then proceed as for Dirichlet boundary conditions, resolving
additional difficulties due to the shift of the spectrum as in the
proof of
Proposition~\ref{propspatiallift0}.

\section{The continuous case}\label{sec5}

In some cases, the covariance function of a Gaussian process $X$ is
given as the cosine transform of some function $f$. For example, this
is the case if the spectral measure of a stationary process has a
density $f$ with respect to the Lebesgue measure; cf. Example~\ref
{exspectralmeasure} and \cite{MR06}, Chapter~5.6. In this case, we
may obtain similar results as for random Fourier series. The key is a
continuous version of Lemma~\ref{lemmarhobound} which we are now
going to present. For a (symmetric) function $f \in L^1(\R)$, let
$\hat f$ denote its (real) Fourier transform. Then the following holds:
%
\begin{lemma}
\label{lemmarhoboundcontcase} Let $f \colon\R\to\R$ be
symmetric and in $L^1(\R) \cap C^2(\R\setminus\{0\})$. Assume $\hat
f \in L^1(\R)$ and
\begin{eqnarray*}
\lim_{\xi\to\infty} \bigl\llvert\xi^3 f^{\prime\prime}(
\xi)\bigr\rrvert+ \bigl\llvert\xi^2 f'(\xi)\bigr
\rrvert+ \bigl\llvert\xi f(\xi)\bigr\rrvert= 0
\end{eqnarray*}
and that there is an $x_0 \in(0,\infty]$ such that
\begin{eqnarray*}
\limsup_{R \to\infty} \int_{0}^R
\frac{\partial^2}{\partial\xi
^2}\bigl(f(\xi) \xi^2\bigr) F_{\xi}(x) \,d\xi
\leq0,
\end{eqnarray*}
for all $x\in(0,x_0)$ where $F_{\xi}(x) = \frac{1-\cos(\xi x)}{x^2}$
denotes the F\'{e}jer kernel. Then $\hat f$ is a convex function on $[0,x_0)$.
\end{lemma}

\begin{pf}
Since the proof is very similar to Lemma~\ref{lemmarhobound} we just
sketch it briefly. By F\'{e}jer's theorem for Fourier transforms (cf.
\cite{K89}, Theorem 49.3),
\begin{eqnarray*}
\lim_{R\to\infty} \frac{1}{2 \pi} \int_{-R}^R
\biggl(1 - \frac
{\llvert \xi\rrvert }{R} \biggr) \hat{g}(\xi) e^{ix\xi} \,d\xi= g(x),
\end{eqnarray*}
for all $x$ provided $g\in C(\R)\cap L^1(\R)$. Setting $g = \hat
{f}$, we obtain from
Fourier inversion
\begin{eqnarray*}
\sigma_R(x):= \int_{-R}^R
\biggl(1 - \frac{\llvert \xi\rrvert }{R} \biggr) f(\xi) e^{ix\xi} \,d\xi
\to\hat f(x)
\qquad\mbox{for } R \to\infty.
\end{eqnarray*}
Applying integration by parts twice, our assumptions imply that
\begin{eqnarray*}
\liminf_{R \to\infty} \sigma''_R(x)
\geq0
\end{eqnarray*}
for all $x\in(0,x_0)$. This implies the claim.
\end{pf}

\begin{remark}
Note that for a given $f\in L^{1}(\R)$, it does not follow that also
$\hat{f}%
\in L^{1}(\R)$. However, Bernstein's theorem states that the Fourier transform
of functions $f$ in the Sobolev space $H^{s}(\R)$ is contained in
$L^{1}(\R)$ for
all $s>\frac{1}{2}$; cf.~\cite{H83}, Corollary 7.9.4.
\end{remark}

%
\begin{example}\label{exfracOUcalc}
Consider the covariance $R$ of a fractional Ornstein--Uhlen\-beck process
with Hurst parameter $H\in(0,1/2]$; cf. Example~\ref
{exspectralmeasure}. Then $R(s,t) = K(t-s)$ with
\begin{eqnarray*}
K(x) = \int f(\xi) \cos(\xi x) \,d\xi,\qquad f(\xi) = c_H
\frac
{\llvert \xi\rrvert ^{1-2H}}{\lambda^2 + \xi^2},\qquad \lambda> 0.
\end{eqnarray*}
We prove that there is an $x_0 > 0$ such that $K$ is convex on
$[0,x_0)$. Since $f(\xi) = O(\xi^{-1-2H})$, $f \in H^s$ for any $s <
2H + 1/2$ and Bernstein's theorem\vspace*{1pt} implies that $\hat{f} \in L^1$ for
any $H>0$. An easy calculation shows that $g:= \partial^2_{\xi,\xi
}(f(\xi) \xi^2)$ is\vspace*{1pt} nonpositive on $[\xi_0,\infty)$ for some $\xi
_0 > 0$ and that $g(\xi) = O(-\xi^{-1-2H})$. It follows that
\begin{eqnarray*}
\limsup_{R \to\infty} \int_{0}^R
\frac{\partial^2}{\partial\xi
^2}\bigl(f(\xi) \xi^2\bigr) F_{\xi}(x) \,d\xi=
\int_{0}^{\infty} \frac{\partial
^2}{\partial\xi^2}\bigl(f(\xi)
\xi^2\bigr) F_{\xi}(x) \,d\xi.
\end{eqnarray*}
Note first that $\int_0^{\xi_0} g(\xi) F_{\xi}(x) \,d\xi$ is
uniformly bounded for $x\searrow0$. Furthermore $\lim_{x \to0}
F_{\xi}(x) = \xi^2/2$, and Fatou's lemma gives
\begin{eqnarray*}
\liminf_{x \to0} \int_{\xi_0}^{\infty} -
g(\xi) F_{\xi}(x) \,d\xi\geq- \frac{1}{2} \int
_{\xi_0}^{\infty} \xi^2 g(\xi) \,d\xi= + \infty.
\end{eqnarray*}
Hence
\begin{eqnarray*}
\lim_{x \to0} \int_0^{\infty} g(
\xi) F_{\xi}(x) \,d\xi= - \infty.
\end{eqnarray*}
Thus there is an $x_0$ such that $\int_0^{\infty} g(\xi) F_{\xi}(x)
\,d\xi\leq0$ for all $x \in(0,x_0)$, and we can apply Lemma~\ref
{lemmarhoboundcontcase} to conclude that $K$ is convex on $[0,x_0)$.
\end{example}

%
\begin{example}
Consider the SPDE
\[
d\Psi_t = -\bigl((-\Delta)^{\alpha} + \l\bigr)
\Psi_t \,dt + dW_t\qquad\mbox{on } \R,
\]
for some $\alpha\in(0,1]$, $\l> 0$. The stationary solution can be
written down explicitly (cf. \cite{W86}), namely
\[
\Psi_t(x) = \int_{-\infty}^t \int
_{\R} K_{t-s}(x,y) W(ds,dy),
\]
where $K$ is the fractional heat kernel operator associated to $%
-((-\Delta)^{\alpha} + \l)$ with Fourier transform given by
\[
\hat{K}_t(\xi) = e^{-t\llvert \xi\rrvert ^{2\alpha} - \l t}.
\]
After some calculations, one sees that the covariance $R$ of the spatial
process $x \mapsto\Psi_t(x)$ for every time point $t$ is given by $%
R(x,y)=K(x-y)$ where
\[
K(x) = \int f(\xi) \cos(\xi x) \,d\xi, \qquad f(\xi) = \frac
{1}{2\llvert \xi\rrvert ^{2\alpha} + 2\l}.
\]
With a similar calculation as in Example~\ref{exfracOUcalc}, one can
see that $K$ is convex in a neighborhood\vspace*{1pt} around $0$. It is easy to see
that $\sigma^2(x) = 2(K(x) - K(0)) = O(\llvert x\rrvert ^{2\alpha-
1})$ (using, e.g., \cite{MR06}, Theorem 7.3.1).
Hence we are in the situation of Example~\ref{StInII}, and we may
conclude that
\[
V_{1,\rho}\bigl(R_X;[x,y]\bigr) = O\bigl(\llvert y-x\rrvert
^{2\alpha- 1}\bigr)
\]
for $\llvert y-x\rrvert $ small enough. Applying \cite{FV10}, Theorem~35, we see that
$\Psi_t$ can be lifted, for every fixed time point $t$, to a process
$\bolds{\Psi}_t$ with
sample paths in $C_0^{0,\beta\mbox{-H\"{o}l}}([x,y],G^{[1/\beta]}(\R
^d)) $,
every $\beta< \alpha- 1/2$, provided $\alpha> 3/4$ and $\llvert y -
x\rrvert $ is
small enough. By concatenation one has the existence of spatial rough
paths lifts on all compact intervals in $\R$.
\end{example}

\section{Application to non-Markovian H\"{o}rmander theory}\label{sec6}

Consider a (rough) differential equation
%
\begin{equation}
\label{eqnRDE} dY=V ( Y ) \,d\mathbf{X}
\end{equation}
driven by a (Gaussian) rough path $\mathbf{X}$ along a vector field $%
V= ( V_{1},\ldots,V_{d} ) $, started at $Y_{0}=y_{0}\in\R^{e}$.
Assume $V$ to be bounded with bounded derivatives of all orders such
that H%
\"{o}rmander's condition $\operatorname{Lie} ( V_{1},\ldots,V_{d} )
\mid_{y_{0}}=\R^{e}$ holds.\footnote{%
We may also include a drift vector $V_{0}$, in which case we mean the
weak H%
\"{o}rmander condition.} If $X$ is sufficiently nondegenerate (e.g., fBm)
one can hope for a density of $Y_{t}$ at positive times $t>0$. This has been
achieved in a series of papers starting with Baudoin and Hairer \cite{BH07}
(with $\mathbf{X}$ fBm for $H>1/2$), followed by \cite{CF10,HP11,CHLT12}
which dealt, respectively, with general Gaussian signals $\mathbf{X}$
($\rho
<2$ subject to CYR\footnote{%
Complementary Young regularity for Cameron--Martin paths $h$: that is, $h$
has finite $q$-variation and $X ( \omega) $ has finite $p$%
-variation a.s. with $1/p+1/q>1$.}), fBm for $H>1/4$ and then again general
Gaussian signals ($\rho<2$ subject to CYR), now with a smoothness result.
The general case \cite{CF10,CHLT12} requires a number of assumptions
on $%
\mathbf{X}$ that are not always easy to check. To wit, even if $X$ is
fBm-like, in the sense that
\[
\sigma^{2} ( s,t ) =F ( t-s )\ge0
\]
with $F$ being concave and $F ( t ) =O ( t^{2H}
) $,
already the CYR is unclear in the aforementioned references \cite
{CF10,CHLT12}. Indeed, CYR for fBm (in case $H>1/4$) relies on the variation
embedding theorem \cite{FV06} which is not applicable in this more general
situation. Our results provide a convenient way to check the
assumptions of
\cite{CHLT12}. Let us illustrate how to proceed by the concrete
example of
an RDE driven by a (Gaussian) process (with i.i.d. components) with
stationary increments.

%
\begin{theorem}
\label{teoCHLT} Assume $F ( t ) =O ( t^{1/\rho
} ) $ with
$\rho\in{}[1,2)$ as $t\downarrow0$, $F$ concave and nonzero. Then
\[
F_{-}' ( T ) >0\qquad\mbox{for some }T>0,
\]
and $Y_{t}$ in \eqref{eqnRDE} has a smooth density for all $t\in(0,T]$.

This applies in particular when $X$ is given as a random Fourier series
as in Example~\ref{exrfsstationary} (with $\rho<2$) or as a fractional
Ornstein--Uhlenbeck process with Hurst parameter $H\in(1/4,1/2]$ as in
Example~\ref{exspectralmeasure}.
\end{theorem}

\begin{pf}
By assumption, $F$ is not identically equal to zero. In order to see
that $%
F_{-}' ( T ) >0$ for some $T$ small enough, assume
the opposite,
that is, $F_{-}' ( t ) \leq0$ for all $t>0$. Then
\[
\frac{F ( t+h ) -F ( t ) }{h}\leq F_{-}^{\prime
} ( t ) \leq0\qquad\forall
h,t>0,
\]
and thus $F$ is nonincreasing. Since $F ( 0 ) =0$ and
$F\geq0$,
this implies that $F$ is trivial and gives the desired contradiction.
We now
proceed by checking the conditions from \cite{CHLT12}. Condition 1 (CYR;
\cite{CF10,CHLT12}) follows from Theorem~\ref{teomain}, applied as in
Example~\ref{StInII} which yields mixed $ ( 1,\rho) $-variation
and thus (cf. part 1) complementary Young regularity. For Condition 2
from \cite{CHLT12}\footnote{%
For the reader's convenience we recall (the essence of) Condition 2 in
\cite{CHLT12}: there exists $c,\alpha>0$ such that $\operatorname
{Var}(X_{s,t}\mid F_{0,s}\vee F_{t,T})\geq c ( t-s ) ^{\alpha}$
for all $0 \leq s < t \leq T$.} we
first note that, leaving the imminent estimate (\ref{VarGE}) to the
end of the proof,
%
\begin{eqnarray}\label{VarGE}
&& 2\operatorname{Var}(X_{s,t}\mid\mcF_{0,s}\vee
\mcF_{t,T})\nonumber
\\
&&\qquad \geq2R\pmatrix{ s,t
\cr
0,T} \nonumber
\\
&&\qquad = \sigma^{2}(0,t)-\sigma^{2}(0,s)+\sigma^{2}(s,T)-
\sigma^{2}(t,T)
\\
&&\qquad = F(t)-F(s)+F(T-s)+F(T-t)\nonumber
\\
&&\qquad \geq 2F_{-}' ( T ) (t-s),
\nonumber
\end{eqnarray}
where we used (thanks to concavity) that the left-hand side derivative
of $F$ at $%
T$ exists and
\[
\inf_{0\leq s<t\leq T}\frac{F ( t ) -F ( s )
}{t-s}%
=F_{-}'
( T ).
\]
By assumption, $F_{-}' ( T ) >0$, and so Condition
2 holds
with $\alpha=1$. Also note that $F_{-}' ( T ) $ is
nonincreasing in $T$; thus Condition 2 remains valid upon decreasing $T$.

Next, we prove that (\cite{CHLT12}, Condition 4, page 10) is satisfied.
Due to concavity of~$F$
and Lemma~\ref{lemincrcorr}, $X$ has nonpositively correlated
increments, and it is enough to show that (cf. \cite{CHLT12}, page 11)
%
\begin{equation}
\label{EX0SXst} \E X_{0,S}X_{s,t}=R\pmatrix{ 0,S
\cr
s,t}
\geq0\qquad\forall [ s,t]\subseteq{}[ 0,S]\subseteq{}[0,T],
\end{equation}
which is clear from our condition (B.ii) of Theorem~\ref{teomain}, which
was seen to be verified in the present situation in Example~\ref
{StInII}. We
also note that (\cite{CHLT12}, Condition~3, page 10) is implied by
Condition 4
(cf. \cite{CHLT12}, Corollary 6.8). In conclusion, \cite{CHLT12}, Theorem 3.5, implies the claim.

It remains to prove estimate (\ref{VarGE}). To this end, let $\mcG:=\mcF
_{0,s}\vee\mcF_{t,T}$. Since $X$ is Gaussian, $\operatorname
{Var}%
(X_{s,t}\mid\mcG)$ is deterministic, and by a simple argument (detailed
in \cite{CHLT12}, Lemma 4.1) one has
%
\begin{eqnarray}
\operatorname{Var}(X_{s,t}\mid\mcG) = \inf_{Y\in L^{2}(\O,\mcG,\P)}
\llVert X_{s,t}-Y\rrVert_{2}^{2}, \label{VarXEINF}
\end{eqnarray}
where the $\inf$ is achieved at $Y= \E[ X_{s,t} \mid\mcG]$, an element
in the first Wiener--It\^o chaos, that is, the $L^2$-closure of $\{
X_t\dvtx 0 \leq t \leq T \}$ and of course $\mcG$-measurable. As a
consequence, $\E[ X_{s,t} \mid\mcG] = \lim_n Y_n$ in $L^2$ for
suitable ``simple'' approximations of the form, with
$[t_{i}^{n},t_{i+1}^{n}] \subset[0,s] \cup[t,T]$,
\[
Y_{n}=\sum_{i=1}^{k_{n}}a_{i}^{n}X_{t_{i}^{n},t_{i+1}^{n}},
\]
and we can replace the $\inf$ in (\ref{VarXEINF}) by the $\inf$
taken over such simple elements.
In what follows let us write $(\tilde t^n_i)$ for the\vspace*{1pt} dissection
obtained from $(t^n_i\dvtx1\le i < k_n) \cup\{ s, t \}$. This way,
we may condense the expression $X_{s,t}-\sum
_{i}a_{i}^{n}X_{t_{i}^{n},t_{i+1}^{n}}$ to $\sum_{i}\td
a_{i}^{n}X_{\td t_{i}^{n},\td t_{i+1}^{n}}$.
Using elementary estimates such as $\alpha_i \alpha_j \le(\alpha
_i^2 + \alpha_j^2)/2$ and symmetry of $R$, we find
\begin{eqnarray*}
\llVert X_{s,t}-Y_{n}\rrVert_{2}^{2}& =&\E\biggl\llvert X_{s,t}-%
\sum_{i}a_{i}^{n}X_{t_{i}^{n},t_{i+1}^{n}}
\biggr\rrvert^{2}
 =\E\biggl\llvert\sum_{i}\td
a_{i}^{n}X_{\td t_{i}^{n},\td t_{i+1}^{n}}\biggr\rrvert^{2}
\\
& =&\sum_{i,j}\td a_{i}^{n}
\td a_{j}^{n}\E X_{\td t_{i}^{n},\td%
t_{i+1}^{n}}X_{\td t_{j}^{n},\td t_{j+1}^{n}}
 =\sum_{i,j}\td a_{i}^{n}
\td a_{j}^{n}R\pmatrix{ \td t_{i}^{n},
\td t_{i+1}^{n}
\cr
\td t_{j}^{n},
\td t_{j+1}^{n}}
\\
& \geq& -\sum_{i}\sum_{j\neq i}
\bigl\llvert\td a_{i}^{n}\bigr\rrvert\bigl\llvert\td
a_{j}^{n}\bigr\rrvert \biggl\llvert R\pmatrix{ \td
t_{i}^{n},\td t_{i+1}^{n}
\cr
\td
t_{j}^{n},\td t_{j+1}^{n}} \biggr
\rrvert+\sum_{i}\bigl(\td a_{i}^{n}
\bigr)^{2}R \pmatrix{ \td t_{i}^{n},\td
t_{i+1}^{n}
\cr
\td t_{i}^{n},\td
t_{i+1}^{n}}
\\
& \geq& -\sum_{i}\sum_{j\neq i}
\bigl(\td a_{i}^{n}\bigr)^{2}\biggl\llvert R
\pmatrix{ \td t_{i}^{n},\td t_{i+1}^{n}
\cr
\td t_{j}^{n},\td t_{j+1}^{n}}
\biggr\rrvert+\sum_{i}\bigl(\td a_{i}^{n}
\bigr)^{2}R\pmatrix{ \td t_{i}^{n},\td
t_{i+1}^{n}
\cr
\td t_{i}^{n},\td
t_{i+1}^{n}}.
\end{eqnarray*}
Due to nonpositively correlated increments we may drop the minus and
absolute values in the line above and %
combine both sums. Thanks to (\ref{EX0SXst}), we can then finish the
desired estimate,
\begin{eqnarray*}
\llVert X_{s,t}-Y_{n}\rrVert_{2}^{2}&
\ge&\sum_{i}\bigl(\td a_{i}^{n}
\bigr)^{2}\sum_{j}R%
\pmatrix{
\td t_{i}^{n},\td t_{i+1}^{n}
\cr
\td t_{j}^{n},\td t_{j+1}^{n}}
 =\sum_{i}\bigl(\td a_{i}^{n}
\bigr)^{2}R\pmatrix{ \td t_{i}^{n},\td
t_{i+1}^{n}
\cr
0,T}
\\
& =&R\pmatrix{ s,t
\cr
0,T} +\sum_{i}
\bigl(a_{i}^{n}\bigr)^{2}R \pmatrix{
t_{i}^{n},t_{i+1}^{n}
\cr
0,T}
 \geq R\pmatrix{ s,t
\cr
0,T}.
\end{eqnarray*}\upqed
\end{pf}

\section*{Acknowledgments}
Peter K. Friz,
Benjamin Gess  and
Sebastian Riedel
would like to thank the Mathematisches Forschungsinstitut
Oberwolfach where parts of this work was first presented in August 2012 as
part of the workshop ``Rough paths and PDEs.'' They would also like to
thank F. Russo for the references \mbox{\cite{KRT07,KR10}} which led to
Example \ref{examplebifbm}. Peter K. Friz would like to thank S. Tindel for pointing out
similar aspects of the present conditions with those from~\cite{CHLT12};
this led to Section~\ref{sec6}. At last, all authors would like to thank the
referee for detailed comments.



%

\printaddresses
\end{document}